\documentclass[a4paper,12pt]{article}
\usepackage{amssymb,amscd,amsfonts,amsbsy}
\usepackage{enumerate}
\usepackage{epsfig}
\input amssymb.sty
\input euscript.sty

\newcommand{\RR}{{\mathbb{R}}}
\newcommand{\NN}{{\mathbb{N}}}
\newcommand{\ZZ}{{\mathbb{Z}}}
\newcommand{\CC}{{\mathbb{C}}}
\newcommand{\eps}{\varepsilon}
\newcommand{\bp}{\noindent {\it Proof}.\,\,\,}
\newcommand{\ep}{\hfill$\Box$ \vskip 0.08in}
\newcommand{\dint}{\int\!\!\!\int}

\newtheorem{proposition}{Proposition}[section]
\newtheorem{theorem}[proposition]{Theorem}
\newtheorem{lemma}[proposition]{Lemma}
\newtheorem{corollary}[proposition]{Corollary}

\begin{document}

\title{Hardy and BMO spaces associated to divergence form
elliptic operators}

\author{Steve Hofmann (hofmann@math.missouri.edu)\\ Svitlana Mayboroda (svitlana@math.ohio-state.edu)\thanks{
2000 {\it Math Subject Classification:} 42B30, 42B35, 42B25, 35J15.
}}

\date{ }

\maketitle

\tableofcontents

\newpage

\begin{abstract}
Consider a second order divergence form elliptic operator $L$
with complex bounded coefficients. In general, operators
related to it (such as the Riesz transform or square function) lie
beyond the scope of the Calder\'{o}n-Zygmund theory. They need not
be bounded in the classical Hardy, BMO and even some $L^p$ spaces.

In this work we develop a
theory of Hardy and BMO spaces associated to $L$, which includes,
in particular, molecular decomposition, maximal and square function
characterizations, duality of Hardy and BMO spaces,  and John-Nirenberg
inequality.
\end{abstract}

\section{Introduction and statement of main results}
\setcounter{equation}{0}

Extensive study of classical real-variable Hardy spaces in $\RR^n$ began in the
early 60's with the fundamental paper of E.M.\,Stein and G.\,Weiss
\cite{StWe}. Since then these classes of functions have played an
important role in harmonic analysis, naturally continuing the scale
of $L^p$ spaces to the range of $p\leq 1$. Although many
real-variable methods have been developed (see especially the work of Fefferman and Stein \cite{FeSt}), the
theory of Hardy spaces is intimately connected with the properties
of harmonic functions and with the properties of Laplacian.

For instance, Hardy space $H^1(\RR^n)$ can be viewed as a
collection of functions $f\in L^1(\RR^n)$ such that the Riesz
transform $\nabla \Delta^{-1/2}f$ belongs to $L^1(\RR^n)$.
One also has alternative characterizations of $H^1(\RR^n)$ by the
square function and the non-tangential maximal function associated
to the Poisson semigroup generated by Laplacian. To be
precise, fix a family of non-tangential cones
$\Gamma(x):=\{(y,t)\in\RR^n\times(0, \infty):\,|x-y|<t\}$,
$x\in\RR^n$, and define

\begin{eqnarray}\label{eq1.1}
S^\Delta f(x)&=&\left(\dint_{\Gamma(x)}\left|t\nabla e^{-t\sqrt
\Delta} f(y)\right|^2\,\frac
{dydt}{t^{n+1}}\right)^{1/2}, \\[4pt]\label{eq1.2}
{\cal N}^\Delta f(x)&=&\sup_{(y,t)\in\Gamma(x)}\left|e^{-t\sqrt
\Delta}f(y)\right|.
\end{eqnarray}

\noindent Then $\|{\cal N}^\Delta f\|_{L^1(\RR^n)}$ and $\|S^\Delta
f\|_{L^1(\RR^n)}$ give equivalent norms in the  space $H^1(\RR^n)$,
that is

\begin{equation}\label{eq1.3}
\|{\cal N}^\Delta f\|_{L^1(\RR^n)}\approx \|S^\Delta
f\|_{L^1(\RR^n)}\approx \|f\|_{H^1(\RR^n)}.
\end{equation}

Consider now a general elliptic operator in divergence form with
complex bounded coefficients. Let $A$ be an $n\times n$ matrix with
entries

\begin{equation}\label{eq1.4}
a_{jk}:L^{\infty}(\RR^n)\longrightarrow \CC,\quad j=1,...,n, \quad
k=1,...,n,
\end{equation}

\noindent satisfying the ellipticity condition

\begin{equation}\label{eq1.5}
\lambda |\xi|^2\leq \Re e A\xi\cdot\bar \xi \quad \mbox{and} \quad
|A\xi\cdot\bar \zeta|\leq \Lambda |\xi||\zeta|, \quad
\forall\,\xi,\zeta\in \CC^n,
\end{equation}

\noindent for some constants $0<\lambda\leq \Lambda <\infty$. Then
the second order divergence form operator is given by

\begin{equation}\label{eq1.6}
Lf:=-{\rm div}(A\nabla f),
\end{equation}

\noindent interpreted in the weak sense via a sesquilinear form.

Unfortunately the classical Hardy spaces turn out to be fairly
useless when studying the problems connected with the general
operator $L$ defined in (\ref{eq1.4})--(\ref{eq1.6}). For example,
the corresponding Riesz transform $\nabla L^{-1/2}$ need not be
bounded from $H^1(\RR^n)$ to $L^1(\RR^n)$. Indeed, we know from
\cite{KatoMain} that the Riesz transform is bounded in $L^2(\RR^n)$.
In fact, we also know that there exists $p_L, 1\leq p_L<2n/(n+2),$
such that the Riesz transform is bounded in $L^p(\RR^n)$ for $p_L
<p\leq 2$ (see \cite{HoMa}, \cite{BK}, \cite{A}, \cite{AT1},
\cite{AT} and \S{2} for more details; see also \cite{BK2} for
related theory). If, in addition, it mapped $H^1(\RR^n)$ to
$L^1(\RR^n)$, then by interpolation $\nabla
L^{-1/2}:L^p(\RR^n)\rightarrow L^p(\RR^n)$ for all $1<p\leq 2$,
which, in general, is not true (i.e., the best possible  $p_L$ can
be strictly greater than $1$). Similar statements apply to the heat
semigroup $e^{-tL}$ (for which $L^p$ boundedness may fail also for
$p$ finite, but sufficiently large), as well as to the square
function (see (\ref{eq1.12}) below). This failure of $L^p$ bounds
for some $p \in (1,2)$ (and in $(2,\infty)$), which relies on
counterexamples built in \cite{MNP}, \cite{ACT}, \cite{Davies},
along with some observations in \cite{AuscherSurvey}, illustrates
the significant difference between the operators associated to $L$
and those arising in the classical Calder\'{o}n-Zygmund theory.  We
note that it has been shown in \cite{AuscherSurvey} that the
intervals of $p\leq 2$ such that the heat semigroup and Riesz
transform are $L^p$-bounded have the same interior.  In the sequel,
we shall denote by $(p_L,\widetilde{p}_L)$ the interior of the
interval of $L^p$ boundedness of the semigroup, and we recall
\cite{AuscherSurvey} that $\widetilde{p}_L > 2n/(n-2)$.

One of the two goals of this paper is to generalize the classical
theory of Hardy spaces in order to ameliorate the deficiencies described in the
previous paragraph:  given an  elliptic operator $L$ as above, we construct an
$H^1$ space adapted to $L$, which, for example, is mapped into $L^1$ by
the Riesz transforms $\nabla L^{-1/2}$.
A primary virtue of the Hardy space $H^1$ is its many useful characterizations, thus, we aim to provide
analogues to most of them, including
the atomic (or molecular) decomposition and characterizations by
square and non-tangential maximal functions.  We note that one problem which remains open
is that of finding a Riesz transform characterization of the adapted $H^1$ space;
i.e., we do not yet know whether $L^1$ bounds for $\nabla L^{-1/2} f$, or some suitable substitute, imply that $f$ belongs to our Hardy space (as mentioned above, we do prove the converse). On the other hand, there is an analogous $H^p$ theory, $p \neq 1$, which we plan to present in
a forthcoming joint paper with A. McIntosh, and in which, in contrast to the classical setting, even the case $p>1$ may involve Hardy spaces which are strictly smaller than $L^p$.  For at least some such spaces,
we do have a Riesz transform characterization.  Finally, we remark that our Hardy spaces belong to an interpolation scale which includes $L^p, \, p_L < p < \widetilde{p}_L.$

Let us begin with the matter of atomic/molecular decompositions. The decomposition into simple
building blocks, atoms, originally proved by R.\,Coifman  for
$n=1$ (\cite{Co}) and by R.\,Latter in higher dimensions
(\cite{La}), is a key feature of $H^1(\RR^n)$. In this paper we
take as our starting point an analogous molecular decomposition, in the spirit of
the one introduced in the classical setting by Taibleson and Weiss \cite{TW}.

Throughout the paper for cube $Q\subset\RR^n$ we denote by $l(Q)$
the sidelength of $Q$ and set

\begin{equation}\label{eq1.7}
S_0(Q):=Q, \quad Q_i=2^iQ,\quad \mbox{ and }\quad
S_i(Q):=2^iQ\setminus 2^{i-1}Q \mbox{ for } i=1,2,...,
\end{equation}

\noindent where $2^iQ$ is cube with the same center as $Q$ and
sidelength $2^il(Q)$.

Let $1 \leq p_L<2n/(n+2)$ and $\widetilde{p}_L > 2n/(n-2)$ retain the same significance as above
(that is, they are the endpoints of the interval of $L^p$ boundedness for the heat semigroup).
A function $m\in L^p(\RR^n)$,
$p_L<p<\widetilde{p}_L$, is called a $(p,\eps,M)$-{\it molecule},
$\eps>0$ and $M\in\NN, M > n/4$, if there exists a cube $Q\subset\RR^n$ such
that
\begin{eqnarray}
&&\hskip -0.7cm (i)\,\,\, \|m\|_{L^p(S_i(Q))}\leq \,2^{-i(n-n/p
+\eps)}|Q|^{1/p-1},
\,\, i=0,1,2,..., \label{eq1.8}\\[4pt]
&&\hskip -0.7cm (ii)\, \|(l(Q)^{-2}L^{-1})^k m\|_{L^p(S_i(Q))}\leq
\,2^{-i(n-n/p +\eps)}|Q|^{1/p-1}, \,\, i=0,1,2,..., \,\,
k=1,...,M.\label{eq1.9}
\end{eqnarray}

\noindent Having fixed some $p$, $\eps$ and $M$, we will often use
the term molecule rather than $(p,\eps,M)$-molecule in the sequel.
Then the Hardy space can be defined as

\begin{equation}\label{eq1.10}
H^1_L(\RR^n)\equiv \left\{\sum_{i=0}^\infty
\lambda_jm_j:\,\{\lambda_j\}_{j=0}^\infty\in \ell^1 \mbox{ and
}\,m_j \mbox{ are  molecules}\right\},
\end{equation}

\noindent with the norm given by

\begin{equation}\label{eq1.11}
\|f\|_{H^1_L(\RR^n)}=
\inf\Bigl\{\sum_{j=0}^\infty|\lambda_j|:\,\,f=\sum_{j=0}^\infty\lambda_j
m_j,\,\, \{\lambda_j\}_{j=0}^\infty \in\ell^1 \mbox{ and
}m_j\textrm{ are molecules} \Bigr\}.
\end{equation}
We shall prove in Section 4 that any fixed choice of $p, \eps$ and $M$ within the
allowable parameters stated above, will generate the same space.

We remark that our molecules are, in particular, classical $H^1$ molecules, so that
$H^1_L \subseteq H^1$ (and this containment is proper for some $L$).  We also
observe that one may construct examples of $H^1_L$ molecules by hand:  given a cube
$Q$, let $f \in L^2(Q),$ with $\|f\|_2 \leq |Q|^{-1/2}$, and set
\begin{equation}\label{eq1.mol} m= (\ell (Q)^2 L)^M e^{-\ell (Q)^2
L} f, \,\,\,\, \tilde{m} = (I - (I + \ell(Q)^2 L)^{-1})^M f.\end{equation}
One may then easily check, using the ``Gaffney estimate" (\ref{eq2.8}) below,
that up to a suitable normalizing constant, $m, \tilde{m}$ are $(2,\eps,M)$ molecules for all
$\eps > 0.$

Molecules have appeared in the $H^1$ theory as an analogue of
atoms not having compact support but rapidly decaying away from
some cube $Q$ (see, e.g., \cite{GaCu}). However, the classical
vanishing moment condition ($\int_{\RR^n} m(x)\,dx=0$) does not
interact well with the operators we have in mind because it does not
provide appropriate cancellation. Instead, we have imposed a new
requirement that the molecule ``absorbs" properly normalized
negative powers of the operator $L$ -- the condition made precise
in (\ref{eq1.9}). In such a setting it can be proved, for instance,
that the Riesz transform

\begin{equation}\label{eq1.17}
\nabla L^{-1/2}: H^1_L(\RR^n)\longrightarrow L^1(\RR^n),
\end{equation}

\noindent as desired.

Next, given an operator $L$ as above and function $f\in L^2(\RR^n)$,
consider the following quadratic and maximal operators associated to the heat
semigroup generated by $L$

\begin{eqnarray}\label{eq1.12}
S_hf(x)&:=&\left(\dint_{\Gamma(x)}|t^2Le^{-t^2L}f(y)|^2\,\frac
{dydt}{t^{n+1}}\right)^{1/2},\\[4pt]
{\cal N}_h f(x)&:=&
\sup_{(y,t)\in\Gamma(x)}\left(\frac{1}{t^n}\int_{B(y,t)}
|e^{-t^2L}f(z)|^2\,dz\right)^{1/2}, \label{eq1.13}
\end{eqnarray}

\noindent where $B(y,t)$, $y\in\RR^n$, $t\in(0,\infty)$, is a ball
in $\RR^n$ with center at $y$ and radius $t$ and $x\in\RR^n$.
These are the natural modifications of
(\ref{eq1.1})--(\ref{eq1.2}). We use an extra averaging in the
space variable for the non-tangential maximal function in order to
compensate for the lack of pointwise estimates on the heat semigroup (an
idea originating in \cite{KePi}).

Alternatively, one can consider the Poisson semigroup generated by
the operator $L$ and the operators

\begin{eqnarray}\label{eq1.14}
S_Pf(x)&:=&\left(\dint_{\Gamma(x)}|t\nabla e^{-t\sqrt
L}f(y)|^2\,\frac {dydt}{t^{n+1}}\right)^{1/2},
\\[4pt]\label{eq1.15}
{\cal N}_P f(x)& := &\sup_{(y,t)\in\Gamma(x)} \left( \frac{1}{t^n}
\int_{B(y,t)} |e^{-t\sqrt L} f(x) |^2 dx\right)^{1/2},
\end{eqnarray}

\noindent with $x\in\RR^n$,  $f\in L^2(\RR^n)$.

Define the spaces $H^1_{S_h}(\RR^n)$, $H^1_{{\cal N}_h}(\RR^n)$,
$H^1_{S_P}(\RR^n)$ and $H^1_{{\cal N}_P}(\RR^n)$ as the completion
of $L^2(\RR^n) \cap L^1(\RR^n)$ in the norms given by the $L^1$ norm of the
corresponding square or maximal function, plus the
$L^1$ norm of the function itself;  e.g.,
$$\|f\|_{H^1_{S_h}(\RR^n)} \equiv \|S_h f\|_{L^1(\RR^n)} + \|f\|_{L^1(\RR^n)}.$$
Then the following result
holds.

\begin{theorem}\label{t1.1} Suppose $p_L<p<\widetilde{p}_L$, $\eps>0$ and $M>n/4$ in
(\ref{eq1.8})--(\ref{eq1.9}). For an operator $L$ given by
(\ref{eq1.4})--(\ref{eq1.6}), the Hardy spaces $H^1_{L}(\RR^n)$,
$H^1_{S_h}(\RR^n)$, $H^1_{{\cal N}_h}(\RR^n)$, $H^1_{S_P}(\RR^n)$,
$H^1_{{\cal N}_P}(\RR^n)$ coincide. In other words,

\begin{equation}\label{eq1.16}
\|f\|_{H^1_L(\RR^n)}\approx \|f\|_{\widetilde{H}^1_L(\RR^n)}\approx\|S_hf\|_{L^1(\RR^n)} \approx \|{\cal
N}_hf\|_{L^1(\RR^n)}\approx \|S_Pf\|_{L^1(\RR^n)}\approx \|{\cal
N}_Pf\|_{L^1(\RR^n)},
\end{equation}

\noindent for every $f\in H^1_L(\RR^n)$.
\end{theorem}

The second half of the paper is devoted to the analogue of the
space of functions of bounded mean oscillation.

The $BMO$ space originally introduced by F.\,John and
L.\,Nirenberg in \cite{JN} in the context of partial differential
equations, has been identified as the dual of classical $H^1$ in the
work by C.\,Fefferman and E.\,Stein \cite{FeSt}. Analogously to
the Hardy space, $BMO(\RR^n)$ plays a significant role in analysis
often substituting for $L^\infty$ as an endpoint of the $L^p$ scale.
And for reasons similar to those outlined above, the classical
$BMO$ is not always applicable for work with the operator $L$ as
well as many other operators beyond the Calder{\'o}n-Zygmund
class.

The second goal of this paper is to generalize the classical
notion of $BMO$. We introduce the space suitable for the operator
$L$ and prove that it is equipped with the characteristic
properties of $BMO$, in particular, it is tied up with the Hardy
space theory via duality.

We define our adapted $BMO$ spaces as follows.  For $\eps >0$
and $M \in \mathbb{N}$ we introduce the norm
$$\|\mu\|_{{\bf M}^{2,\eps,M}_0} \equiv  \sup_{j\geq 0}
2^{j(n/2+\eps)}\sum_{k=0}^M\|L^{-k}\mu\|_{L^2(S_j(Q_0))},$$
where $Q_0$ is the unit cube centered at $0$,
and set
$${\bf M}^{2,\eps,M}_0 \equiv \{ \mu : \|\mu\|_{{\bf M}^{2,\eps,M}_0} < \infty\}.$$
Implicitly, of course, this space depends upon $L$, and we shall write
${\bf M}^{2,\eps,M}_0(L)$ when we need to indicate this dependence explicitly.
We note that if $ \mu \in {\bf M}^{2,\eps,M}_0$ with norm $1$, then
$\mu$ is a $(2,\eps,M)$ molecule adapted to $Q_0$.  Conversely, if $m$ is a
$(2,\eps,M)$ molecule adapted to any cube, then $m \in {\bf M}^{2,\eps,M}_0$
(this follows from the fact that, given any two cubes $Q_1$ and $Q_2$, there exists
integers $K_1$ and $K_2$, depending upon $\ell(Q_1), \ell(Q_2)$ and dist$(Q_1,Q_2)$, such that
$2^{K_1}Q_1 \supseteq Q_2$ and $2^{K_2}Q_2 \supseteq Q_1$).  Let $({\bf M}^{2,\eps,M}_0)^*$
be the dual of ${\bf M}^{2,\eps,M}_0$, and let $A_t$ denote either $(I + t^2L)^{-1}$ or $ e^{-t^2L}$.
We claim that if $f \in ({\bf M}^{2,\eps,M}_0)^*$, then $(I-A_t^*)^M f$ is globally
well defined in the sense of distributions, and belongs to $L^2_{loc}.$  Indeed, if $\varphi \in L^2(Q)$
for some cube $Q$, it follows from the Gaffney estimate (\ref{eq2.8}) below that $(I-A_t)^M \varphi
\in {\bf M}^{2,\eps,M}_0$ for every $\eps > 0$.  Thus,
\begin{equation}\label{eq1.loc}
\langle (I-A_t^*)^M f, \varphi \rangle \equiv  \langle  f, (I-A_t)^M \varphi \rangle
\leq C_{t,\ell (Q), \textrm{dist}(Q,0)} \|f\|_{({\bf M}^{2,\eps,M}_0)^*}\|\varphi\|_{L^2(Q)}.\end{equation}
Since $Q$ was arbitrary, the claim follows.  Similarly, $(t^2L^*)^M A_t^* f \in L^2_{loc}.$

We are now ready to define our adapted $BMO$ spaces.  We suppose now and in the sequel that $M>n/4$.
An element
\begin{equation}\label{eq1.def}f \in \cap_{\eps > 0}({\bf M}^{2,\eps,M}_0)^*\equiv ({\bf M}^{2,M}_0)^*\end{equation}
is said to belong to $BMO_{L^*}(\RR^n)$
if

\begin{equation}\label{eq1.18}
\|f\|_{BMO_{L^*}(\RR^n)}:=\sup_{Q\subset\RR^n}\left(\frac{1}{|Q|}\int_Q
\left|(I-e^{-l(Q)^2L^*})^Mf(x)\right|^2\,dx\right)^{1/2}<\infty,
\end{equation}

\noindent where $M>n/4$ and $Q$ stands for a cube in $\RR^n$.  Eventually, we shall see
that this definition is independent of the choice of $M > n/4$ (up to ``modding out"
elements in the null space of the operator $(L^*)^M$, as these are annihilated by
$(I-e^{-l(Q)^2L^*})^M$; we thank Lixin Yan for this observation). Clearly, we can define
$BMO_L$ by interchanging the roles of $L$ and $L^*$ in the preceding discussion.
Using the ``Gaffney" estimate (\ref{eq2.8}) below, and the fact that
$e^{-l(Q)^2L} 1 = 1$, one may readily verify that $BMO \subseteq BMO_L.$
Compared to the classical definition, in (\ref{eq1.18}) the heat
semigroup $e^{-l(Q)^2L}$ plays the role of
averaging over the cube, and an extra power $M>n/4$
provides the necessary cancellation.

As usual, we have some
flexibility in the choice of $BMO$ norm guaranteed by an
appropriate version of the John-Nirenberg inequality (see
\cite{JN} for the case of usual $BMO$), although in our case we obtain a more restricted range of equivalence.
To be more precise, let us denote by $BMO_L^p(\RR^n)$,
$1<p<\infty$, the set of elements of $({\bf M}^{2,M}_0)^*$ with the property that

\begin{equation}\label{eq1.21}
\|f\|_{BMO_L^p(\RR^n)}:=\sup_{Q\subset\RR^n}\left(\frac{1}{|Q|}\int_Q
\left|(I-e^{-l(Q)^2L})^Mf(x)\right|^p\,dx\right)^{1/p}
\end{equation}

\noindent is finite. Clearly, $BMO_L^2(\RR^n)=BMO_L(\RR^n)$. Then
the following result holds.

\begin{theorem}\label{t1.3} For all $p$ such that
$p_L<p<\widetilde{p}_L$ the spaces $BMO_L^p(\RR^n)$ coincide.
\end{theorem}

Another important feature of classical $BMO$ is its Carleson measure
characterization. Roughly speaking, we shall establish the following analogue of
the Fefferman-Stein criterion:

\begin{equation}\label{eq1.20}
f\in BMO_L(\RR^n)\quad\Longleftrightarrow\quad
\left|(t^2L)^Me^{-t^2L}f(y)\right|^2\,\frac{dydt}{t}\quad \mbox{
is a Carleson measure};
\end{equation}

\noindent see Theorem~\ref{t9.1} for the precise statement, in which, as in the classical case,
a certain ``controlled growth" hypothesis is needed to prove the
$\Leftarrow$ direction.

And finally, we prove the desired duality with the Hardy spaces.

\begin{theorem}\label{t1.2} Suppose $p_L<p<\widetilde{p}_L$, $\eps>0$ and $M>n/4$ in
(\ref{eq1.8})--(\ref{eq1.9}), (\ref{eq1.18}). Then for an operator
$L$ given by (\ref{eq1.4})--(\ref{eq1.6})

\begin{equation}\label{eq1.19}
\left(H^1_L(\RR^n)\right)^*=BMO_{L^*}(\RR^n).
\end{equation}
\end{theorem}

As we have mentioned, our results lie beyond the classical
Calder{\'o}n-Zygmund setting. Moreover, the methods we have at our
disposal are substantially restricted. For instance, no analogue
of the subaveraging property of harmonic functions, no maximum
principle, no regularity or pointwise bounds for the kernel of the
heat or Poisson semigroup are available. The operators we work
with do not even possess a kernel in the regular sense. In fact,
we employ only certain estimates in $L^2$ and $L^p$ with $p$ close
to 2, controlling the growth of the heat semigroup and the
resolvent.

We should mention that some of the problems discussed above were
earlier considered under different restrictions on the operator
and underlying space, but always assuming pointwise bounds for the
heat semigroup (see \cite{AR}, \cite{DL}, \cite{DL2}).

The layout of the paper is as follows: Section 2 contains a few
preliminary results, regarding general square functions and
properties of the operator $L$. Section 3 is devoted to behavior
of the sublinear operators, in particular, Riesz transform, acting
on the functions from $H^1_L(\RR^n)$. Sections 4, 5, 6 and 7 cover
the characterizations of Hardy space announced in
Theorem~\ref{t1.1}. Finally, in Sections 8, 9 and 10 we discuss
the spaces $BMO_L(\RR^n)$, duality, connection with Carleson
measures and John-Nirenberg inequality.

While this work was in preparation, we learned that some of the
questions concerning Hardy spaces have been considered
independently by P.\,Auscher, A.\,McIntosh and E.\,Russ.  Their
results are stated in the context of a Dirac operator on a
Riemannian manifold, but some arguments carry over to the present
setting as well. To the best of our knowledge, the theory of
$BMO_L(\RR^n)$ spaces is unique to this paper.

\vskip 0.08in

\noindent{\it Acknowledgements}. This work was conducted at the
University of Missouri at Columbia and the Centre for Mathematics
and its Applications at the Australian National University. The
second author thanks ANU for their warm hospitality.  We particularly thank
Pascal Auscher for showing us the proof of the $L^2$ boundedness of the non-tangential maximal
function of the Poisson semigroup, for pointing out a gap in the proof of our duality result in the original draft of this manuscript, and for providing us with the argument which we have formalized here as
Lemma \ref{l3.sub}.  The first author also thanks Emmanuel Russ for a helpful conversation, and especially for
explaining how to prove the completeness of molecular $H^1_L$.  Finally, we thank Chema Martell and
Lixin Yan for helpful comments.

\section{Notation and preliminaries}
\setcounter{equation}{0}

Let $\Gamma^{\alpha}$, $\alpha>0$, be the cone of aperture $\alpha$,
i.e. $\Gamma^{\alpha}(x):=\{(y,t)\in\RR^n\times(0,
\infty):\,|x-y|<\alpha t\}$ for $x\in\RR^n$. Then for a closed set
$F\in\RR^n$ we define a saw-tooth region ${\cal
R}^\alpha(F):=\bigcup_{x\in F}\Gamma^\alpha(x)$. For simplicity we
will often write $\Gamma$ in place of $\Gamma^1$ and ${\cal R}(F)$
instead of ${\cal R}^1(F)$.

Suppose $F$ is a closed set in $\RR^n$ and $\gamma\in(0,1)$ is
fixed. We set

\begin{equation}\label{eq2.1}
F^*:=\left\{x\in\RR^n:\, \mbox{for every} \,B(x), \mbox{ ball in
$\RR^n$ centered at $x$},\quad  \frac{|F\cap B(x)|}{|B(x)|}\geq
\gamma\right\},
\end{equation}

\noindent and every $x$ as above is called a point having global
$\gamma$-density with respect to $F$. One can see that $F^*$ is
closed and $F^*\subset F$. Also,

\begin{equation}\label{eq2.2}
\,^cF^*=\{x\in\RR^n:\,{\cal M}(\chi_{\,^cF})(x)>1-\gamma\},
\end{equation}

\noindent which implies $|\,^cF^*|\leq C|\,^cF|$ with $C$
depending on $\gamma$ and the dimension only.

  Here and throughout the paper we denote by $\,^cF$
the complement of $F$, $\chi_F$ is the characteristic function of
$F$, and ${\cal M}$ is the Hardy-Littlewood maximal operator, i.e.

\begin{equation}\label{eq2.3}
{\cal M}f(x):=\sup_{r>0}\frac{1}{|B(x,r)|}\int_{B(x,r)}|f(y)|\,dy,
\end{equation}

\noindent where $f$ is a locally integrable function and $B(x,r)$
stands for the ball with radius $r$ centered at $x\in\RR^n$.

\begin{lemma}\label{l2.1}\,(\cite{CMS}) Fix some $\alpha>0$. There exists $\gamma\in(0,1)$,
sufficiently close to 1, such that for every closed set $F$ whose
complement has finite measure and every non-negative function
$\Phi$ the following inequality holds:

\begin{equation}\label{eq2.4}
\int_{{\cal R}^\alpha(F^*)}\Phi(x,t)\,t^n\,dxdt\leq
C(\alpha,\gamma)\int_F\left[\int_{\Gamma(x)}\Phi(y,t)\,dydt\right]\,dx,
\end{equation}

\noindent where $F^*$ denotes the set of points of global
$\gamma$-density with respect to $F$.

Conversely,

\begin{equation}\label{eq2.5}
\int_F\left[\int_{\Gamma^{\alpha}(x)}\Phi(y,t)\,dydt\right]\,dx\leq
C(\alpha) \int_{{\cal R}^\alpha(F)}\Phi(x,t)\,t^n\,dxdt,
\end{equation}

\noindent for every closed set $F\subset\RR^n$ and every
non-negative function $\Phi$.
\end{lemma}

\begin{lemma}\label{l2.2}\, (\cite{CMS}) Consider the
operator

\begin{equation}\label{eq2.6}
S^\alpha F(x):=
\left(\dint_{\Gamma^\alpha(x)}|F(y,t)|^2\frac{dydt}{t^{n+1}}\right)^{1/2},
\end{equation}

\noindent where $\alpha>0$. There exists a constant $C>0$
depending on the dimension only such that

\begin{equation}\label{eq2.7}
\|S^\alpha F\|_{L^1(\RR^n)}\leq C \|S^1 F\|_{L^1(\RR^n)}.
\end{equation}
\end{lemma}

Both lemmas above are proved in \cite{CMS}.

\hskip 0.08in

Turning to the properties of the differential operator $L$ we
start with the off-diagonal estimates.  We say that the family of operators $\{S_t\}_{t>0}$ satisfies
$L^2$ off-diagonal estimates (``Gaffney estimates") if there exist
some constants $c, C,\beta >0$ such that for arbitrary closed sets
$E,F\subset \RR^n$

\begin{equation}\label{eq2.8}
\|S_tf\|_{L^2(F)}\leq C\,e^{-\left(\frac{{\rm
dist}\,(E,F)^2}{ct}\right)^\beta}\,\|f\|_{L^2(E)},
\end{equation}

\noindent for every $t>0$ and every $f\in L^2(\RR^n)$ supported in
$E$.

\begin{lemma}\label{l2.3}\, (\cite{HoMa})
If two families of operators, $\{S_t\}_{t>0}$ and $\{T_t\}_{t>0}$,
satisfy Gaffney estimates (\ref{eq2.8}) then so does
$\{S_tT_t\}_{t>0}$. Moreover, there exist $c, C>0$ such that for
arbitrary closed sets $E,F\subset \RR^n$

\begin{equation}\label{eq2.9}
\|S_sT_tf\|_{L^2(F)}\leq C\,e^{-\left(\frac{{\rm
dist}\,(E,F)^2}{c\max\{t,s\}}\right)^\beta}\,\|f\|_{L^2(E)},
\end{equation}

\noindent for all $t,s>0$ and all $f\in L^2(\RR^n)$ supported in
$E$.
\end{lemma}

\begin{lemma}\label{l2.4}\, (\cite{HoMa}, \cite{KatoMain})
The families

\begin{equation}\label{eq2.10}
\{e^{-tL}\}_{t>0},\qquad \{tLe^{-tL}\}_{t>0},\qquad \{t^{1/2}\nabla
e^{-tL}\}_{t>0},
\end{equation}

\noindent as well as

\begin{equation}\label{eq2.11}
\{(1+tL)^{-1}\}_{t>0},\qquad \{t^{1/2}\nabla (1+tL)^{-1}\}_{t>0},
\end{equation}

\noindent satisfy Gaffney estimates with $c, C>0$ depending on
$n$, $\lambda$ and $\Lambda$ only. For the operators in (\ref{eq2.10}),
$\beta = 1$, and in (\ref{eq2.11}), $\beta = 1/2$.\end{lemma}
We remark that it is well known from functional calculus and ellipticity (accretivity) that the operators in (\ref{eq2.10})--(\ref{eq2.11}) are bounded from $L^2(\RR^n)$ to
$L^2(\RR^n)$ uniformly in $t$.

\begin{lemma}\label{l2.5}
There exists $p_L, 1 \leq p_L <\frac{2n}{n+2}$ and $\widetilde{p}_L, \frac{2n}{n-2}< \widetilde{p}_L
<\infty,$ such that for every $p$ and $q$
with $p_L<p\leq q<\widetilde{p}_L$,
the family $\{e^{-tL}\}_{t>0}$
satisfies $L^p-L^q$ off-diagonal estimates, i.e. for arbitrary
closed sets $E,F\subset \RR^n$

\begin{equation}\label{eq2.12}
\|e^{-tL}f\|_{L^q(F)}\leq Ct^{\frac 12\left(\frac nq-\frac
np\right)}\,e^{-\frac{{\rm dist}\,(E,F)^2}{ct}}\,\|f\|_{L^p(E)},
\end{equation}

\noindent for every $t>0$ and every $f\in L^p(\RR^n)$ supported in
$E$. The operators $e^{-tL}$, $t>0$, are bounded from $L^p(\RR^n)$
to $L^q(\RR^n)$ with the norm $Ct^{\frac 12\left(\frac nq-\frac
np\right)}$ and from $L^p(\RR^n)$ to $L^p(\RR^n)$ with the norm
independent of $t$.

The statement of the Lemma remains valid with $\{e^{-tL}\}_{t>0}$
replaced by $\{(1+tL)^{-1}\}_{t>0}$, and with exponent $\beta = 1/2$ in the
exponential decay expression.
\end{lemma}

\bp For the heat semigroup the proof of the Lemma can be found in
\cite{AuscherSurvey} and the result for the resolvent can be
obtained following similar ideas. \ep

\noindent {\it Remark}.\,It has been shown in \cite{AuscherSurvey}
that the interval of $p$ such that the heat semigroup
is $L^p$-bounded, and the interval of $p,q$ such that it
enjoys the off-diagonal bound (\ref{eq2.12}), have the same interior. In
particular, there is no inconsistency between the definitions of
$p_L$ in the Introduction and in Lemma~\ref{l2.5}. We will
preserve this notation for $p_L$ and $\widetilde{p}_L$ throughout the paper.

\begin{lemma}\label{l2.6} The
operator

\begin{equation}\label{eq2.13}
S_h^Kf(x):=\left(\dint_{\Gamma(x)}|(t^2L)^K e^{-t^2L}f(y)|^2\,\frac
{dydt}{t^{n+1}}\right)^{1/2},\quad x\in\RR^n, \quad K\in\NN,
\end{equation}

\noindent is bounded in $L^p(\RR^n)$ for
$p\in\left(p_L,\widetilde{p}_L\right)$.
\end{lemma}

\bp The proof closely follows an analogous argument for vertical
square function (see \cite{AuscherSurvey}). We omit the details. \ep

Finally, the solutions of strongly parabolic and elliptic systems satisfy the
following versions of Caccioppoli inequality.

\begin{lemma}\label{l2.7} Suppose $Lu=0$ in
$B_{2r}(x_0)=\{x\in\RR^n:\,|x-x_0|<2r\}$. Then there exists
$C=C(\lambda, \Lambda)>0$ such that

\begin{equation}\label{eq2.14}
\int_{B_r(x_0)}|\nabla u(x)|^2\,dx\leq
\frac{C}{r^2}\int_{B_{2r}(x_0)}|u(x)|^2\,dx.
\end{equation}

\end{lemma}

\begin{lemma}\label{l2.8} Suppose $\partial_t u=-Lu$ in
$I_{2r}(x_0,t_0)$, where $I_r(x_0,t_0)=B_r(x_0)\times
[t_0-cr^2,t_0]$, $t_0>4cr^2$ and $c>0$. Then there exists
$C=C(\lambda, \Lambda, c)>0$ such that

\begin{equation}\label{eq2.15}
\dint_{I_r(x_0,t_0)}|\nabla u(x,t)|^2\,dxdt\leq
\frac{C}{r^2}\dint_{I_{2r}(x_0,t_0)}|u(x,t)|^2\,dxdt.
\end{equation}

\end{lemma}

\section{Sublinear operators in Hardy spaces}
\setcounter{equation}{0}

\begin{theorem}\label{t3.1}
Let $p_L<p\leq 2$ and assume that the sublinear operator

\begin{equation}\label{eq3.1}
T:L^p(\RR^n)\longrightarrow L^p(\RR^n)
\end{equation}

\noindent satisfies the following estimates. There exists
$M\in\NN$, $M>n/4$, such that for all closed sets $E$, $F$ in
$\RR^n$ with ${\rm dist }(E,F)>0$ and every $f\in L^p(\RR^n)$
supported in $E$

\begin{eqnarray}\label{eq3.2}
\|T(I-e^{-tL})^{M}f\|_{L^p(F)}&\leq & C\,\left(\frac{t}{{\rm
dist}\,(E,F)^2}\right)^{M}\,\|f\|_{L^p(E)},\qquad\forall\,t>0,\\[4pt]
\|T(tLe^{-tL})^{M}f\|_{L^p(F)}&\leq & C\,\left(\frac{t}{{\rm
dist}\,(E,F)^2}\right)^{M}\,\|f\|_{L^p(E)},\qquad\forall\,t>0.\label{eq3.3}
\end{eqnarray}

\noindent Then

\begin{equation}\label{eq3.4}
T:H^1_L(\RR^n)\longrightarrow L^1(\RR^n).
\end{equation}

\end{theorem}

\vskip 0.08 in


\noindent {\it Remark}.\,The proof below shows that
(\ref{eq3.1})--(\ref{eq3.3}) imply (\ref{eq3.4}) with the Hardy
space $H^1_L(\RR^n)$ defined by linear combinations of $(p,\eps,
M)$-molecules for the same values of $p$ and $M$ as in
(\ref{eq3.1})--(\ref{eq3.3}). We do not emphasize this fact as the
space $H^1_L(\RR^n)$ does not depend on the choice (within the stated allowable parameters)
of $p$, $\eps$
and $M$ in (\ref{eq1.8})--(\ref{eq1.11}) -- see
Corollary~\ref{c4.3}. \vskip 0.08 in

\bp Suppose that $T: L^p \rightarrow L^p$ is sublinear.  We claim that
for every $(p,\eps,M)$ molecule $m$, we have
\begin{equation}\label{eq3.6}
\|Tm\|_{L^1(\RR^n)}\leq C
\end{equation}
with constant $C$ independent of $m$.  Let us take this statement for granted
momentarily.  The conclusion of the theorem is then an immediate consequence of the following lemma,
whose proof we learned from P. Auscher.

\begin{lemma}\label{l3.sub} Fix $p,\epsilon,M$ within the allowable parameters,
with $p_L < p  \leq 2.$
Suppose that $T$ is a sublinear operator, bounded on
$L^p$, which
satisfies (\ref{eq3.6}) for all $(p,\epsilon,M)$
molecules.  Then $T$ extends to a bounded operator on $H^1_L$, and
$$\|Tf\|_{L^1} \leq C \|f\|_{H^1_L}.$$
\end{lemma}
\bp
For $f \in L^p \cap L^1$, define a norm
$$\|f\|_{T,1} \equiv \|Tf\|_{L^1} + \|f\|_{L^1},$$
and define an auxiliary space $\mathcal{H}_T^1$ as the completion of
$L^p$ with respect to this norm. Suppose that $f \in H^1_L$, so that
we may write $f = \sum_{i=0}^\infty \lambda_i m_i$ where the $m_i$
are $(p,\epsilon,M)$ molecules, and where $\sum_{i=0}^\infty
|\lambda_i| \approx \|f\|_{H^1_L}$.  Set $f_k = \sum_{i=0}^k
\lambda_i m_i$. Then $$\|f_k\|_{T,1} \leq C \sum_{i=0}^k |\lambda_i|
\leq C \|f\|_{H^1_L}.$$ Moreover, if $k' \leq k,$ then $$\|f_k -
f_{k'}\|_{T,1} \leq C \sum_{i=k'}^k |\lambda_i|\to 0$$ as $k,k' \to
\infty.$  By completeness, there exists $g \in \mathcal{H}^1_T$ with
$f_k \to g $ in $\mathcal{H}^1_T.$  In particular, $f_k \to g$ in
$L^1$, so that $g=f$ a.e. Thus, $f \in \mathcal{H}^1_T$, with
$$\|f\|_{T,1} = \lim_{k \to \infty} \|f_k\|_{T,1} \leq C \|f\|_{H^1_L}.$$
For $f \in L^p \cap H^1_L$ (which is dense in $H^1_L$), this implies that
$$\|Tf\|_{L^1} \leq C \|f\|_{H^1_L},$$
and we may extend T to all of $H^1_L$ by continuity.
\ep

We now turn to the proof of (\ref{eq3.6}).  To begin, we
decompose $\|Tm\|_{L^1(\RR^n)}$ in the following way

\begin{displaymath}
\|Tm\|_{L^1(\RR^n)}\leq C \sum_{i=0}^\infty
\|T(I-e^{-l(Q)^2L})^{M}(m\chi_{S_i(Q)})\|_{L^1(\RR^n)}
+C\|T[I-(I-e^{-l(Q)^2L})^{M}]m\|_{L^1(\RR^n)},
\end{displaymath}

\noindent where the family of annuli $\{S_i(Q)\}_{i=0}^\infty$ is
taken with respect to the cube $Q$ associated with $m$. Going
further,

\begin{eqnarray}\nonumber
&&\hskip -0.7cm\|T(I-e^{-l(Q)^2L})^{M}(m\chi_{S_i(Q)})\|_{L^1(\RR^n)}\\[4pt]
\nonumber &&\qquad \leq C\sum_{j=0}^\infty (2^{i+j}l(Q))^{n-\frac
np}
\|T(I-e^{-l(Q)^2L})^{M}(m\chi_{S_i(Q)})\|_{L^p(S_j(Q_i))} \\[4pt]
\nonumber &&\qquad \leq C\sum_{j=2}^\infty(2^{i+j}l(Q))^{n-\frac
np} \left(\frac{l(Q)^2}{{\rm
dist}\,(S_j(Q_i),S_i(Q))^2}\right)^{M}\,\|m\|_{L^p(S_i(Q))}\\[4pt]
&&\quad\qquad+\,C(2^{i}l(Q))^{n-\frac np}\|m\|_{L^p(S_i(Q))},
\label{eq3.7}
\end{eqnarray}

\noindent where the last inequality follows from
(\ref{eq3.1})--(\ref{eq3.2}) and uniform boundedness of
$\{e^{-tL}\}_{t>0}$ in $L^p(\RR^n)$.
Then by the properties of $(p,\eps,M)$-molecules the expression
above is bounded by

\begin{equation} \label{eq3.8}
C\sum_{j=2}^\infty
2^{i(-2M-\eps)}2^{j(n-n/p-2M)}+C\,2^{-i\eps}\leq C\,2^{-i\eps}.
\end{equation}

\noindent As for the remaining part,

\begin{equation}\label{eq3.9}
I-(I-e^{-l(Q)^2L})^{M}=\sum_{k=1}^{M}C^{M}_k(-1)^{k+1}e^{-kl(Q)^2L},
\end{equation}

\noindent where $C^{M}_k=\frac{{M}!}{({M}-k)!k!}$, $k=1,..., M$,
are binomial coefficients. Therefore,

\begin{eqnarray}\nonumber
&&\|T[I-(I-e^{-l(Q)^2L})^{M}]m\|_{L^1(\RR^n)}\leq C
\sup_{1\leq k\leq M}\|Te^{-kl(Q)^2L}m\|_{L^1(\RR^n)}\\[4pt]
&&\qquad \quad  \leq C \sup_{1\leq k\leq M}\left\|T\left(\frac
k{M}\, l(Q)^{2}Le^{-\frac
k{M}l(Q)^2L}\right)^{M}(l(Q)^{-2}L^{-1})^{M}m\right\|_{L^1(\RR^n)}.\label{eq3.10}
\end{eqnarray}

\noindent At this point we proceed as in
(\ref{eq3.7})--(\ref{eq3.8}) with $(l(Q)^{-2}L^{-1})^{M}m$ in
place of $m$, \break $\left(\frac k{M}\, l(Q)^{2}Le^{-\frac
k{M}l(Q)^2L}\right)^{M}$ in place of $(I-e^{-l(Q)^2L})^{M}$, and
using estimates (\ref{eq3.3}) and (\ref{eq1.9}). \ep

\vskip 0.08in

\noindent {\it Remark.}\,The result of Theorem~\ref{t3.1} holds
for $p\in(2,\widetilde{p}_L)$ as well. However, in that case one has to take
$M>\frac 12\,\left(n-\frac np\right)$.

\vskip 0.08in

For $f\in L^2(\RR^n)$ consider the following vertical version of
square function:

\begin{equation}\label{eq3.11}
g_hf(x):=\left(\int_0^\infty |t^2Le^{-t^2L}f(y)|^2\,\frac
{dt}{t}\right)^{1/2}.
\end{equation}

\begin{theorem}\label{t3.2}
The operators $g_h$ and $\nabla L^{-1/2}$ satisfy (\ref{eq3.2}) --
(\ref{eq3.3}) for $p=2$, $M>n/4$, and map $H^1_L(\RR^n)$ to
$L^1(\RR^n)$.
\end{theorem}

\bp The proof rests on ideas leading to Gaffney-type estimates for
Riesz transform (cf. Lemma~2.2 in \cite{HoMa}) and
Theorem~\ref{t3.1}.

First of all, the operators under consideration are obviously
sublinear. Further, $L^2$ boundedness of Riesz transform has been
proved in \cite{KatoMain} and boundedness of $g_h$ in $L^2(\RR^n)$
follows from the quadratic estimates for operators having bounded
holomorphic calculus (cf. \cite{ADM}).

Let us now address the inequalities (\ref{eq3.2})--(\ref{eq3.3})
for $p=2$ and the operator $T=g_h$. An argument for the Riesz
transform, viewed as

\begin{equation}\label{eq3.12}
\nabla L^{-1/2}=C \int_0^{\infty}\nabla e^{-sL}f\,\frac{ds}{\sqrt
s},
\end{equation}

\noindent  is completely analogous (cf. also \cite{HoMa}). Write

\begin{eqnarray}
&& \left\|g_h(I-e^{-tL})^{M}f\right\|_{L^2(F)}= C
\left\|\left(\int_0^\infty |sLe^{-sL}(I-e^{-tL})^{M}f|^2\,\frac
{ds}{s}\right)^{1/2}\right\|_{L^2(F)}\nonumber \\[4pt]
&& \qquad\leq C \left\|\left(\int_0^\infty
|sLe^{-s(M+1)L}(I-e^{-tL})^{M}f|^2\,\frac
{ds}{s}\right)^{1/2}\right\|_{L^2(F)}\nonumber \\[4pt]
&&\qquad \leq C \left(\int_0^t\left\|
sLe^{-s(M+1)L}(I-e^{-tL})^{M}f\right\|_{L^2(F)}^2\,\frac
{ds}{s}\right)^{1/2}  \nonumber \\[4pt]
&&\qquad \,\,+\,C\left(\int_t^\infty\left\|
sLe^{-s(M+1)L}(I-e^{-tL})^{M}f\right\|_{L^2(F)}^2\,\frac
{ds}{s}\right)^{1/2}=:I_1+I_2\label{eq3.13}.
\end{eqnarray}

\noindent We will analyze $I_1$ and $I_2$ separately. Expanding
$(I-e^{-tL})^{M}$ by binomial formula, one can see that

\begin{eqnarray}
&&\hskip -.7cm I_1 \leq C \left(\int_0^t\left\|
sLe^{-s(M+1)L}f\right\|_{L^2(F)}^2\,\frac {ds}{s}\right)^{1/2} + C
\sup_{1\leq k\leq M} \left(\int_0^t\left\|
sLe^{-s(M+1)L}e^{-ktL}f\right\|_{L^2(F)}^2\,\frac
{ds}{s}\right)^{1/2} \nonumber \\[4pt]
&&\hskip -.7cm\,\,\leq C \left(\int_0^t\left\|
sLe^{-s(M+1)L}f\right\|_{L^2(F)}^2\,\frac {ds}{s}\right)^{1/2} + C
\sup_{1\leq k\leq M} \left(\int_0^t\left\|
e^{-s(M+1)L}ktLe^{-ktL}f\right\|_{L^2(F)}^2\,\frac{sds}{t^2}\right)^{1/2} \nonumber \\[4pt]
&&\hskip -.7cm\,\,\leq C \left(\int_0^t e^{-\frac{{\rm
dist}\,(E,F)^2}{cs}}\,\frac {ds}{s}\right)^{1/2}\,\|f\|_{L^2(E)} +
C \sup_{1\leq k\leq M} \left(\frac{1}{t^2}\,e^{-\frac{{\rm
dist}\,(E,F)^2}{ct}}\, \int_0^t sds\right)^{1/2}\,\|f\|_{L^2(E)},
\nonumber
\end{eqnarray}

\noindent where we used Lemma~\ref{l2.4} and Lemma~\ref{l2.3}
repeatedly for the families $ktLe^{-ktL}$ and $m+1$ copies of
$e^{-sL}$ when estimating $\left\|
e^{-s(m+1)L}ktLe^{-ktL}f\right\|_{L^2(F)}$. It is not hard to see
that the expression above is bounded by $C\,\bigl(\frac{t}{{\rm
dist}\,(E,F)^2}\bigr)^{M}\,\|f\|_{L^2(E)}$, as desired.

Turning to the second integral,

\begin{equation}
I_2 \leq C\left(\int_t^\infty\left\|
sLe^{-sL}(e^{-sL}-e^{-(s+t)L})^{M}f\right\|_{L^2(F)}^2\,\frac
{ds}{s}\right)^{1/2} \label{eq3.14}.
\end{equation}

\noindent It was observed in \cite{HoMa} that

\begin{eqnarray}\nonumber
&&\left\| \frac st
(e^{-sL}-e^{-(s+t)L})g\right\|_{L^2(F)}=\left\|\frac st
\int_0^t\partial_r e^{-(s+r)L}g\,dr\right\|_{L^2(F)}\\[4pt]
&&\quad \leq C \frac st\int_0^t \left\|(s+r)L
e^{-(s+r)L}g\right\|_{L^2(F)}\,\frac{dr}{s+r} \nonumber\\[4pt]
&&\quad \leq  C \,\|g\|_{L^2(E)} \,\left(\frac st\int_0^t
e^{-\frac{{\rm
dist}\,(E,F)^2}{c(s+r)}}\,\frac{dr}{s+r}\right).\label{eq3.15}
\end{eqnarray}

\noindent But $s+r\approx s$ for $s\geq t$ and $r\in(0,t)$,
therefore the expression above does not exceed

\begin{equation}\label{eq3.16}
C \,\|g\|_{L^2(E)} e^{-\frac{{\rm dist}\,(E,F)^2}{cs}} \,\left(\frac
st\int_0^t \frac{dr}{s+r}\right) \leq C e^{-\frac{{\rm
dist}\,(E,F)^2}{cs}}\|g\|_{L^2(E)}.
\end{equation}

\noindent Now we multiply and divide the integrand in
(\ref{eq3.14}) by $\left(\frac st\right)^{2M}$ and use
Lemma~\ref{l2.3} for $sLe^{-sL}$ and $M$ copies of $\frac st
(e^{-sL}-e^{-(s+t)L})$ to get

\begin{equation}
I_2 \leq C\left(\int_t^\infty e^{-\frac{{\rm
dist}\,(E,F)^2}{cs}}\left(\frac ts\right)^{2M} \, \frac
{ds}{s}\right)^{1/2} \|f\|_{L^2(E)} \label{eq3.17}.
\end{equation}

\noindent Next, making a change of variables $r:=\frac{{\rm
dist}\,(E,F)^2}{cs}$, we can control the expression above by

\begin{equation}
I_2 \leq C\left(\int_0^\infty e^{-r}\left(\frac {tr}{{\rm
dist}\,(E,F)^2}\right)^{2M} \, \frac {dr}{r}\right)^{1/2}
\|f\|_{L^2(E)}\leq C\left(\frac {t}{{\rm dist}\,(E,F)^2}\right)^M
\|f\|_{L^2(E)}, \label{eq3.18}
\end{equation}

\noindent as desired.

This finishes the proof of (\ref{eq3.2}) for operator $g_h$. The
argument for (\ref{eq3.3}) follows essentially the same path. More
precisely, one needs to estimate the integrals $I_1$ and $I_2$ with
$(I-e^{-tL})^{M}$ replaced by $(tLe^{-tL})^{M}$. As for the first
one,

\begin{eqnarray}\nonumber
&& \left(\int_0^t\left\|
sLe^{-s(M+1)L}(tLe^{-tL})^{M}f\right\|_{L^2(F)}^2\,\frac
{ds}{s}\right)^{1/2}\\[4pt]\nonumber
&&\qquad \leq C \left(\frac{1}{t^2}\int_0^t\Bigl\|
e^{-s(M+1)L}(tLe^{-tL})^{M-1}
\left(\frac t2 Le^{-\frac{t}{2}L}\right)^{2}f\Bigr\|_{L^2(F)}\,sds\right)^{1/2}\\[4pt]
&&\qquad \leq C e^{-\frac{{\rm dist}\,(E,F)^2}{ct}}\,
\|f\|_{L^2(E)}\leq C\,\left(\frac{t}{{\rm
dist}\,(E,F)^2}\right)^{M}\,\|f\|_{L^2(E)}, \label{eq3.19}
\end{eqnarray}

\noindent by Lemma~\ref{l2.3}. Concerning to the analogue of
$I_2$, one can write

\begin{eqnarray}\nonumber
&&\left(\int_t^\infty\left\|
sLe^{-s(M+1)L}(tLe^{-tL})^{M}f\right\|_{L^2(F)}^2\,\frac
{ds}{s}\right)^{1/2}\\[4pt]
&&\qquad \leq \left(\int_t^\infty\Bigl(\frac ts\Bigr)^{2M}\left\|
sLe^{-sL}(sLe^{-(t+s)L})^{M}f\right\|_{L^2(F)}^2\,\frac
{ds}{s}\right)^{1/2}.\label{eq3.20}
\end{eqnarray}

\noindent At this stage, similarly to (\ref{eq3.15}),

\begin{eqnarray}\nonumber
&&\hskip -.7cm \left\|
sLe^{-(s+t)L}g\right\|_{L^2(F)}=\left\|\frac st\, e^{-sL}
\int_0^t\partial_r(rL e^{-rL})g\,dr\right\|_{L^2(F)}\\[4pt]
\nonumber &&\hskip -.7cm\quad \leq C \left\|\frac st \,e^{-sL}
\int_0^t(Le^{-rL} - rL^2e^{-rL})g\,dr\right\|_{L^2(F)}\\[4pt]
\nonumber &&\hskip -.7cm\quad \leq C \frac st
\int_0^t\|(r+s)Le^{-(r+s)L}g\|_{L^2(F)}\,\frac{dr}{r+s} + C \frac
st
\int_0^t\|(r+s)^2L^2e^{-(r+s)L}g\|_{L^2(F)}\,\frac{rdr}{(r+s)^2}\\[4pt]
&&\hskip -.7cm  \quad\leq C \,\|g\|_{L^2(E)} \,\left(\frac
st\int_0^t e^{-\frac{{\rm
dist}\,(E,F)^2}{c(s+r)}}\,\frac{dr}{s+r}\right)\leq C
e^{-\frac{{\rm dist}\,(E,F)^2}{cs}}\|g\|_{L^2(E)} \label{eq3.21},
\end{eqnarray}

\noindent and following (\ref{eq3.17})--(\ref{eq3.18}) we complete
the proof. \ep

\vskip 0.08in

\noindent {\it Remark.}\, Using the same line of reasoning and
$L^p-L^2$ off-diagonal estimates ($p_L<p\leq 2$) one can
show for all closed sets $E$, $F$ in $\RR^n$ with ${\rm dist
}(E,F)>0$ and every $f\in L^p(\RR^n)$ supported in $E$

\begin{eqnarray}\label{eq3.22}
\|g_h(I-e^{-tL})^{M}f\|_{L^2(F)}&\leq & C\,t^{\frac 12\left(\frac
n2-\frac np\right)}\,\left(\frac{t}{{\rm
dist}\,(E,F)^2}\right)^{M}\,\|f\|_{L^p(E)},\qquad\forall\,t>0,\\[4pt]
\|g_h(tLe^{-tL})^{M}f\|_{L^2(F)}&\leq & C\,t^{\frac 12\left(\frac
n2-\frac np\right)}\,\left(\frac{t}{{\rm
dist}\,(E,F)^2}\right)^{M}\,\|f\|_{L^p(E)},\qquad\forall\,t>0.\label{eq3.23}
\end{eqnarray}

\section{Characterization by the square function associated to the heat semigroup}
\setcounter{equation}{0}

To begin, let us note that the space $H^1_L$, defined by means of
molecular decompositions, is complete. We learned the following
proof of this fact from E. Russ. We first require a well known
completeness criterion from functional analysis.

\begin{lemma}\label{l1.complete}
Let $X$ be a normed space which enjoys the property that $\sum x_k$
converges in $X$, whenever $\sum \|x_k\| < \infty.$  Then, $X$ is
complete.
\end{lemma}

The lemma is well known, and we omit the proof.

Let us now use the lemma to establish completeness of $H^1_L$.  To
this end, we suppose that $f_k \in H^1_L$, and that $\sum
\|f_k\|_{H^1_L} < \infty$. Given the former fact, there exists for
each $k$ a molecular decomposition $f_k = \sum_{i=0}^\infty
\lambda_i^k m_i^k$, with $\sum_{i=0}^\infty |\lambda_i^k| \approx
\|f_k\|_{H^1_L}$. Thus,
$$\sum_{i,k} |\lambda_i^k| \approx \sum_k \|f_k\|_{H^1_L} < \infty.$$
Consequently, the sum $ \sum f_k = \sum_{k,i} \lambda_i^k m_i^k $
converges in $H^1_L,$ as desired. \ep

For technical reasons related to the characterization of the dual of
our Hardy space, we shall need to work also with a modified version
of the molecular representations. Given $p\in(p_L,\widetilde{p}_L)$,
$\eps>0$, $M>n/4$, and $\delta
> 0,$ we say that $f = \sum \lambda_j m_j$ is a
$\delta$-representation of $f$ if $\{\lambda_j \}_{j=0}^\infty \in
\ell^1$ and each $m_j$ is a $(p,\eps,M)$-molecule adapted to a cube
$Q_j$ of side length at least $\delta$. We set
$$H^1_{L,\delta}(\RR^n) \equiv \{ f \in L^1(\RR^n): f \textrm{ has a
$\delta$-representation}\}.$$ Observe that a $\delta$-representation
is also a $\delta^{'}$-representation for all $\delta^{'} < \delta.$
Thus, $H^1_{L,\delta} \subseteq H^1_{L,\delta^{'}}$  for $0<
\delta^{'} < \delta.$ Set $$ \widehat{H}^1_L(\RR^n) \equiv
\cup_{\delta>0} H^1_{L,\delta}(\RR^n),$$ and define
\begin{equation}\label{eq1.10a}\|f\|_{\widehat{H}^1_L(\RR^n)} \equiv \inf \,\left\{ \sum_{j=0}^\infty |\lambda_j| :
f = \sum_{j=0}^\infty \lambda_j m_j \textrm{ is a
$\delta$-representation for some } \delta > 0
\right\}.\end{equation} We note that for $f \in \widehat{H}^1_L,$
$\|f\|_{L^1} \leq C \|f\|_{\widehat{H}^1_L},$ since in particular,
$\|m\|_{L^1} \leq C$ for every molecule $m$ (of course, a similar
statment is true for $H^1_L$). Now for $f \in \widehat{H}^1_L$, set
\begin{equation}\label{eq1.11a}\|f\|_{\widetilde{H}^1_L(\RR^n)} \equiv \|f\|_{\widehat{H}^1_L(\RR^n)}
+ \|f\|_{L^1(\RR^n)} \approx
\|f\|_{\widehat{H}^1_L(\RR^n)},\end{equation} and define
$\widetilde{H}^1_L$ as the completion of $\widehat{H}^1_L$ with
respect to this norm. We shall show that
$$H^1_L = \widetilde{H}^1_L,$$ for any choice of $p, \eps$ and $M$,
within the allowable parameters.  By definition,
\begin{equation}\label{eq1.compare}
\widetilde{H}^1_L \subseteq H^1_L, \,\, \textrm{ and } \,\,
\|f\|_{H^1_L} \leq \|f\|_{\widetilde{H}^1_L}
\end{equation}
whenever $f \in \widetilde{H}^1_L$.  The non-trivial issue will be
to verify that these statements can be reversed (up to a
multiplicative constant).

\begin{theorem}\label{t4.1}  The spaces $H^1_L$, $\widetilde{H}^1_L$ and $H^1_{S_h}$ are the same;  in particular,
\begin{equation}\label{eq4.compare}\|f\|_{H^1_L(\RR^n)}\approx\|f\|_{\widetilde{H}^1_L(\RR^n)}\approx\|S_h f\|_{L^1(\RR^n)} + \|f\|_{L^1(\RR^n)}.\end{equation}
\end{theorem}

In light of (\ref{eq1.compare}), the theorem is an immediate
consequence of the following lemma.

\begin{lemma}\label{lmain}  We have the containments $H^1_L \subseteq H^1_{S_h} \subseteq \widetilde{H}^1_L $.
Moreover,

(i) If $f \in L^2 \cap H^1_{S_h}$, then $f$ is the limit in
$\widetilde{H}^1_L$ of $f_N \in \widehat{H}_L^1$.  Furthermore, for
every $\eps>0$, $M\in\NN$, $p_L<p<\widetilde{p}_L$ there exists a
family of $(p,\eps,M)$-molecules $\{m_i\}_{i=0}^\infty$ and a
sequence of numbers $\{\lambda_i\}_{i=0}^\infty$ such that $f$ can
be represented in the form $f=\sum_{i=0}^\infty \lambda_i m_i$, with

\begin{equation}\label{eq4.1}
\|f\|_{\widetilde{H}^1_L(\RR^n)}\leq C \sum_{i=0}^\infty |\lambda_i|\leq
C\|f\|_{H^1_{S_h}(\RR^n)}.
\end{equation}

(ii) Conversely, given $\eps>0$, $M>n/4$ and
$p_L<p<\widetilde{p}_L$, let $f \equiv \sum_{i=0}^\infty \lambda_i
m_i$, where $\{m_i\}_{i=0}^\infty$ is a family of
$(p,\eps,M)$-molecules and $\sum_{i=0}^\infty |\lambda_i|<\infty.$
Then the series $\sum_{i=0}^\infty \lambda_i m_i$ converges in
$H^1_{S_h}(\RR^n)$ and
$$\Bigl\|\sum_{i=0}^\infty \lambda_i m_i\Bigr\|_{H^1_{S_h}(\RR^n)}\leq
C\sum_{i=0}^\infty |\lambda_i|.$$
In particular, we have that
\begin{equation}\label{eq4.2}\|f\|_{H^1_{S_h}(\RR^n)}\leq C \|f\|_{H^1_L(\RR^n)}.\end{equation}

\end{lemma}

The ideas in this proof go back to the decomposition of the
tent spaces in \cite{CMS} as well as to the atomic decomposition of the
classical Hardy spaces (see \cite{Wilson}).

\bigskip

\bp
\noindent {\bf Step I}. Let $f\in H^1_{S_h}(\RR^n)\cap L^2(\RR^n)$.
We shall construct
a family $\{f_N\}_{N=1}^\infty \subseteq \widehat{H}^1_L$ such that
$f_N \rightarrow f$ in $\widetilde{H}^1_L$ and in $H^1_{S_h}$, with
$$\sup_N\|f_N\|_{\widetilde{H}^1_L(\RR^n)} \leq C \|f\|_{H^1_{S_h}(\RR^n)}.$$  In particular, this will show that
$f \in \widetilde{H}^1_L$, with norm controlled by $\|f\|_{H^1_{S_h}}$.  The
claimed molecular decomposition will be established in the course
of the proof.

We start with a suitable version of the Calder\'{o}n reproducing
formula. By $L^2$ functional calculus, for every $f\in L^2(\RR^n)$
one can write

\begin{eqnarray}\label{eq4.3}
f&=&C_M\,\int_0^\infty (t^2Le^{-t^2L})^{M+2}f\,\frac{dt}{t}\nonumber\\&=& C_M\lim_{N\rightarrow \infty}
\int_{1/N}^N (t^2Le^{-t^2L})^{M+2}f\,\frac{dt}{t}\equiv \lim_{N\rightarrow \infty} f_N,\label{eq4.3}
\end{eqnarray}

\noindent with the integral converging in $L^2(\RR^n)$.

Now define the family of sets ${O}_k:= \{x\in\RR^n:\,S_hf(x)>2^k\}$,
$k\in\ZZ$, and consider ${O}_k^*:=\{x\in\RR^n:\, {\mathcal
M}(\chi_{O_k})>1-\gamma\}$ for some fixed $0<\gamma<1$. Then
$O_k\subset O_k^*$ and $|O_k^*|\leq C(\gamma)|O_k|$ for every
$k\in\ZZ$. Next let $\{Q_k^j\}_j$ be a Whitney decomposition of
$O_k^*$ and $\widehat O_k^*$ be a tent region, that is

\begin{equation}\label{eq4.4}
\widehat O_k^*:=\{(x,t)\in\RR^n\times(0, \infty):\,{\rm
dist}(x,\,^cO_k^*)\geq t\}.
\end{equation}

\noindent For every $k,j\in\ZZ$ we define

\begin{equation}\label{eq4.5}
T_k^j:=\left(Q_k^j\times(0,\infty)\right)\cap \widehat O_k^* \cap
\,^c\widehat O_{k+1}^*,
\end{equation}

\noindent and then recall the formula (\ref{eq4.3}) to write

\begin{equation}\label{eq4.6}
f_N=C_M\,\sum_{j,k\in\ZZ}\int_{1/N}^N
(t^2Le^{-t^2L})^{M+1}\left(\chi_{T_k^j}t^2Le^{-t^2L}\right)f\,\frac{dt}{t}
=:\sum_{j,k\in\ZZ}\lambda_{k}^jm_k^j(N),
\end{equation}

\noindent where $\lambda_k^j=C_M2^k|Q_k^j|$ and

\begin{equation}\label{eq4.7}
m_k^j(N)=\frac{1}{\lambda_k^j}\int_{1/N}^N
(t^2Le^{-t^2L})^{M+1}\left(\chi_{T_k^j}t^2Le^{-t^2L}\right)f\,\frac{dt}{t}.
\end{equation}

We claim that, up to a normalization by a harmless multiplicative constant, the $m_k^j(N)$ are molecules.  Assuming the claim,
we note that by definition of $T_k^j$, $m_k^j(N) = 0$ if $C\ell(Q_k^j) < 1/N,$
so that (\ref{eq4.6}) is a $\delta$-representation with $\delta \approx 1/N$.  Thus,
once the claim is established, we shall have

\begin{eqnarray}
\sup_N\|f_N\|_{\widetilde{H}^1_L(\RR^n)}\leq
C\sum_{j,k\in\ZZ}{\lambda_k^j}&=&C\sum_{j,k\in\ZZ}2^k|Q_k^j| \leq
C\sum_{k\in\ZZ}2^k|O_k^*|\nonumber\\[4pt]&\leq&
C\sum_{k\in\ZZ}2^k|O_k|\leq C\|S_hf\|_{L^1(\RR^n)},\label{eq4.8}
\end{eqnarray}

\noindent as desired. Let us now prove the claim:  that is, we will show that for every
$j,k\in\ZZ$, and $N\in\NN$ the function $C^{-1}m_k^j(N)$ is a $(p,\eps,M)$- molecule
associated with the cube $Q_k^j$, $2\leq p<\widetilde{p}_L$, for some harmless constant $C$. The case $p\leq
2$ follows from this one by H\"older's inequality.

To this end, fix $j,k\in\ZZ$ and $i\in\NN\cup\{0\}$ and consider
some $h\in L^{p'}(S_i(Q_k^j))$ such that
$\|h\|_{L^{p'}(S_i(Q_k^j))}=1$, $\frac 1p+\frac 1{p'}=1$. Then

\begin{eqnarray}\nonumber
&&\hskip -.5cm \left|\int_{\RR^n}m_k^j(N)(x)\overline{h(x)}\,dx\right|\leq
\frac{1}{\lambda_k^j}\dint_{T_k^j}
\left|t^2Le^{-t^2L}f(x)\left((t^2Le^{-t^2L})^{M+1}\right)^*h(x)\right|\,\frac{dt}{t}\,dx\nonumber\\[4pt]
&&\quad \leq \frac{C}{\lambda_k^j}\int_{\,^c
O_{k+1}}\dint_{\Gamma(x)}
\chi_{\left(Q_k^j\times(0,\infty)\right)\cap \widehat O_k^*}(y,t)
\left|t^2Le^{-t^2L}f(y)\left((t^2Le^{-t^2L})^{M+1}\right)^*h(y)\right|\,\frac{dtdy}{t^{n+1}}\,dx\nonumber\\[4pt]
&&\quad \leq \frac{C}{\lambda_k^j}\left(\int_{\,^c O_{k+1}\cap
cQ_k^j}\left(\dint_{\Gamma(x)}\bigl|t^2Le^{-t^2L}f(y)\bigr|^2\,\frac{dydt}{t^{n+1}}\right)^{\frac
p2 }\,dx\right)^{\frac 1p}\times\nonumber\\[4pt]
&&\quad\qquad \times\left(\int_{\,^c O_{k+1}\cap
cQ_k^j}\left(\dint_{\Gamma(x)\cap
\left(Q_k^j\times(0,\infty)\right)\cap \widehat
O_k^*}\left|\left((t^2Le^{-t^2L})^{M+1}\right)^*h(y)\right|^2\,\frac{dydt}{t^{n+1}}\right)^{\frac
{p'}{2}}\,dx\right)^{\frac{1}{p'}}
\nonumber\\[4pt]
&&\quad =: I_1\times I_2.\label{eq4.9}
\end{eqnarray}

\noindent Some comments are in order here. For the second
inequality we used Lemma~\ref{l2.1} with $F=\,^c O_{k+1}$ (so that
$F^*=\,^cO_{k+1}^*$ and ${\cal R}^1(F^*)=\,^c\widehat O_{k+1}^*$)
and

\begin{equation}\label{eq4.10}
\Phi(x,t):=\chi_{\left(Q_k^j\times(0,\infty)\right)\cap \widehat
O_k^*}(y,t)
\left|t^2Le^{-t^2L}f(x)\left((t^2Le^{-t^2L})^{M+1}\right)^*h(x)\right|t^{-n-1}.
\end{equation}

\noindent The third estimate above is based on H\"older's
inequality and the fact that whenever \break $y\in \Gamma(x)\cap
\left(Q_k^j\times(0,\infty)\right)\cap \widehat O_k^*$ we have
$x\in c\,Q_k^j$, where the constant $c$ is related to the implicit
constant in Whitney decomposition. Without loss of generality we
will assume that $c\leq 3$.

Observe now that

\begin{equation}\label{eq4.11}
I_1\leq C\frac{1}{\lambda_k^j}\left(\int_{cQ_k^j\cap \,^c
O_{k+1}}(S_hf(x))^p\,dx\right)^{1/p},
\end{equation}

\noindent and $S_hf(x)$ is bounded by $2^{k+1}$ for every $x\in
\,^c O_{k+1}$. Therefore

\begin{equation}\label{eq4.12}
I_1 \leq C \frac{1}{\lambda_k^j}\,2^{k+1}\,|Q_k^j|^{\frac 1p}\leq
C\, |Q_k^j|^{\frac 1p-1}.
\end{equation}

Turning to $I_2$, recall that ${\rm supp}\,\,h\subset S_i(Q_k^j)$.
Then to handle $i\leq 4$ it is enough to notice that

\begin{equation}\label{eq4.13}
I_2\leq C\left\|S^{M+1}_{h}h\right\|_{L^{p'}(\RR^n)}\leq
C\|h\|_{L^{p'}(S_i(Q_k^j))}\leq C,
\end{equation}

\noindent using Lemma~\ref{l2.6}. Then

When $i\geq 5$, we proceed as follows. We invoke H\"older's
inequality to estimate $L^{p'}$ norm by $L^2$ norm and then apply
(\ref{eq2.5}) with $F=\,^c O_{k+1}\cap c\,Q_k^j$ and

\begin{equation}\label{eq4.14}
\Phi(y,t)=\chi_{(0,\,c\,l(Q_k^j))}(t)\left|\left((t^2Le^{-t^2L})^{M+1}\right)^*h(y)\right|^2
t^{-n-1}.
\end{equation}

\noindent Then

\begin{equation}\label{eq4.15}
I_2\leq C|Q_k^j|^{\frac 1{p'}-\frac
12}\left(\int_{3c\,Q_k^j}\int_0^{c\,l(Q_k^j)}
\left|\left((t^2Le^{-t^2L})^{M+1}\right)^*h(x)\right|^2\,\frac{dxdt}{t}\right)^{\frac{1}{2}}
\end{equation}

By Lemmas~\ref{l2.5} and \ref{l2.3} applied to the operator $L^*$
the expression in (\ref{eq4.15}) is bounded by

\begin{eqnarray}\nonumber
&&\hskip -1cm C|Q_k^j|^{\frac 1{p'}-\frac
12}\left(\int_0^{c\,l(Q_k^j)}e^{-\frac{{\rm dist}\,(3cQ_k^j,
S_i(Q_k^j))^2}{ct^2}}\,t^{2\left(\frac n2-\frac n{p'}\right)}\,
\frac{dt}{t}\right)^{1/2}\|h\|_{L^{p'}(S_i(Q_k^j))}\\[4pt]
&&\hskip -1cm \qquad \leq C|Q_k^j|^{\frac 1{p'}-\frac
12}\left(\int_0^{c\,l(Q_k^j)}\left(\frac{t}{2^il(Q_k^j)}\right)^{2\left(\frac
n{p'}+\eps\right)}\,t^{2\left(\frac n2-\frac n{p'}\right)}\,
\frac{dt}{t}\right)^{1/2}\leq C 2^{-i\left(n-\frac
n{p}+\eps\right)}.\label{eq4.16}
\end{eqnarray}

All in all,

\begin{equation}\label{eq4.17}
\left|\int_{\RR^n}m_k^j(N)(x)\overline{h(x)}\,dx\right|\leq C2^{-i(n-n/p
+\eps)}|Q_k^j|^{1/p-1}
\end{equation}

\noindent for every $h\in L^{p'}(S_i(Q_k^j))$ with
$\|h\|_{L^{p'}(S_i(Q_k^j))}=1$. Taking supremum over all such $h$
we arrive at (\ref{eq1.8}).

The condition (\ref{eq1.9}) can be verified directly applying
$(l(Q_k^j)^{-2}L^{-1})^k$, $1\leq k\leq M$, to the molecule and
arguing along the lines (\ref{eq4.9})--(\ref{eq4.17}). A few
modifications relate solely to $I_2$ which should be majorized by

\begin{displaymath}
\left(\int_{\,^c O_{k+1}\cap cQ_k^j}\left(\dint_{\Gamma(x)\cap
\left(Q_k^j\times(0,\infty)\right)\cap \widehat
O_k^*}\left|\left(e^{-kt^2L}(t^2Le^{-t^2L})^{M+1-k}\right)^*h(y)\right|^2\,\frac{dydt}{t^{n+1}}\right)^{\frac
{p'}{2}}\,dx\right)^{\frac{1}{p'}},
\end{displaymath}

\noindent the rest of the argument follows verbatim.  We have therefore established that
$f_N \in \widehat{H}^1_L$, and that $f_N$ satisfies the desired bounds in $\widetilde{H}^1_L$, uniformly
in $N$.  It remains to verify that $f_N \rightarrow f$ in $\widetilde{H}^1_L$ and in $H^1_{S_h}$, and also that
$f = \sum \lambda_k^j m_k^j$, where the $m_k^j = \lim_{N\to \infty} m_k^j (N)$ exist and are molecules
(up to a harmless normalization).  We defer consideration of this matter for the moment, and proceed to:

\vskip 0.08 in \noindent {\bf Step II}. Let us now move to part
$(ii)$ of the lemma. To begin, observe that it is enough to
consider the case $p\leq 2$.

We will employ ideas similar to those in the proof of
Theorem~\ref{t3.1}. We recall that $S_h:L^p \to L^p, p_L < p <
\tilde{p}_L$ (Lemma~\ref{l2.6}).  By Lemma \ref{l3.sub}, it will
therefore be enough to show that $S_h$ maps allowable $(p,\eps,M)$
molecules uniformly into $L^1$. In fact the geometry of the operator
$S_h$ does not allow to obtain the estimates
(\ref{eq3.2})--(\ref{eq3.3}) for arbitrary closed sets $E$ and $F$,
however one can observe that we do not use
(\ref{eq3.2})--(\ref{eq3.3}) in full generality either. More
specifically, these inequalities are required to estimate

\begin{equation}\label{eq4.18}
\|S_h(I-e^{-l(Q)^2L})^{M}(m\chi_{S_i(Q)})\|_{L^p(S_j(Q_i))}
\end{equation}

\noindent in (\ref{eq3.7}) and

\begin{equation}\label{eq4.19}
\sup_{1\leq k\leq M}\left\|S_h\left(\frac k{M}\,
l(Q)^{2}Le^{-\frac k{M}l(Q)^2L}\right)^{M}
[\chi_{S_i(Q)}(l(Q)^{-2}L^{-1})^{M}m]\right\|_{L^p(S_j(Q_i))}
\end{equation}

\noindent later on ($i,j=0,1,...$). In what follows we will
analyze (\ref{eq4.18})--(\ref{eq4.19}) directly and insert the
result into the proof of Theorem~\ref{t3.1}.

Since $S_h$ and $(I-e^{-l(Q)^2L})^{M}$ are bounded in
$L^p$ (in the latter case with constant independent of $l(Q)$), we can write for
$j=0,1,2$,

\begin{equation}\label{eq4.20}
\|S_h(I-e^{-l(Q)^2L})^{M}(m\chi_{S_i(Q)})\|_{L^p(S_j(Q_i))}\leq C
\|m\|_{L^p(S_i(Q))}.
\end{equation}

\noindent Assume now that $j\geq 3$. By H\"older's inequality and
Lemma~\ref{l2.1}

\begin{eqnarray}\nonumber
&&\hskip -0.7cm \|S_h(I-e^{-l(Q)^2L})^{M}(m\chi_{S_i(Q)})\|_{L^p(S_j(Q_i))}^2\\[4pt]\nonumber
&&\hskip -0.7cm\,\,\, \leq C (2^{i+j}l(Q))^{2\left(\frac np-\frac
n2\right)}\dint_{{\cal R}(S_j(Q_i))}
|t^2Le^{-t^2L}(I-e^{-l(Q)^2L})^{M}(m\chi_{S_i(Q)})(x)|^2\,\frac{dtdx}{t}\\[4pt]\nonumber
&&\hskip -0.7cm\,\,\, \leq C (2^{i+j}l(Q))^{2\left(\frac np-\frac
n2\right)}\int_{\RR^n\setminus Q_{j-2+i}}\int_0^\infty
|t^2Le^{-t^2L}(I-e^{-l(Q)^2L})^{M}(m\chi_{S_i(Q)})(x)|^2\,\frac{dtdx}{t}+\\[4pt]\nonumber
&&\hskip -0.7cm\,\,\, +C(2^{i+j}l(Q))^{2\left(\frac np-\frac
n2\right)}\sum_{k=0}^{j-2}\int_{S_k(Q_i)}\int_{(2^{j-1}-2^k)2^il(Q)}^\infty
|t^2Le^{-t^2L}(I-e^{-l(Q)^2L})^{M}(m\chi_{S_i(Q)})(x)|^2\,\frac{dtdx}{t}\\[4pt]
&&\hskip -0.7cm\,\,\, =:I+\sum_{k=0}^{j-2}I_k. \nonumber
\end{eqnarray}

\noindent We observe that

\begin{eqnarray}\nonumber
I &\leq & C (2^{i+j}l(Q))^{2\left(\frac np-\frac n2\right)}
\|g_h(I-e^{-l(Q)^2L})^{M}(m\chi_{S_i(Q)})\|_{L^2(\RR^n\setminus
Q_{j-2+i})}^2\\[4pt]
&\leq & C\,(2^{i+j}l(Q))^{2\left(\frac np-\frac
n2\right)}l(Q)^{2\left(\frac n2-\frac
np\right)}\,\left(\frac{l(Q)^2}{{\rm dist}\,(S_i(Q),\RR^n\setminus
Q_{j-2+i})^2}\right)^{2M}\,\|m\|_{L^p(S_i(Q))}^2\nonumber\\[4pt]
&\leq & C\,\Bigl(\frac{1}{2^{i+j}}\Bigr)^{4M+2\left(\frac n2-\frac
np\right)}\,\|m\|_{L^p(S_i(Q))}^2, \label{eq4.21}
\end{eqnarray}

\noindent where the second inequality follows from (\ref{eq3.22}).
Turning to $I_k$, $k=0,1,...,j-2$, we make a change of variables
$s:=t^2/(m+1)$, so that $s\geq
[(2^{j-1}-2^k)2^il(Q)]^2/(m+1)\approx [2^{i+j}l(Q)]^2$ and

\begin{eqnarray}
I_k&\leq& C (2^{i+j}l(Q))^{2\left(\frac np-\frac
n2\right)}\int_{c\,[2^{i+j}l(Q)]^2}^\infty
\|sLe^{-(m+1)sL}(I-e^{-l(Q)^2L})^{M}(m\chi_{S_i(Q)})\|^2_{L^2(S_k(Q_i))}\,\frac{ds}{s}\nonumber\\[4pt]
\nonumber &\leq& C (2^{i+j}l(Q))^{2\left(\frac np-\frac
n2\right)}\int_{c\,[2^{i+j}l(Q)]^2}^\infty
\Bigl(\frac{l(Q)^2}{s}\Bigr)^{2M}\times\\[4pt]
&&\qquad \times
\left\|sLe^{-sL}\Bigl[\frac{s}{l(Q)^2}(e^{-sL}-e^{-(l(Q)^2+s)L})\Bigr]^{M}
(m\chi_{S_i(Q)})\right\|^2_{L^2(S_k(Q_i))}\,\frac{ds}{s}.
\label{eq4.22}
\end{eqnarray}

\noindent At this point we apply (\ref{eq3.15})--(\ref{eq3.16})
with $t=l(Q)^2$ combined with $L^p-L^2$ off-diagonal estimates for
$sLe^{-sL}$ and obtain

\begin{eqnarray}\nonumber
I_k&\leq& C (2^{i+j}l(Q))^{2\left(\frac np-\frac n2\right)}
\|m\|_{L^2(S_i(Q))}^2 \\[4pt]
&&\qquad \times \int_{c\,[2^{i+j}l(Q)]^2}^\infty s^{\left(\frac
n2-\frac np\right)}
 e^{-\frac{{\rm
dist}\,(S_i(Q),S_k(Q_i))^2}{cs}}\,
\Bigl(\frac{l(Q)^2}{s}\Bigr)^{2M}\,\frac {ds}{s}\, .
\label{eq4.23}
\end{eqnarray}

\noindent Therefore,

\begin{eqnarray}
I_k&\leq& C \left(\int_{c\,[2^{i+j}l(Q)]^2}^\infty
\Bigl(\frac{l(Q)^2}{s}\Bigr)^{2M}\,\frac {ds}{s}\right)\,
\|m\|_{L^p(S_i(Q))}^2 \leq C \Bigl(\frac{1}{2^{i+j}}\Bigr)^{4M}
\|m\|_{L^p(S_i(Q))}^2. \label{eq4.24}
\end{eqnarray}

\noindent  Combining (\ref{eq4.21}) and (\ref{eq4.24}), one can
see that

\begin{equation}
\|S_h(I-e^{-l(Q)^2L})^{M}(m\chi_{S_i(Q)})\|_{L^p(S_j(Q_i))}^2\leq
C j\left(\frac {1}{2^{i+j}}\right)^{4M+2\left(\frac n2-\frac
np\right)} \|m\|_{L^p(S_i(Q))}^2. \label{eq4.25}
\end{equation}

\noindent And finally, substituting (\ref{eq4.25}) to
(\ref{eq3.7}),

\begin{eqnarray}\nonumber
&&\hskip -0.7cm\|S_h(I-e^{-l(Q)^2L})^{M}(m\chi_{S_i(Q)})\|_{L^1(\RR^n)}\\[4pt]
\nonumber && \leq C\sum_{j=3}^\infty(2^{i+j}l(Q))^{n-\frac
np}\,\sqrt j\, \left(\frac {1}{2^{i+j}}\right)^{2M+\left(\frac
n2-\frac np\right)} \|m\|_{L^p(S_i(Q))}+C(2^{i}l(Q))^{n-\frac
np}\|m\|_{L^p(S_i(Q))}
\\[4pt]
&& \leq C2^{i\left(\frac np-\frac n2
-2M-\eps\right)}\,\sum_{j=3}^\infty \sqrt
j\,2^{j(n/2-2M)}+C\,2^{-i\eps}\leq C\,2^{-i\eps}.\nonumber
\end{eqnarray}

Similar argument provides an estimate for
$\|S_h[I-(I-e^{-l(Q)^2L})]^{M}m\|_{L^1(\RR^n)}$. This time we use
(\ref{eq3.23}) to control an analogue of (\ref{eq4.21}), and
(\ref{eq3.21}) instead of (\ref{eq3.15}) at the step corresponding
to (\ref{eq4.22})--(\ref{eq4.23}).  Thus, we have established part $(ii)$ of the lemma.

To complete part $(i)$, it remains to show that $f_N \to f$ in $\widetilde{H}^1_L$ and in $H^1_{S_h}$,
and that $f$ has the claimed molecular decomposition.  To this end, we recall that $f_N \to f$ in $L^2$.
We claim that $\{f_N\}$ is a Cauchy sequence in $\widetilde{H}^1_L$.  Let us
postpone establishing this claim until
the end of the proof.  Assuming the claim, we see that there exists $g\in \widetilde{H}^1_L$ such that  $f_N \to g$ in $\widetilde{H}^1_L$, and, in particular, in $L^1$.  By taking subsequences which converge a.e., we see that the $L^1$ and $L^2$ limits are the same, i.e., $g=f$, and
$f_N \to f$ in $\widetilde{H}^1_L$.

Since we have already established part $(ii)$, we may use
(\ref{eq4.2}) to extend
$S_h$ to all of $H^1_L$ (thus in particular by (\ref{eq1.compare}) to $\widetilde{H}^1_L$) by continuity.
Let us momentarily call this extension
$\widetilde{S}_h$.  Then by (\ref{eq4.2}) and (\ref{eq1.compare}), we have
$$\|\widetilde{S}_h(f_N -
f)\|_{L^1(\RR^n)} \to 0.$$ Thus, by sublinearity, $S_hf_N \to
\widetilde{S}_hf$ in $L^1$. But also, $S_hf_N \to S_hf$ in $L^2$,
so $S_hf = \widetilde{S}_hf$ almost everywhere.  Therefore, using
(\ref{eq4.8}), we have that
$$\|f\|_{\widetilde{H}^1_L(\RR^n)} = \lim_{N\to \infty} \|f_N\|_{\widetilde{H}^1_L(\RR^n)} \leq C \sum
\lambda_k^j \leq \|S_h f\|_{L^1(\RR^n)} =
\|\widetilde{S}_hf\|_{L^1(\RR^n)} \leq C \|f\|_{H^1_L(\RR^n)}.$$

Next, we show that $m_k^j(N)$ converges weakly in each $L^p,
p_L<p<\widetilde{p}_L,$ and that the limits are molecules (up to a
harmless multiplicative constant).  Indeed, let $h \in L^{p'} \cap
L^2$. Then \begin{eqnarray*} \langle m_k^j(N),h\rangle &=&
\frac{1}{\lambda_k^j} \int_{\RR^n} \left(\int_{1/N}^N
(t^2Le^{-t^2L})^{M+1} \chi_{T_k^j} t^2L e^{-t^2L} f(x) \frac{dt}{t}
\right) \overline{h(x)} \, dx \\
&=& \frac{1}{\lambda_k^j}\int_{1/N}^N \langle \chi_{T_k^j} t^2L e^{-t^2L} f, (t^2L^*e^{-t^2L^*})^{M+1}h\rangle \frac{dt}{t}\\
&\to & \frac{1}{\lambda_k^j}\int_{0}^\infty \langle \chi_{T_k^j} t^2L e^{-t^2L} f, (t^2L^*e^{-t^2L^*})^{M+1}h\rangle \frac{dt}{t},
\end{eqnarray*}
by dominated convergence, since the square functions $$\left(\int_0^\infty |t^2Le^{-t^2L} f|^2 \frac{dt}{t} \right)^{1/2}, \,\,\, \left(\int_0^\infty |(t^2L^*e^{-t^2L^*})^{M+1} h|^2 \frac{dt}{t} \right)^{1/2},$$
belong, in particular, to $L^2$.  Similarly, but even more crudely, we may obtain existence of
$$\lim_{N\to \infty} \langle (\ell(Q_k^j)^2L)^{-i} m_k^j(N),h\rangle, \,\,\,\, 1 \leq i\leq M,$$
since $t\leq C \ell (Q_k^j)$ in $T_k^j$.  On the other hand, we have shown that, up to a multiplicative constant, the $m_k^j(N)$ are molecules, i.e., the bounds (\ref{eq1.8}) and (\ref{eq1.9})
hold uniformly in $N$, for $m_k^j(N)$, with $Q=Q_k^j$.  In particular,
$$\sup_{N} \|m_k^j(N)\|_{L^p(\RR^n)} \leq C |Q_k^j|^{1/p -1}.$$
Taking a supremum over $h \in L^{p'} (\RR^n)$, with norm $1$, we therefore obtain by the Riesz representation theorem that the weak limit $m_k^j$ belongs to $L^p$.  The desired bounds
(\ref{eq1.8}) and (\ref{eq1.9}) follow by taking $h \in L^{p'}(S_i(Q_k^j))$, and using the corresponding uniform bounds for $m_k^j(N).$  Thus, up to a multiplicative constant, the $m_k^j$ are molecules.

We now show that $f=\sum \lambda_k^j m_k^j $.  Let $\varphi \in C_0^\infty .$  Then, using absolute convergence, and the fact that $m_k^j(N) \to m_k^j$ weakly in $L^p$, we obtain
\begin{eqnarray*}
\int \varphi \left(\sum \lambda_k^j m_k^j\right)&=& \sum \lambda_k^j \int \varphi \, m_k^j\\
&=& \sum \lambda_k^j \lim_{N\to \infty}\int \varphi \,m_k^j(N)\\
&=& \lim_{N\to \infty}\sum \lambda_k^j \int \varphi \,m_k^j(N) \,\,\,\textrm{ (by dominated convergence)}\\
&=& \lim_{N\to \infty}\int \varphi \left(\sum \lambda_k^j m_k^j(N)\right)\\
&\equiv& \lim_{N\to \infty}\int \varphi f_N = \int \varphi f.
\end{eqnarray*}
Since this equality holds for all $\varphi \in C_0^\infty$, we have that $f = \sum \lambda_k^j m_k^j$
almost everywhere.

To complete the proof of the lemma, it remains to show that
$\{f_N\}$ is a Cauchy sequence in $\widetilde{H}^1_L$.  We recall that $f_N =
\sum \lambda_k^j m_k^j(N)$ where $\lambda_k^j = C 2^k |Q_k^j|,$
and
$$m_k^j(N) \equiv \frac{1}{\lambda_k^j} \int_{1/N}^N (t^2 L e^{-t^2 L} )^{M+1} \chi_{T_k^j} t^2 L e^{-t^2L} f \frac{dt}{t}.$$
For $K \in \NN ,$ we write \begin{eqnarray*}f_N &=& \sum_{j+k\leq K} \lambda_k^j m_k^j(N) \,+\,
\sum_{j+k > K} \lambda_k^j m_k^j(N)\\
&=& \sigma_K(N) + R_K(N).
\end{eqnarray*}
Then \begin{equation}\label{eq4.26} \sup_N
\|R_K(N)\|_{\widetilde{H}^1_L(\RR^n)} \leq \sum_{j+k > K}
|\lambda_k^j|\rightarrow 0\end{equation} as $K \to \infty .$ Thus,
it suffices to consider
$$\|\sigma_K(N) -\sigma_K(N')\|_{\widetilde{H}^1_L(\RR^n)} =
\left\|\sum_{j+k\leq K} \lambda_k^j \left(m_k^j(N) -
m_k^j(N')\right)\right\|_{\widetilde{H}^1_L(\RR^n)}.$$  Let $\eta
>0$ be given, and choose $K$ so that (\ref{eq4.26}) is bounded by
$\eta$. It is enough to show that for all $p\in
(p_L,\widetilde{p}_L)$, for every $\eps
>0$, $M>n/4$, and every $K\in \NN$, there exists an integer $N_1 =
N_1(\eta,K,p,\eps ,M)$ such that
\begin{equation}\label{eq4.27}\max_{j+k\leq K} \|m_k^j(N) - m_k^j(N')\|_{p,\eps, M,Q_k^j}
< \eta,\end{equation} whenever $N' \geq N \geq N_1$, where the
``$(p,\eps ,M)$-molecular norm adapted to $Q$" is defined as
$$\|\mu\|_{p,\eps ,M,Q} \equiv  \sup_{i\geq 0} 2^{i(n -
n/p+\eps)}|Q|^{1-1/p}\sum_{\nu=0}^M\|(\ell(Q)^2L)^{-\nu}\mu\|_{L^p(S_i(Q))}.$$
To this end, we note that, for $N,N'$ sufficiently large,
$$\mu_k^j(N,N') \equiv m_k^j(N') - m_k^j(N) = \frac{1}{\lambda_k^j} \int_{1/N'}^{1/N} (t^2 L e^{-t^2 L} )^{M+1} \chi_{T_k^j} t^2 L e^{-t^2L} f \frac{dt}{t},$$
since $t \leq C \ell(Q_k^j)$ in $T_j^k.$
Let $h\in L^{p'}(S_i(Q_k^j)),$ with $\|h\|_{p'} = 1.$  Then, following the argument from
(\ref{eq4.9}) to (\ref{eq4.17}), we obtain that
\begin{eqnarray*}|\langle \mu_k^j(N,N'),h\rangle | &=& \left|\frac{1}{\lambda_k^j} \int_{1/N'}^{1/N} \langle \chi_{T_k^j} t^2 L e^{-t^2L} f,(t^2 L^* e^{-t^2 L^*} )^{M+1} h \rangle\frac{dt}{t}\right|\\&\leq&
C \frac{1}{\lambda_k^j} \left(\int_{CQ_k^j \cap \,
^cO_{k+1}}(S_h^{1/N} f)^p dx \right)^{1/p} 2^{-i(n-n/p + \eps
)},\end{eqnarray*} where $$S_h^{1/N} f \equiv \left( \int
\!\!\int_{|x-y|<t<1/N} |t^2Le^{-t^2L} f(y)|^2 \frac{dy
dt}{t^{n+1}}\right)^{1/2}.$$ Now, since $f \in L^2$, $S_h^{1/N}f \to
0$ in $L^2$. We choose $N_1$ so large that
\begin{equation}\label{eq4.28}\|S_h^{1/N} f \|_{L^2(\RR^n)} \leq
\frac{1}{C}\, \eta^R \min_{j+k \leq K} \lambda_k^j
\,|Q_k^j|^{-1/2},\end{equation} whenever $N \geq N_1,$ where $R$
will be chosen depending on $p$. If $p=2$, taking $R=1$, and
taking the surpemum over $h$ as above, we obtain immediately that
$$\|\mu_k^j(N,N')\|_{L^2(S_i(Q_k^j))} \leq  \eta 2^{-i(n/2 + \eps )} |Q_k^j|^{-1/2}. $$
The case $p_L<p<2$ follows by H\"{o}lder's inequality.  For $2<p<\widetilde{p}_L,$ we choose
$r\in (p,\widetilde{p}_L),$ and using that $S_hf \leq 2^{k+1}$ on $^cO_{k+1}$ by definition, we interpolate between
(\ref{eq4.28}) and the crude bound
$$\left(\int_{CQ_k^j \cap \, ^cO_{k+1}}(S_h^{1/N} f(x))^r dx \right)^{1/r} \leq C 2^k |Q_k^j|^{1/r},$$
to deduce that
$$|\langle\mu_k^j(N,N'),h\rangle | \leq  \eta 2^{-i(n - n/p + \eps )} |Q_k^j|^{1/p-1}, $$
for $R$ chosen large enough depending on $p$.
We now obtain (\ref{eq4.27}), by applying $(\ell(Q_k^j)^2L)^{-\nu}$ to $\mu_k^j(N,N')$, and then
repeating the previous argument with
minor changes.   It follows that $\{f_N\}$ is a Cauchy sequence in $\widetilde{H}^1_L$.
This concludes the proof of Lemma~\ref{lmain} and therefore also that of Theorem~\ref{t4.1}.
 \ep

\vskip 0.08 in

We conclude this section with

\begin{corollary}\label{c4.3} The spaces $H^1_L(\RR^n)$ coincide for
different choices of $\eps>0$, $p_L<p<\widetilde{p}_L$ and $M\in\NN$
such that $M>n/4$.
\end{corollary}

\section{Characterization by the square function associated to the Poisson semigroup}
\setcounter{equation}{0}

We start with the following auxiliary result.

\begin{lemma}\label{l5.1} Fix $K\in\NN$.
For all closed sets $E$, $F$ in
$\RR^n$ with ${\rm dist }(E,F)>0$

\begin{equation}\label{eq5.1}
\left\|(t\sqrt L)^{2K}e^{-t\sqrt L}\right\|_{L^2(F)}\leq
C\,\left(\frac{t}{{\rm
dist}\,(E,F)}\right)^{2K+1}\,\|f\|_{L^2(E)},\qquad\forall\,t>0,\\[4pt]
\end{equation}

\noindent if $f\in L^2(\RR^n)$ is supported in $E$.
\end{lemma}

\bp The subordination formula

\begin{equation}\label{eq5.2} e^{-t\sqrt
L}f=C\int_0^\infty\frac{e^{-u}}{\sqrt u}\,\,
e^{-\frac{t^2L}{4u}}f\,du
\end{equation}

\noindent allows us to write

\begin{eqnarray}\nonumber
\left\|(t\sqrt L)^{2K}e^{-t\sqrt L}\right\|_{L^2(F)} &\leq & C
\int_0^\infty \frac{e^{-u}}{\sqrt
u}\,\left\|\Bigl(\frac{t^2L}{4u}\Bigr)^K
e^{-\frac{t^2L}{4u}}f\right\|_{L^2(F)}
u^K\,du\\[4pt]
& \leq & C \|f\|_{L^2(E)}\int_0^\infty e^{-u}\,e^{-\frac{{\rm
dist}\,(E,F)^2}{cu^2}}\,u^{K-1/2}\,du.\label{eq5.3}
\end{eqnarray}

\noindent Then we make the change of variables $u\mapsto
s:=u\frac{{\rm dist}\,(E,F)^2}{t^2}$ to bound (\ref{eq5.3}) by

\begin{eqnarray}\nonumber
&&C \|f\|_{L^2(E)}\int_0^\infty e^{-s\frac{t^2}{{\rm
dist}\,(E,F)^2}}\,e^{-s}\,\Bigl(s\frac{t^2}{{\rm
dist}\,(E,F)^2}\Bigr)^{K-1/2}\frac{t^2}{{\rm
dist}\,(E,F)^2}\,\,ds\\[4pt]\nonumber
&&\qquad \leq C\|f\|_{L^2(E)}\left(\frac{t}{{\rm
dist}\,(E,F)}\right)^{2K+1}\int_0^\infty e^{-s}\,s^{K-1/2}\,ds\\[4pt]
&&\qquad \leq C \left(\frac{t}{{\rm
dist}\,(E,F)}\right)^{2K+1}\|f\|_{L^2(E)},\label{eq5.4}
\end{eqnarray}

\noindent as desired. \ep

\begin{theorem}\label{t5.2} Consider the operator

\begin{equation}\label{eq5.5}
S_P^Kf(x):=\left(\dint_{\Gamma(x)}|(t\sqrt L)^{2K}e^{-t\sqrt
L}f(y)|^2\,\frac {dydt}{t^{n+1}}\right)^{1/2}, \quad
x\in\RR^n,\quad f\in L^2(\RR^n).
\end{equation}
Suppose $K\in\NN$, $M\in\NN$, $M+K>n/4-1/2$, and
$\eps=2M+2K+1-n/2$. If $f\in L^2(\RR^n)$, with
$\|S_P^Kf\|_{L^1(\RR^n)}<\infty$, then $f \in H^1_L.$  Furthermore, there exists a family of
$(2,\eps,M)$-molecules $\{m_i\}_{i=0}^\infty$ and a sequence of
numbers $\{\lambda_i\}_{i=0}^\infty$ such that $f$ can be
represented in the form $f=\sum_{i=0}^\infty \lambda_i m_i$, with

\begin{equation}\label{eq5.6}
\|f\|_{H^1_L(\RR^n)} \leq C \sum_{i=0}^\infty |\lambda_i|\leq
C\|S_P^Kf\|_{L^1(\RR^n)}.
\end{equation}

\end{theorem}

\bp The lemma can be proved following the argument of
Theorem~\ref{t4.1} with minor modifications. To be more precise,
we use the Calder\'{o}n reproducing formula in the form

\begin{equation}
f=C\int_0^\infty \Bigl((t^2L)^{M+K}e^{-t\sqrt
L}\Bigr)^2f\,\frac{dt}{t}=C\int_0^\infty (t^2L)^{2M+K}e^{-t\sqrt
L}(t^2L)^{K}e^{-t\sqrt L}f\,\frac{dt}{t}, \label{eq5.7}
\end{equation}

\noindent for $f\in L^2(\RR^n)$. To be completely rigorous, we should truncate and approximate by
$f_N$ as in the proof of Theorem~\ref{t4.1}.  As the details are similar in the present case,
we shall merely sketch a formal proof, and leave the details of the limiting arguments to the reader.

To begin, we define

\begin{equation}\label{eq5.8}
{O}_k:=\{x\in\RR^n:\,S_P^Kf(x)>2^k\},
\end{equation}

\noindent and

\begin{equation}\label{eq5.9}
m_k^j=\frac{1}{\lambda_k^j}\int_0^\infty (t^2L)^{2M+K}e^{-t\sqrt
L}\left(\chi_{T_k^j}(t^2L)^{K}e^{-t\sqrt L}\right)f\,\frac{dt}{t},
\end{equation}

\noindent with $T_k^j$, $k,j\in\ZZ$, given analogously to
(\ref{eq4.5}). The rest of the proof follows the same path, using
Lemma~\ref{l5.1} instead of Gaffney relations, which allows to
derive the estimate

\begin{equation}\label{eq5.10}
\|m_k^j\|_{L^2(S_i(Q_k^j))}\leq C\,2^{-i(4M+2K+1)}|Q_k^j|^{-1/2},
\qquad i=0,1,2,...,
\end{equation}

\noindent for all $k,j\in\ZZ$.

As for the vanishing moment condition,

\begin{equation}\label{eq5.11}
(l(Q_k^j)^{-2}L^{-1})^{M}m_k^j=\frac{1}{\lambda_k^j}\int_0^\infty
\Bigl(\frac{t}{l(Q_k^j)}\Bigr)^{2M} (t^2L)^{M+K}e^{-t\sqrt
L}\left(\chi_{T_k^j}(t^2L)^{K}e^{-t\sqrt L}\right)f\,\frac{dt}{t},
\end{equation}

\noindent and hence

\begin{equation}\label{eq5.12}
\|(l(Q_k^j)^{-2}L^{-1})^{M}m_k^j\|_{L^2(S_i(Q_k^j))}\leq
C\,2^{-i(2M+2K+1)}|Q_k^j|^{-1/2}, \qquad i=0,1,2,....
\end{equation}

\noindent Combined with (\ref{eq5.10}), this finishes the argument.
\ep

\begin{theorem}\label{t5.3}
Let $ \eps>0$ and $M>n/4$. Then for every representation
$\sum_{i=0}^\infty \lambda_i m_i$, where $\{m_i\}_{i=0}^\infty$ is a family of
$(2,\eps,M)$-molecules  and
$\sum_{i=0}^\infty |\lambda_i|<\infty$, the series
$\sum_{i=0}^\infty \lambda_i m_i$ converges in $H^1_{S_P}(\RR^n)$
and

\begin{equation}\label{eq5.13}
\Bigl\|\sum_{i=0}^\infty \lambda_i m_i\Bigr\|_{H^1_{S_P}(\RR^n)}\leq
C\sum_{i=0}^\infty |\lambda_i|.
\end{equation}
\end{theorem}

\bp We will follow the argument of Theorem~\ref{t4.1}, Step II, and
mention only necessary changes.

First, by Lemma \ref{l3.sub}, it will be enough to establish a
uniform $L^1$ bound on molecules.  To this end, we observe that the operator

\begin{equation}\label{eq5.14}
g_Pf:=\left(\int_0^\infty |t\nabla e^{-t \sqrt L}f|^2\,\frac
{dt}{t}\right)^{1/2},
\end{equation}

\noindent is bounded in $L^2(\RR^n)$. This follows from the
estimates on the operators having bounded holomorphic functional
calculus in $L^2$ (see \cite{ADM}) and integration by parts. Then
$S_P$ is bounded in $L^2(\RR^n)$, since

\begin{eqnarray}\nonumber
&&\hskip -.7cm \|S_Pf\|_{L^2(\RR^n)}\leq C
\left(\dint_{\RR^n\times(0,\infty)}
\left[\int_{\RR^n}\chi_{\{x:\,|x-y|< t\}}(x) \,dx\right] |t\nabla
e^{-t \sqrt L}f(y)|^2\,\frac {dydt}{t^{n+1}}\right)^{1/2}
\\[4pt]\label{eq5.15}
&&\leq C \left(\dint_{\RR^n\times(0,\infty)} |t\nabla e^{-t \sqrt
L}f(y)|^2\,\frac {dydt}{t}\right)^{1/2} \leq
C\|g_Pf\|_{L^2(\RR^n)}\leq C\|f\|_{L^2(\RR^n)}.
\end{eqnarray}

Therefore, for $j=0,1,2$,

\begin{equation}\label{eq5.16}
\|S_P(I-e^{-l(Q)^2L})^{M}(m\chi_{S_i(Q)})\|_{L^2(S_j(Q_i))}\leq C
\|m\|_{L^2(S_i(Q))}.
\end{equation}

\noindent Turning to the case $j\geq 3$, we write

\begin{equation}\label{eq5.17}
\|S_P(I-e^{-l(Q)^2L})^{M}(m\chi_{S_i(Q)})\|_{L^2(S_j(Q_i))}^2\leq
I+\sum_{k=0}^{j-2}I_k,
\end{equation}

\noindent where

\begin{equation}\label{eq5.18}
I=C \int_{\RR^n\setminus Q_{j-2+i}}\int_0^\infty |t\nabla e^{-t
\sqrt L}(I-e^{-l(Q)^2L})^{M}(m\chi_{S_i(Q)})(x)|^2\,\frac{dtdx}{t}
\end{equation}

\noindent and

\begin{equation}\label{eq5.19}
I_k=C\int_{S_k(Q_i)}\int_{(2^{j-1}-2^k)2^il(Q)}^\infty |t\nabla
e^{-t \sqrt
L}(I-e^{-l(Q)^2L})^{M}(m\chi_{S_i(Q)})(x)|^2\,\frac{dtdx}{t}
\end{equation}

\noindent for $k=0,...,j-2$. Then

\begin{eqnarray}\nonumber
&&\hskip -0.7cm I^{1/2} \leq  C \int_0^\infty e^{-u}
\left(\int_{\,^c Q_{j-2+i}}\int_0^\infty
\left|\frac{t}{\sqrt{4u}}\nabla e^{-\frac{t^2 L}{4u}}
(I-e^{-l(Q)^2L})^{2M}(m\chi_{S_i(Q)})(x)\right|^2\,\frac{dtdx}{t}\right)^{1/2}du\\[4pt]\nonumber
&&\hskip -0.7cm\qquad\leq  C \int_0^\infty e^{-u} \left(\int_{\,^c
Q_{j-2+i}}\int_0^\infty \left|s\nabla e^{-s^2L}
(I-e^{-l(Q)^2L})^{4M}(m\chi_{S_i(Q)})(x)\right|^2\,\frac{dsdx}{s}\right)^{1/2}du,
\end{eqnarray}

\noindent where we made the change of variables $t\mapsto
s:=\frac{t}{\sqrt{4u}}$. However, it can be proved along the lines
of Theorem~\ref{t3.2} that $\tilde g_hf:=\left(\int_0^\infty
|t\nabla e^{-t^2L}f|^2\,\frac {dt}{t}\right)^{1/2}$ (similarly to
$g_hf$) satisfies the estimates (\ref{eq3.2}) with $p=2$ and
therefore

\begin{eqnarray}\nonumber
I &\leq & C\,\left(\frac{l(Q)^2}{{\rm dist}\,(S_i(Q),\RR^n\setminus
Q_{j-2+i})^2}\right)^{2M}\,\|m\|_{L^2(S_i(Q))}^2\leq
C\,\Bigl(\frac{1}{2^{i+j}}\Bigr)^{4M}\,\|m\|_{L^2(S_i(Q))}^2.
\end{eqnarray}

Concerning $I_k$, $k=0,1,...,j-2$, we use subordination formula once
again to write

\begin{eqnarray}\nonumber
I_k^{1/2}&\leq & C\int_0^\infty
e^{-u}\Bigl(\int_{S_k(Q_i)}\int_{(2^{j-1}-2^k)2^il(Q)}^\infty
\Bigl|\frac{t}{\sqrt{4u}}\nabla e^{-\frac{t^2
L}{4u}}\\[4pt]
&&\qquad\qquad\qquad\times(I-e^{-l(Q)^2L})^{M}
(m\chi_{S_i(Q)})(x)\Bigr|^2\,\frac{dtdx}{t}\Bigr)^{1/2}du.\label{eq5.20}
\end{eqnarray}

\noindent  Then one can make a change of variables
$s:=\frac{t^2}{4u(m+1)}$, so that following
(\ref{eq4.22})--(\ref{eq4.23})

\begin{eqnarray}
I_k^{1/2}&\leq& C \int_0^\infty e^{-u}\Bigl(
\int_{\frac{[(2^{j-1}-2^k)2^il(Q)]^2}{cu}}^\infty \|\sqrt
s\,\nabla e^{-(m+1)sL}\nonumber\\[4pt]
&& \qquad\qquad\qquad\times (I-e^{-l(Q)^2L})^{M}
(m\chi_{S_i(Q)})\|^2_{L^2(S_k(Q_i))}\,\frac{ds}{s}\Bigr)^{1/2}du\nonumber\\[4pt]
&\leq& C \int_0^\infty
e^{-u}\left(\int_{\frac{[(2^{j-1}-2^k)2^il(Q)]^2}{cu}}^\infty
 e^{-\frac{{\rm
dist}\,(S_i(Q),S_k(Q_i))^2}{cs(1+u)}}\,
\Bigl(\frac{l(Q)^2}{s}\Bigr)^{2M}\,\frac
{ds}{s}\right)^{1/2}du\nonumber\\[4pt]
 && \qquad\qquad\qquad\times \|m\|_{L^2(S_i(Q))}, \label{eq5.21}
\end{eqnarray}

\noindent and hence

\begin{equation}
I_k \leq C \Bigl(\frac{1}{2^{i+j}}\Bigr)^{4M}
\|m\|_{L^2(S_i(Q))}^2. \label{eq5.22}
\end{equation}

\noindent  Then

\begin{eqnarray}\nonumber
&&\hskip -0.7cm\|S_P(I-e^{-l(Q)^2L})^{M}(m\chi_{S_i(Q)})\|_{L^1(\RR^n)}\\[4pt]
\nonumber &&\leq C\sum_{j=3}^\infty(2^{i+j}l(Q))^{n/2}\,\sqrt j\,
\left(\frac {1}{2^{i+j}}\right)^{2M}
\|m\|_{L^2(S_i(Q))}+C(2^{i}l(Q))^{n/2}\|m\|_{L^2(S_i(Q))}\leq
C\,2^{-i\eps}.
\end{eqnarray}

Similar argument provides an estimate for
$\|S_P[I-(I-e^{-l(Q)^2L})]^{M}m\|_{L^1(\RR^n)}$ and finishes the
proof.  \ep

\begin{lemma}\label{l5.4}
For all $f\in L^2(\RR^n)$

\begin{equation}\label{eq5.23}
\|S_P^1f\|_{L^1(\RR^n)} \leq C \|S_Pf\|_{L^1(\RR^n)}.
\end{equation}
\end{lemma}

\bp To start, let us define the family of truncated cones

\begin{equation}\label{eq5.24}
\Gamma^{\eps,R,\alpha}(x):=\{(y,t)\in\RR^n\times(\eps,
R):\,|x-y|<t\alpha\}, \qquad x\in\RR^n.
\end{equation}

\noindent Then for every function $\eta\in C_0^\infty
(\Gamma^{\eps/2,2R,2}(x))$ such that $\eta\equiv 1$ on
$\Gamma^{\eps,R,1}(x)$ and $0\leq\eta\leq 1$

\begin{eqnarray}
&&\left(\dint_{\Gamma^{\eps,R,1}(x)}|t^2Le^{-t\sqrt
L}f(y)|^2\,\frac {dydt}{t^{n+1}}\right)^{1/2} \nonumber\\[4pt]
&&\qquad \leq \left(\dint_{\Gamma^{\eps/2,2R,2}(x)}t^2Le^{-t\sqrt
L}f(y)\,{\overline{t^2Le^{-t\sqrt L}f(y)}}\eta(y,t)\,\frac
{dydt}{t^{n+1}}\right)^{1/2}\nonumber\\[4pt]
&&\qquad \leq \left(\dint_{\Gamma^{\eps/2,2R,2}(x)}tA\nabla
e^{-t\sqrt L}f(y)\cdot t \nabla\left[{\overline{t^2Le^{-t\sqrt
L}f(y)}}\right]\eta(y,t)\,\frac
{dydt}{t^{n+1}}\right)^{1/2}\nonumber\\[4pt]
&&\qquad \quad+ \left(\dint_{\Gamma^{\eps/2,2R,2}(x)}t A\nabla
e^{-t\sqrt L}f(y)\cdot{\overline{t^2Le^{-t\sqrt
L}f(y)}}\,t\nabla\eta(y,t)\,\frac
{dydt}{t^{n+1}}\right)^{1/2}.\label{eq5.25}
\end{eqnarray}

\noindent We can always assume that $\|\nabla
\eta\|_{L^\infty(\Gamma^{\eps/2,2R,2}(x))}\leq 1/t$, so that the
expression above is bounded by

\begin{eqnarray}
&&\hskip -0.7 cm \left(\dint_{\Gamma^{\eps/2,2R,2}(x)}|t\nabla
e^{-t\sqrt L}f(y)|^2 \frac{dydt}{t^{n+1}}\right)^{1/4}
\left(\dint_{\Gamma^{\eps/2,2R,2}(x)}|t \nabla t^2Le^{-t\sqrt
L}f(y)|^2\,\frac
{dydt}{t^{n+1}}\right)^{1/4}\nonumber\\[4pt]
&&\hskip -0.7 cm \qquad + \left(\dint_{\Gamma^{\eps/2,2R,2}(x)}|t
\nabla e^{-t\sqrt L}f(y)|^2\,\frac
{dydt}{t^{n+1}}\right)^{1/4}\left(\dint_{\Gamma^{\eps/2,2R,2}(x)}
|t^2Le^{-t\sqrt L}f(y)|^2\,\frac
{dydt}{t^{n+1}}\right)^{1/4}.\nonumber
\end{eqnarray}

\noindent Consider now the decomposition of the set
$\Gamma^{\eps/2,2R,2}(x)$ given by collection of balls \break
$\{B(z_k,r_k)\}_{k=0}^\infty$ in $\RR^{n+1}$, such that

\begin{eqnarray}\nonumber
&& \qquad \cup_{k=0}^\infty
B(z_k,r_k)=\cup_{k=0}^\infty
B(z_k,2r_k)=\Gamma^{\eps/2,2R,2}(x),\\[4pt]
&& c_1 \,{\rm dist}(z_k, \,^c(\Gamma^{\eps/2,2R,2}(x))) \leq r_k
\leq c_2\, {\rm dist}(z_k,  \,^c(\Gamma^{\eps/2,2R,2}(x))),
\label{eq5.26}
\end{eqnarray}

\noindent every point $z\in  \Gamma^{\eps/2,2R,2}(x)$ belongs at
most to $c_3$ balls of the form $B(z_k,2r_k)$, $0<c_1<c_2<1$ and
$c_3\in\NN$ are some fixed constants. Such collection
$\{B(z_k,r_k)\}_{k=0}^\infty$ can always be built on the basis of
Whitney decomposition (for the latter see \cite{CW}, \cite{St}).
Then we use Caccioppoli inequality (Lemma~\ref{l2.7}) for the
operator $\tilde Lf=-{\rm{div}}_{y,t}(B\nabla_{y,t} f)$ in place of
$L$. Here $\tilde L$ is understood in the usual weak sense, $B$ is
the $(n+1)\times(n+1)$ block diagonal matrix with components $1$ and
$A$ and ${\rm{div}}_{y,t}$, $\nabla_{y,t}$ denote, respectively,
divergence and gradient taken in space and time variables. Clearly,
$\tilde L e^{-t\sqrt L}f=0$. We obtain

\begin{eqnarray}\nonumber
&&\hskip -1.2cm \dint_{\Gamma^{\eps/2,2R,2}(x)}|t \nabla t^2Le^{-t\sqrt
L}f(y)|^2\,\frac
{dydt}{t^{n+1}}\leq \sum_{k=0}^\infty
\dint_{B(z_k,r_k)}|t \nabla_{y,t} t^2Le^{-t\sqrt
L}f(y)|^2\,\frac {dydt}{t^{n+1}}\\[4pt]
&&\hskip -1.2cm\leq C\sum_{k=0}^\infty
\dint_{B(z_k,2r_k)}\left|\frac{t}{r_k} \,t^2Le^{-t\sqrt
L}f(y)\right|^2\,\frac {dydt}{t^{n+1}}\leq
C\dint_{\Gamma^{\eps/2,2R,2}(x)} |t^2Le^{-t\sqrt L}f(y)|^2\,\frac
{dydt}{t^{n+1}}.\label{eq5.27}
\end{eqnarray}

\noindent Combining this with the formulae above and passing to the limit as
$\eps\to 0$ and $R\to\infty$, we arrive at

\begin{equation}\label{eq5.28}
S_P^1 f(x)\leq C (S_Pf(x))^{1/2}(S_P^1f(x))^{1/2},
\end{equation}

\noindent and hence (\ref{eq5.23}), as desired. \ep

\vskip 0.08in

\begin{corollary}\label{c5.5}
$H^1_L(\RR^n)=H^1_{S_P}(\RR^n)$, in particular,
$\|f\|_{H^1_L(\RR^n)}\approx\|S_P f\|_{L^1(\RR^n)} + \|f\|_{L^1(\RR^n)}$.
\end{corollary}

\bp  The left-to-the-right inclusion follows from
Theorem~\ref{t5.3}, the converse from Theorem~\ref{t5.2} combined
with Lemma~\ref{l5.4}; in addition, we use
Corollary~\ref{c4.3} to guarantee that the molecular spaces $H^1_L(\RR^n)$
coincide for different choices of $p\in (p_L,\widetilde{p}_L)$ and
$\eps>0$, thus removing the constraints on $\eps$ and $p$ in Theorem~\ref{t5.2}.
We omit the details.\ep

Finally, consider two more versions of the square function:

\begin{eqnarray}\label{eq5.29}
\bar S_Pf(x)&:=& \left(\dint_{\Gamma(x)}|t\sqrt L e^{-t
\sqrt L}f(y)|^2\,\frac {dydt}{t^{n+1}}\right)^{1/2},\\[4pt]\label{eq5.30}
\widehat S_P f(x)&:=& \left(\dint_{\Gamma(x)}|t\nabla_{x,t} e^{-t
\sqrt L}f(y)|^2\,\frac {dydt}{t^{n+1}}\right)^{1/2},
\end{eqnarray}

\noindent where $f\in L^2(\RR^n)$, $x\in\RR^n$ and $\nabla_{y,t}$
stands for the gradient in space and time variables.

\begin{theorem}\label{t5.6}
We have the equivalence

\begin{equation}\label{eq5.31}
\|f\|_{H^1_L(\RR^n)}\approx\|\bar S_P f\|_{L^1(\RR^n)}\approx
\|\widehat S_P f\|_{L^1(\RR^n)}.
\end{equation}
\end{theorem}

This result is just a slight modification of the previous ones in
this section. In fact, the argument for $\bar S_P $ follows the
same lines as the proofs of Theorems~\ref{t5.2}--\ref{t5.3} ones
we observe that $t\sqrt L e^{-t\sqrt L}=-t\partial_t e^{-t\sqrt
L}$, and the result for $\widehat S_P$ is a combination of those
for $S_P$ and $\bar S_P$.

\section{Characterization by the non-tangential maximal function associated to the heat semigroup}
\setcounter{equation}{0}

\begin{theorem}\label{t6.1}
For every $f\in L^2(\RR^n)$

\begin{equation}\label{eq6.1}
\|S_hf\|_{L^1(\RR^n)}\leq C\|{\cal N}_hf\|_{L^1(\RR^n)}.
\end{equation}
\end{theorem}

\bp The idea of the proof is based on the analogous argument for
the maximal and square functions associated to the Poisson
semigroup for Laplacian that appeared in \cite{FeSt}. Similar
ideas have also been used in \cite{AR}.

To begin, notice that the argument of Lemma~\ref{l5.4} also
provides the estimate

\begin{equation}\label{eq6.2}
\|S_hf\|_{L^1(\RR^n)}\leq C\|\tilde S_hf\|_{L^1(\RR^n)},
\end{equation}

\noindent for $\widetilde S_h=\widetilde S_h^1$, where

\begin{equation}\label{eq6.3}
\widetilde S_h^\beta f(x):= \left(\dint_{\Gamma^\beta(x)}|t\nabla
e^{-t^2 L}f(y)|^2\,\frac {dydt}{t^{n+1}}\right)^{1/2},\quad f\in
L^2(\RR^n), \quad x\in\RR^n.
\end{equation}

\noindent Therefore, it is enough to prove (\ref{eq6.1}) with
$\widetilde S_h$ in place of $S_h$. Also, recall the definition of
truncated cone (\ref{eq5.24}) and denote

\begin{equation}\label{eq6.4}
\widetilde S_h^{\eps, R, \beta} f(x):= \left(\dint_{\Gamma^{\eps, R,
\beta}(x)}|t\nabla e^{-t^2 L}f(y)|^2\,\frac
{dydt}{t^{n+1}}\right)^{1/2},\quad f\in L^2(\RR^n), \quad x\in\RR^n.
\end{equation}

In what follows we will work with $\widetilde S_h^{\eps, R, \beta}$
rather than $\widetilde S_h^{\beta}$ and then pass to the limit as
$\eps\to 0$, $R\to\infty$, all constants in estimates will not
depend on $\eps$ and $R$ unless explicitly stated.

Consider the non-tangential maximal function

\begin{equation}\label{eq6.5}
{\cal N}_h^\beta
f(x):=\sup_{(y,t)\in\Gamma^\beta(x)}\left(\frac{1}{(\beta
t)^n}\int_{B(y,\beta t)} |e^{-t^2L}f(z)|^2\,dz\right)^{1/2},\quad
 f\in L^2(\RR^n),
\end{equation}

\noindent  where $\Gamma^\beta(x)$, $x\in\RR^n$, $\beta>0$, is the
cone of aperture $\beta$. Let us introduce the following sets:

\begin{equation}\label{eq6.6}
E:=\{x\in\RR^n:\,{\cal N}_h^\beta f(x)\leq \sigma\}, \quad
\sigma\in\RR,
\end{equation}

\noindent  where $\beta$ is some fixed constant to be determined
later, and

\begin{equation}\label{eq6.7}
E^*:=\left\{x\in\RR^n:\, \mbox{for every} \,B(x), \mbox{ ball in
$\RR^n$ centered at $x$},\quad  \frac{|E\cap B(x)|}{|B(x)|}\geq
\frac 12 \right\},
\end{equation}

\noindent the set of points having global $1/2$ density with
respect to $E$. Also,

\begin{equation}\label{eq6.8}
B:=\,^cE,\qquad B^*:=\,^cE^*.
\end{equation}

\noindent Finally, denote

\begin{eqnarray}\label{eq6.9}
&&{\cal R}^{\eps, R, \beta}(E^*):=\bigcup_{x\in E^*}\Gamma^{\eps,
R,\beta}(x),\\[4pt]\label{eq6.10}
&& \hskip -1cm {\cal B}^{\eps, R, \beta}(E^*) \mbox{\, -- the
boundary of \,}{\cal R}^{\eps, R, \beta}(E^*),
\end{eqnarray}

\noindent and

\begin{equation}\label{eq6.11}
u(y,t):=e^{-t^2L}f(y), \qquad t\in(0,\infty), \qquad y\in\RR^n.
\end{equation}

It is not hard to see that

\begin{eqnarray}\nonumber
&&\int_{E^*}\left(\widetilde S_h^{2\eps, R,
1/2}f(x)\right)^{2}\,dx \leq  \int_{E^*}\left(\widetilde
S_h^{\alpha\eps, \alpha R, 1/\alpha}f(x)\right)^{2}\,dx\\[4pt]
&&\qquad\quad \leq  C \dint_{{\cal R}^{\alpha\eps, \alpha R,
1/\alpha}(E^*)} t|\nabla u(y,t)|^2\,dydt, \qquad \forall\,\alpha\in
(1,2), \label{eq6.12}
\end{eqnarray}

\noindent by Lemma~\ref{l2.1}. Going further,

\begin{eqnarray}\nonumber
&&\hskip -2cm \dint_{{\cal R}^{\alpha\eps, \alpha R,
1/\alpha}(E^*)} t|\nabla u(y,t)|^2\,dydt  = \dint_{{\cal
R}^{\alpha\eps, \alpha R, 1/\alpha}(E^*)} t\nabla u(y,t)\cdot
{\overline{\nabla
u(y,t)}}\,dydt\\[4pt]\label{eq6.13}
&&\quad \hskip -2cm \leq C \Re e \dint_{{\cal R}^{\alpha\eps,
\alpha R, 1/\alpha}(E^*)} \left[tA(y)\nabla u(y,t)\cdot
{\overline{\nabla u(y,t)}}+\nabla u(y,t)\cdot
t{\overline{A(y)\nabla u(y,t)}}\right]\,dydt,
\end{eqnarray}

\noindent using ellipticity of $A$. Now we integrate by parts to
bound (\ref{eq6.13}) by

\begin{eqnarray}\nonumber
&&\hskip -.7cm C \Re e \dint_{{\cal R}^{\alpha\eps, \alpha R,
1/\alpha}(E^*)} \left[-\,t\,{\rm div}A(y)\nabla u(y,t)
{\overline{u(y,t)}}-u(y,t)\overline{t\,{{\rm div}A(y)\nabla
u(y,t)}}\right]\,dydt+\\[4pt]\nonumber
&&\hskip -.7cm \quad +C\Re e \int_{{\cal B}^{\alpha\eps, \alpha R,
1/\alpha}(E^*)} \left[tA(y)\nabla u(y,t)\cdot
\nu_y(y,t){\overline{u(y,t)}}+u(y,t)\nu_y(y,t)\cdot
t{\overline{A(y)\nabla u(y,t)}}\right]\,d\sigma_{y,t},
\end{eqnarray}

\noindent where $\nu_y(y,t)$, $y\in\RR^n$, $t\in\RR$, is the
projection of normal vector to ${\cal B}^{\alpha\eps, \alpha R,
1/\alpha}(E^*)$ on $\RR^n$ (similarly, $\nu_t$ will denote
projection on $\RR$). However, $u$ given by (\ref{eq6.11}) is a
solution of system $\partial_t u=-2t\,{\rm div} A \nabla u$, and
hence the first integral above can be represented modulo
multiplicative constant as

\begin{eqnarray}\nonumber
&&\dint_{{\cal R}^{\alpha\eps, \alpha R, 1/\alpha}(E^*)}
\left[\partial_t u(y,t)\cdot {\overline{u(y,t)}}+u(y,t)\cdot
{\overline{\partial_t u(y,t)}}\right]\,dydt\\[4pt]
&&\qquad =\dint_{{\cal R}^{\alpha\eps, \alpha R, 1/\alpha}(E^*)}
\partial_t |u(y,t)|^2\,dydt = \int_{{\cal B}^{\alpha\eps, \alpha R, 1/\alpha}(E^*)}
|u(y,t)|^2\nu_t(y,t)\,d\sigma_{y,t}.\label{eq6.14}
\end{eqnarray}

\noindent Combining (\ref{eq6.13})--(\ref{eq6.14}), one can write

\begin{eqnarray}\nonumber
&&\hskip -0.7cm \int_1^2\dint_{{\cal R}^{\alpha\eps, \alpha R,
1/\alpha}(E^*)} t|\nabla
u(y,t)|^2\,dydt\,d\alpha\\[4pt]\nonumber
&&\hskip -0.7cm \quad \leq C\int_1^2\int_{{\cal B}^{\alpha\eps,
\alpha R, 1/\alpha}(E^*)} t|\nabla
u(y,t)|\,|u(y,t)|d\sigma_{y,t}\,d\alpha+ C\int_1^2\int_{{\cal
B}^{\alpha\eps, \alpha R, 1/\alpha}(E^*)}
|u(y,t)|^2\,d\sigma_{y,t}\,d\alpha\\[4pt]\nonumber
&&\hskip -0.7cm \quad \leq C\dint_{\widetilde {\cal B}^{\eps,
R}(E^*)} |\nabla u(y,t)|\,|u(y,t)|\,dydt+ C\dint_{\widetilde {\cal
B}^{\eps, R}(E^*)}|u(y,t)|^2\,\frac{dydt}{t}\\[4pt]\nonumber
&&\hskip -0.7cm \quad \leq C\left(\dint_{\widetilde {\cal
B}^{\eps, R}(E^*)} t|\nabla
u(y,t)|^2\,dydt\right)^{1/2}\,\left(\dint_{\widetilde {\cal
B}^{\eps, R}(E^*)}|u(y,t)|^2\,\frac{dydt}{t}\right)^{1/2}+\\[4pt]
&&\hskip 1.7cm +\,C\dint_{\widetilde {\cal B}^{\eps,
R}(E^*)}|u(y,t)|^2\,\frac{dydt}{t}\label{eq6.15}
\end{eqnarray}

\noindent where

\begin{equation}\label{eq6.16}
\widetilde {\cal B}^{\eps,
R}(E^*):=\{(x,t)\in\RR^n\times(0,\infty):\,(x,t)\in{\cal
B}^{\alpha\eps, \alpha R, 1/\alpha}(E^*)\,\mbox{ for some
}\,1<\alpha<2\}.
\end{equation}

Consider the following three regions:

\begin{eqnarray}\label{eq6.17}
&&\widetilde {\cal B}^{\eps}(E^*):=
\{(x,t)\in\RR^n\times (\eps,2\eps):\,{\mbox{ dist}}(x, E^*)<t\},\\[4pt]\label{eq6.18}
&&\widetilde {\cal B}^{R}(E^*):=
\{(x,t)\in\RR^n\times (R,2R):\,{\mbox{ dist}}(x, E^*)<t\},\\[4pt]\label{eq6.19}
&& \widetilde {\cal B}'(E^*):= \{(x,t)\in B^* \times
(\eps,2R):\,{\mbox{ dist}}(x, E^*)<t<2{\mbox{ dist}}(x, E^*)\},
\end{eqnarray}

\noindent and observe that

\begin{equation}\label{eq6.20}
\widetilde {\cal B}^{\eps, R}(E^*)\,\subset \,\widetilde {\cal
B}^{\eps}(E^*)\cup \widetilde {\cal B}^{R}(E^*) \cup \widetilde
{\cal B}'(E^*).
\end{equation}

\noindent Below we will analyze separately the parts of integrals in
(\ref{eq6.15}) corresponding to the regions
(\ref{eq6.17})--(\ref{eq6.19}).

Let us start with

\begin{equation}\label{eq6.21}
I_1^\eps:=\dint_{\widetilde {\cal
B}^{\eps}(E^*)}|u(y,t)|^2\,\frac{dydt}{t}.
\end{equation}

\noindent For every $(y,t)\in\widetilde {\cal B}^{\eps}(E^*)$ there
exists $y^*\in E^*$ such that $y^*\in B(y,t)$. By definition of
$E^*$ this implies that $|E\cap B(y^*,t)|\geq C|B(y^*,t)|$ and
therefore $|E\cap B(y,2t)|\geq Ct^n$. Then

\begin{eqnarray}\nonumber
I_1^\eps &\leq & C\dint_{\widetilde {\cal
B}^{\eps}(E^*)}\int_{E\cap
B(y,2t)}|u(y,t)|^2\,dz\,\frac{dydt}{t^{n+1}}\\[4pt]\nonumber
&\leq & C\int_\eps^{2\eps}\int_E
\left(\frac{1}{t^n}\int_{B(z,2t)}|u(y,t)|^2\,dy\right)dz\frac{dt}{t}\\[4pt]
&\leq & C\int_\eps^{2\eps}\int_E |{\cal N}_h^\beta
f(z)|^2\,dz\frac{dt}{t}\leq C\int_E |{\cal N}_h^\beta
f(z)|^2\,dz,\label{eq6.22}
\end{eqnarray}

\noindent for every $\beta\geq 2$.

Using similar ideas,

\begin{equation}\label{eq6.23}
I_2^\eps:=\dint_{\widetilde {\cal B}^{\eps}(E^*)}t|\nabla
u(y,t)|^2\,dydt\leq C\int_\eps^{2\eps}\int_E
\left(\frac{1}{t^{n-2}}\int_{B(z,2t)}|\nabla
u(y,t)|^2\,dy\right)dz\frac{dt}{t}.
\end{equation}

\noindent Recall now parabolic Caccioppoli inequality
(\ref{eq2.15}). By definition $u(y,t)=e^{-t^2L}f(y)$, therefore,
making the change of variables in (\ref{eq2.15}), one can see that

\begin{equation}\label{eq6.24}
\int_{t_0-cr}^{t_0}\int_{B(x_0,r)}t|\nabla u(x,t)|^2\,dxdt\leq
\frac{C}{r^2}\int_{t_0-2cr}^{t_0}\int_{B(x_0,2r)}t|u(x,t)|^2\,dxdt,
\end{equation}

\noindent for every $x_0\in\RR^n$, $r>0$, $t_0>2cr$. Here $c>0$ and
the constant $C$ depends on $c$. Next, we divide the integral in
$t\in(\eps, 2\eps)$ from (\ref{eq6.21}) into integrals over $(\eps,
3\eps/2)$ and $(3\eps/2, 2\eps)$, and apply (\ref{eq6.24}) with
$t_0=2\eps$ and $t_0=3\eps/2$, $r=2\eps$, $c=1/4$ to obtain the
bound

\begin{eqnarray}
I_2^\eps &\leq& C\int_{\eps/2}^{2\eps}\int_E
\left(\frac{1}{t^n}\int_{B(z,8\eps)}|u(y,t)|^2\,dy\right)dz\frac{dt}{t} \nonumber\\[4pt]
&\leq & C\int_{\eps/2}^{2\eps}\int_E
\left(\frac{1}{t^n}\int_{B(z,16t)}|u(y,t)|^2\,dy\right)dz\frac{dt}{t}\leq
C\int_E |{\cal N}_h^\beta f(z)|^2\,dz,\label{eq6.25}
\end{eqnarray}
\noindent where $\beta\geq 16$.

Observe that the same argument applies to estimates

\begin{eqnarray}\label{eq6.26}
&&\dint_{\widetilde {\cal B}^{R}(E^*)}|u(y,t)|^2\,\frac{dydt}{t}
\leq C\int_E |{\cal N}_h^\beta f(z)|^2\,dz,\\[4pt]
&&\dint_{\widetilde {\cal B}^{R}(E^*)}t|\nabla u(y,t)|^2\,dydt\leq
C\int_E |{\cal N}_h^\beta f(z)|^2\,dz,\label{eq6.27}
\end{eqnarray}

\noindent with $\beta\geq 16$.

To control the integral over $\widetilde {\cal B}'(E^*)$, we first
decompose $B^*$ into a family of Whitney balls,
$\{B(x_k,r_k)\}_{k=0}^\infty$, such that $\cup_{k=0}^\infty
B(x_k,r_k)=B^*$, $c_1 \,{\rm dist}(x_k, E^*)\leq r_k\leq c_2\,
{\rm dist}(x_k, E^*)$, and every point $x\in B^*$ belongs at most
to $c_3$ balls, $0<c_1<c_2<1$ and $c_3\in\NN$ are some fixed
constants, independent of $B^*$ (see \cite{CW}, \cite{St}). Then

\begin{eqnarray}\nonumber
I_1'&:=&\dint_{\widetilde {\cal
B}'(E^*)}|u(y,t)|^2\,\frac{dydt}{t}\leq \sum_{k=0}^\infty
\int_{r_k(1/c_2-1)}^{2r_k(1/c_1+1)} \int_{B(x_k,r_k)}
|u(y,t)|^2\,\frac{dydt}{t}\\[4pt]
 &\leq & C\sum_{k=0}^\infty
r_k^n\int_{r_k(1/c_2-1)}^{2r_k(1/c_1+1)}
\left[\frac{1}{t^n}\int_{B(x_k,\frac{c_2}{1-c_2}\,t)}
|u(y,t)|^2\,dy\right]\frac{dt}{t}.
 \label{eq6.28}
\end{eqnarray}

\noindent On the other side, $E^*\subset E$, hence, ${\rm
dist}(x_k,E)\leq {\rm dist}(x_k,E^*)\leq
\frac{c_2}{(1-c_2)c_1}\,t$ and the expression in brackets above
can be majorized by the square of non-tangential maximal function
${\cal N}^\beta f(z)$ for some $z\in E$ and
$\beta\geq\frac{c_2}{(1-c_2)c_1}$. Hence,

\begin{equation}\label{eq6.29}
I_1'\leq C\sum_{k=0}^\infty r_k^n \left(\sup_{z\in E}{\cal
N}^\beta_h f(z)\right)^2 \leq C|B^*|\left(\sup_{z\in E}{\cal
N}^\beta_h f(z)\right)^2.
\end{equation}

Similarly to (\ref{eq6.28})--(\ref{eq6.29}) we can prove that there
exists $C_0=C_0(c_1,c_2)>0$ such that

\begin{equation}\label{eq6.30}
I_2':=\dint_{\widetilde {\cal B}'(E^*)}t|\nabla u(y,t)|^2\,dydt \leq
C|B^*| \left(\sup_{z\in E}{\cal N}^\beta_h f(z)\right)^2,
\end{equation}

\noindent for $\beta\geq C_0$, using (\ref{eq6.24}) to control the
gradient of $u$.

Let us choose now

\begin{equation}\label{eq6.31}
\beta:=\max\left\{16,\frac{c_2}{(1-c_2)c_1}, C_0\right\}
\end{equation}

\noindent in (\ref{eq6.6}). Then

\begin{equation}\label{eq6.32}
I_1'\leq C\sigma^2 |B^*|\quad\mbox{and}\quad I_2'\leq C\sigma^2
|B^*|.
\end{equation}

\noindent Combining all the estimates above allows us to write

\begin{equation}\label{eq6.33}
\int_{E^*}\left(\widetilde S_h^{2\eps, R,
1/2}f(x)\right)^{2}\,dx\leq C\sigma^2 |B^*|+C\int_E |{\cal
N}_h^\beta f(z)|^2\,dz,
\end{equation}

\noindent and therefore, passing to the limit as $\eps\to 0$ and
$R\to \infty$,

\begin{equation}\label{eq6.34}
\int_{E^*}\left(\widetilde S_h^{1/2}f(x)\right)^{2}\,dx\leq
C\sigma^2 |B^*|+C\int_E |{\cal N}_h^\beta f(z)|^2\,dz.
\end{equation}

Denote by $\lambda_{{\cal N}_h^\beta f}$ the distribution function
of ${\cal N}_h^\beta f$ and recall that ${\cal N}_h^\beta\leq
\sigma$ on $E$. Then

\begin{equation}\label{eq6.35}
\int_{E^*}\left(\widetilde S_h^{1/2}f(x)\right)^{2}\,dx\leq
C\sigma^2 \lambda_{{\cal N}_h^\beta f}(\sigma)+C\int_0^\sigma
t\lambda_{{\cal N}_h^\beta f}(t)\,dt,
\end{equation}

\noindent since $|B^*|\leq C|B|\leq C\lambda_{{\cal N}_h^\beta
f}(\sigma)$. Next,

\begin{eqnarray}\nonumber
&&\hskip -0.7cm  \lambda_{\widetilde S_h^{1/2}f}(\sigma)\leq
|\{x\in
E^*:\,\widetilde S_h^{1/2}f(x)>\sigma\}|+|\,^cE^*|\\[4pt]
&& \hskip -0.7cm \quad \leq C \frac{1}{\sigma^2}
\int_{E^*}\left(\widetilde S_h^{1/2}f(x)\right)^{2}\,dx+C
\lambda_{{\cal N}_h^\beta f}(\sigma)\leq
C\frac{1}{\sigma^2}\int_0^\sigma t\lambda_{{\cal N}_h^\beta
f}(t)\,dt+C\lambda_{{\cal N}_h^\beta f}(\sigma), \label{eq6.36}
\end{eqnarray}

\noindent and therefore,

\begin{equation}\label{eq6.37}
\|\widetilde S_h^{1/2}f\|_{L^1(\RR^n)}=\int_0^\infty
\lambda_{\widetilde S_h^{1/2}f}(\sigma)\,d\sigma\leq C \|{\cal
N}_h^\beta f\|_{L^1(\RR^n)},
\end{equation}

\noindent for $\beta$ as in (\ref{eq6.31}). In view of
Lemma~\ref{l2.2} and (\ref{eq6.2}) the theorem is proved modulo the
result we present below. \ep

\begin{lemma}\label{l6.2} For all $f\in L^2(\RR^n)$ and $\beta\geq 1$

\begin{equation}\label{eq6.38}
\|{\cal N}_h^\beta f\|_{L^1(\RR^n)}\leq C\beta^n\|{\cal N}_h^1
f\|_{L^1(\RR^n)}.
\end{equation}
\end{lemma}

\bp Fix $\sigma\in (0,\infty)$ and consider the following sets:

\begin{equation}\label{eq6.39}
E_\sigma:=\{x\in\RR^n:\,{\cal N}^1_h
f(x)>\sigma\}\quad\mbox{and}\quad E_\sigma^*:=\{x\in\RR^n:\,{\cal
M}(\chi_{E_\sigma})(x)>C/\beta^n\}.
\end{equation}

\noindent It is not hard to see that $|E_\sigma^*|\leq C\beta^n
|E_\sigma|$.

Assume now that $x\not\!\in E_\sigma^*$. Then $B(y,t)\not\!\subset
E_\sigma$ for every $(y,t)\in \Gamma^{2\beta}(x)$. Indeed, if
$B(y,t)\subset E_\sigma$, then

\begin{equation}\label{eq6.40}
{\cal M}(\chi_{E_\sigma})(x)>C\frac{|B(y,t)|}{|B(x,2\beta t)|}\geq C
/\beta^n,
\end{equation}

\noindent which implies $x\in E_\sigma^*$.

Therefore, there exists $z\in B(y,t)$ such that ${\cal
N}^1_hf(z)\leq\sigma$, in particular,

\begin{equation}\label{eq6.41}
\left(\frac{1}{t^n}\int_{B(y,t)}
|e^{-t^2L}f(z)|^2\,dz\right)^{1/2}\leq \sigma.
\end{equation}

\noindent Recall that the above inequality holds for all $(y,t)\in
\Gamma^{2\beta}(x)$. Now for every $w\in B(x,\beta t)$ one can cover
$B(w,\beta t)$ by $C\beta^n$ balls $B(y_i,t)$, where $y_i\in
\Gamma^{2\beta}(x)$, to prove that

\begin{equation}\label{eq6.42}
\frac{1}{(\beta t)^n}\int_{B(w,\beta t)} |e^{-t^2L}f(z)|^2\,dz\leq
\frac{C}{\beta^{n}} \sum_{i}\frac{1}{t^n}\int_{B(y_i,t)}
|e^{-t^2L}f(z)|^2\,dz\leq C\sigma^2,
\end{equation}

\noindent hence,

\begin{equation}\label{eq6.43}
{\cal N}^\beta_h f(x)\leq C\sigma \quad\mbox{for every} \quad
x\not\!\in E_\sigma^*.
\end{equation}

\noindent Having this at hand, we simply write

\begin{eqnarray}\nonumber
&&\|{\cal N}^\beta_h f\|_{L^1(\RR^n)}\leq C \int_0^\infty
|\{x\in\RR^n:\,{\cal N}^\beta_h f>C\sigma\}|\,d\sigma\\[4pt]
&&\quad\leq C\int_0^\infty |E_\sigma^*|\,d\sigma\leq C\int_0^\infty
\beta^n|E_\sigma|\,d\sigma\leq C\beta^n \|{\cal N}^1_h
f\|_{L^1(\RR^n)}, \label{eq6.44}
\end{eqnarray}

\noindent and finish the argument. \ep

\begin{theorem}\label{t6.3}
Let $\eps>0$ and $M>n/4$.  Then for every representation
$\sum_{i=0}^\infty \lambda_i m_i$, where $\{m_i\}_{i=0}^\infty$ is a family of
$(2,\eps,M)$-molecules  and
$\sum_{i=0}^\infty |\lambda_i|<\infty$, the series
$\sum_{i=0}^\infty \lambda_i m_i$ converges in $H^1_{{\cal
N}_h}(\RR^n)$ and

\begin{equation}\label{eq6.45}
\Bigl\|\sum_{i=0}^\infty \lambda_i m_i\Bigr\|_{H^1_{{\cal
N}_h}(\RR^n)}\leq C\sum_{i=0}^\infty |\lambda_i|.
\end{equation}
\end{theorem}

\bp As usual, by Lemma \ref{l3.sub} we need only establish a uniform $L^1$ bound on molecules.
Consider the following modifications of the non-tangential
maximal function

\begin{eqnarray}\label{eq6.46}
{\cal N}_h^* f(y) & = & \sup_{t>0} \left( \frac{1}{t^n}
\int_{|x-y|<t}
|e^{-t^2 L} f(x) |^2 dx\right)^{1/2},\\[4pt]\label{eq6.47}
{\cal N}_h^{*,M} f(y) & = & \sup_{t>0} \left( \frac{1}{t^n}
\int_{|x-y|<t} |t^{2M} L^M e^{-t^2 L} f(x) |^2 dx\right)^{1/2},
\end{eqnarray}

\noindent where $y\in\RR^n$, $M\in\NN$ and $f\in L^2(\RR^n)$. Both
of the operators above are bounded on $L^2(\RR^n)$.

Indeed,

\begin{eqnarray}
\|{\cal N}_h^* f\|_{L^2(\RR^n)}^2 & \leq &
C\int_{\RR^n}\left[\sup_{t>0}\sum_{j=0}^\infty \left( \frac{1}{t^n}
\int_{B(y,t)} |e^{-t^2 L} (f\chi_{S_j(B(y,t))})(x) |^2
dx\right)^{1/2}\right]^2\,dy\nonumber\\[4pt]
& \leq & C\int_{\RR^n}\left[\sup_{t>0}\sum_{j=0}^\infty
\frac{1}{t^{n/p}}\, e^{-\frac{{\rm
dist}\,(B(y,t),S_j(B(y,t)))^2}{ct^2}}
\|f\|_{L^p(S_j(B(y,t)))}\right]^2\,dy,\label{eq6.48}
\end{eqnarray}

\noindent for every $p_L< p \leq 2$ by $L^p-L^2$ off-diagonal
estimates. Therefore, for every $\eps>0$ and $p<2$ as above

\begin{eqnarray}
\|{\cal N}_h^* f\|_{L^2(\RR^n)}^2 & \leq &
C\int_{\RR^n}\left[\sup_{t>0}\sum_{j=0}^\infty \frac{1}{t^{n/p}}\,
2^{-j(n/p+\eps)}
\|f\|_{L^p(S_j(B(y,t)))}\right]^2\,dy\nonumber\\[4pt]
&\leq &C\int_{\RR^n}\left[{\cal M}(|f|^p)(y)\right]^{2/p}\,dy\leq C
\int_{\RR^n}|f(y)|^2\,dy,\label{eq6.49}
\end{eqnarray}

\noindent using $L^{2/p}(\RR^n)$ boundedness of the Hardy-Littlewood
maximal function.

Along the same lines we can prove $L^2$ boundedness of the function
${\cal N}_h^{*,M}$.

On the other side, by Lemma~\ref{l6.2}

\begin{equation}\label{eq6.50}
\|{\cal N}_hf\|_{L^1(\RR^n)}\leq C \|{\cal
N}_h^{1/2}f\|_{L^1(\RR^n)} \leq C\|{\cal N}_h^*f\|_{L^1(\RR^n)},
\end{equation}

\noindent and therefore, it is enough to show that

\begin{equation}\label{eq6.51}
\|{\cal N}_h^*m\|_{L^1(\RR^n)}\leq C
\end{equation}

\noindent for every $m$ a $(2,\eps,M)$-molecule associated to some
cube $Q$.

To this end, we use the annular decomposition of $\RR^n$ along
with H\"older's inequality  to write

\begin{eqnarray}\nonumber
&&\hskip -0.7cm\|{\cal N}_h^*m\|_{L^1(\RR^n)}\leq C
\sum_{j=0}^\infty(2^jl(Q))^{n/2}\|{\cal
N}_h^*m\|_{L^2(S_j(Q))}\\[4pt]
&&\hskip -0.7cm\qquad \leq C \sum_{j=0}^{10}(2^jl(Q))^{n/2}\|{\cal
N}_h^*m\|_{L^2(S_j(Q))} +C
\sum_{j=10}^\infty(2^jl(Q))^{n/2}\|{\cal N}_h^*m\|_{L^2(S_j(Q))}.
\label{eq6.52}
\end{eqnarray}

\noindent The finite sum above is bounded by some constant in view
of $L^2(\RR^n)$ boundedness of ${\cal N}_h^*$ and (\ref{eq1.8})
condition on molecules.

To handle the second sum in (\ref{eq6.52}), we fix some number
$0<a<1$ such that $n/2-2aM<0$ and split ${\cal N}_h^*m$ according
to whether $t \leq c\, 2^{aj}l(Q)$ or $t \geq c\, 2^{aj}l(Q)$.
Consider the former case first. Set

\begin{equation}
U_j (Q) := 2^{j+3} Q \setminus 2^{j-3} Q, \quad R_j (Q) := 2^{j+5}
Q \setminus 2^{j-5} Q, \quad \mbox{and}\quad E_j(Q)=\,^cR_j(Q),
\label{eq6.53}
\end{equation}

\noindent for every $j\geq 10$ and split $m = m \chi_{R_j(Q)} + m
\chi_{E_j(Q)}.$

For $x \in S_j(Q)$, $|x-y|<t$ and $t \leq c\, 2^{aj}l(Q)$ we have
$y \in U_j(Q)$. Moreover,\break ${\rm dist}(U_j(Q),E_j(Q)) \approx
C 2^jl(Q)$. Then the Gaffney estimates (Lemma~\ref{l2.4})
guarantee that for every such $t$, $y\in\RR^n$, $a<1$ and
$N\in\NN$

\begin{eqnarray}\nonumber
&&\left(\frac{1}{t^n} \int_{|x-y|<t} |e^{-t^2 L} (m
\chi_{E_j(Q)})(x) |^2\,dx\right)^{1/2} \\[4pt]
&&\qquad\leq \frac{C}{t^{n/2}}
\,e^{-\frac{(2^jl(Q))^2}{ct^2}}\,\|m\|_{L^2(E_j(Q))}\leq
\frac{C}{t^{n/2}}
\left(\frac{t}{2^jl(Q)}\right)^N\|m\|_{L^2(\RR^n)},\label{eq6.54}
\end{eqnarray}

\noindent so that

\begin{eqnarray}\nonumber
&&\sum_{j=10}^\infty (2^jl(Q))^{n/2}\left\|\sup_{t\leq c\,
2^{aj}l(Q)} \left( \frac{1}{t^n} \int_{|x-\,\cdot|\,<t} |e^{-t^2 L}
(m
\chi_{E_j(Q)})(x)|^2 dx\right)^{1/2}\right\|_{L^2(S_j(Q))}\\[4pt]
&&\qquad  \leq \sum_{j=10}^\infty 2^{j(1-a)(n/2-N)} \leq
C,\label{eq6.55}
\end{eqnarray}

\noindent when $N>n/2$.

As for the contribution of $m\chi_{R_j(Q)}$, by $L^2$ boundedness
of ${\cal N}_h^*$ we have

\begin{eqnarray}\nonumber
&&\sum_{j=10}^\infty (2^jl(Q))^{n/2}\| {\cal N}_h^* (m
\chi_{R_j(Q)}) \|_{L^2(S_j(Q))}\\[4pt]
&&\qquad \leq C \sum_{j=10}^\infty (2^jl(Q))^{n/2}\| m
\chi_{R_j(Q)}\|_{L^2(R_j(Q))} \leq C \sum_{j=10}^\infty 2^{-j
\eps}\leq C. \label{eq6.56}
\end{eqnarray}

Now we consider the case $t \geq c \,2^{aj}l(Q)$.  For every
$y\in\RR^n$

\begin{eqnarray}\nonumber
&&\sup_{t \geq c \,2^{aj}l(Q)} \left( \frac{1}{t^n} \int_{|x-y|<t}
|e^{-t^2 L} m(x) |^2 dx\right)^{1/2}\\[4pt]\nonumber
&&\qquad  = \sup_{t \geq c \,2^{aj}l(Q)} \left( \frac{1}{t^n}
\int_{|x-y|<t} |(t^{2M} L^Me^{-t^2 L}) (t^{-2M}L^{-M} m )(x) |^2
dx\right)^{1/2}\\[4pt]
&&\qquad  \leq C \,2^{-2aMj} {\cal N}_h^{*,M} (l(Q)^{-2M} L^{-M}
m)(y), \label{eq6.57}
\end{eqnarray}

\noindent so we use the boundedness of ${\cal N}_h^{*,M}$ on
$L^2(\RR^n)$ to finish the argument.\ep

\begin{corollary}\label{c6.4}
$H^1_L(\RR^n)=H^1_{{\cal N}_h}(\RR^n)$, in particular,
$\|f\|_{H^1_L(\RR^n)}\approx\|{\cal N}_h f\|_{L^1(\RR^n)}+\| f\|_{L^1(\RR^n)}$.
\end{corollary}

\bp The right-to-the-left inclusion is a direct consequence of
Theorems~\ref{t4.1} and \ref{t6.1}, the converse follows from
Theorem~\ref{t6.3} and Corollary~\ref{c4.3}. \ep

\section{Characterization by the non-tangential maximal function
associated to the Poisson semigroup} \setcounter{equation}{0}

\begin{theorem}\label{t7.1}
For every $f\in L^2(\RR^n)$

\begin{equation}\label{eq7.1}
\|\widehat S_Pf\|_{L^1(\RR^n)}\leq C\|{\cal N}_Pf\|_{L^1(\RR^n)},
\end{equation}

\noindent where $\widehat S_P$ is the operator defined in
(\ref{eq5.30}).
\end{theorem}

\bp The proof of this fact closely follows the argument of
Theorem~\ref{t6.1}. More precisely, at the step corresponding to
(\ref{eq6.11}) we assign

\begin{equation}\label{eq7.2}
u(y,t):=e^{-t\sqrt L}f(y), \qquad t\in(0,\infty), \qquad y\in\RR^n.
\end{equation}

\noindent The analogue of (\ref{eq6.13})--(\ref{eq6.14}) can be
obtained observing that this time $-{\rm div}_{y,t} B \nabla_{y,t}
u=0$, where as before $B$ is the $(n+1)\times(n+1)$ block diagonal
matrix with entries 1 and $A$ and ${\rm div}_{y,t}$ is divergence
in space and time variables. Concretely, we can write

\begin{eqnarray}\nonumber
&&\hskip -.7cm \dint_{{\cal R}^{\alpha\eps, \alpha R, \frac 1
\alpha}(E^*)}
t|\nabla_{y,t} u(y,t)|^2\,dydt  \\[4pt]\nonumber
&&\hskip -.7cm\leq C \Re e \dint_{{\cal R}^{\alpha\eps, \alpha R,
\frac 1 \alpha}(E^*)} \left[-\,{\rm div}_{y,t}[tB(y)\nabla_{y,t}
u(y,t)] {\overline{u(y,t)}}-u(y,t)\overline{{{\rm
div}_{y,t}[tB(y)\nabla_{y,t}
u(y,t)]}}\right]dydt\\[4pt]\nonumber
&&\hskip -.7cm\quad +C\Re e \int_{{\cal B}^{\alpha\eps, \alpha R,
\frac 1 \alpha}(E^*)} \left[tB(y)\nabla_{y,t} u(y,t)\cdot
\nu(y,t){\overline{u(y,t)}}+u(y,t)\nu(y,t)\cdot
t{\overline{B(y)\nabla_{y,t} u(y,t)}}\right]d\sigma_{y,t}\\[4pt]\nonumber
&&\hskip -.7cm\leq C\Re e\dint_{{\cal R}^{\alpha\eps, \alpha R,
\frac 1 \alpha}(E^*)} \left[-\partial_t u(y,t)\cdot
{\overline{u(y,t)}}-u(y,t)\cdot
{\overline{\partial_t u(y,t)}}\right]dydt\\[4pt]\nonumber
&&\hskip -.7cm\quad +C\Re e \int_{{\cal B}^{\alpha\eps, \alpha R,
\frac 1 \alpha}(E^*)} \left[tB(y)\nabla_{y,t} u(y,t)\cdot
\nu(y,t){\overline{u(y,t)}}+u(y,t)\nu(y,t)\cdot
t{\overline{B(y)\nabla_{y,t} u(y,t)}}\right]d\sigma_{y,t},
\end{eqnarray}

\noindent so that (\ref{eq6.14})--(\ref{eq6.15}) with $\nabla_{y,t}$
in place of space gradient holds.

The rest of the argument is essentially the same as the proof of Theorem~\ref{t6.1},
just employing elliptic instead of parabolic Caccioppoli inequality.
\ep

To handle the converse to (\ref{eq7.1}), we start with two
auxiliary Lemmas.

\begin{lemma}\label{l7.2} Define

\begin{equation}\label{eq7.3}
\bar g_Pf(x):= \left(\int_0^\infty|t\sqrt L e^{-t \sqrt
L}f(x)|^2\,\frac {dt}{t}\right)^{1/2}.
\end{equation}

\noindent Then

\begin{equation}\label{eq7.4}
\bar g_P f(x)\leq C g_h f(x), \qquad x\in\RR^n,
\end{equation}

\noindent for every $f\in L^2(\RR^n)$.
\end{lemma}

\bp By subordination formula (\ref{eq5.2}) and Minkowski
inequality

\begin{eqnarray}\nonumber
\bar g_Pf(x)&=& \left(\int_0^\infty|t\partial_t e^{-t \sqrt
L}f(x)|^2\,\frac {dt}{t}\right)^{1/2}\\[4pt]\label{eq7.5}
&\leq & C\int_0^\infty\frac{e^{-u}}{\sqrt
u}\,\,\left(\int_0^\infty|t\partial_t
e^{-\frac{t^2L}{4u}}f(x)|^2\,\frac {dt}{t}\right)^{1/2}\,du.
\end{eqnarray}

\noindent After the change of variables $t\mapsto
s=\frac{t^2}{4u}$, $s\in (0,\infty)$,
$\frac{dt}{t}=\frac{ds}{2s}$, the expression above can be written
as

\begin{equation}
C\int_0^\infty\frac{e^{-u}}{\sqrt u}\,\,\left(\int_0^\infty|sL
e^{-sL}f(x)|^2\,\frac {ds}{s}\right)^{1/2}\,du\leq Cg_h(x),
\label{eq7.6}\end{equation}

\noindent as desired. \ep

\begin{lemma}\label{l7.3} Define

\begin{equation}\label{eq7.7}
g_P^{aux}f(x):= \left(\int_0^\infty|(e^{-t\sqrt L}-
e^{-t^2L})f(x)|^2\,\frac {dt}{t}\right)^{1/2}.
\end{equation}

\noindent  Then

\begin{equation}\label{eq7.8}
g_P^{aux} f(x)\leq C g_h f(x), \qquad x\in\RR^n,
\end{equation}

\noindent for every $f\in L^2(\RR^n)$.
\end{lemma}

\bp By subordination formula (\ref{eq5.2})

\begin{eqnarray}\nonumber
g_P^{aux}f(x)&=&
C\left(\int_0^\infty\left|\int_0^\infty\frac{e^{-u}}{\sqrt
u}(e^{-\frac{t^2L}{4u}}-
e^{-t^2L})f(x) \,du\right|^2\,\frac {dt}{t}\right)^{1/2}\\[4pt]\label{eq7.9}
&=& C\left(\int_0^\infty\left|\int_0^\infty\frac{e^{-u}}{\sqrt
u}\int_t^{t/\sqrt{4u}} 2rLe^{-r^2L}f(x)\,dr du\right|^2\,\frac
{dt}{t}\right)^{1/2}.
\end{eqnarray}

\noindent We now split the integral in $u$ according to whether
$u<1/4$ or $u>1/4$. In the first case,

\begin{eqnarray}\nonumber
&&\left|\int_0^{1/4}\frac{e^{-u}}{\sqrt
u}\int_t^{t/\sqrt{4u}} 2rLe^{-r^2L}f(x)\,dr du\right|\\[4pt]
&&\quad\leq \int_t^{\infty}\left|\int_0^{t^2/(4r^2)}
\frac{e^{-u}}{\sqrt u}\,du\right|\, |2rLe^{-r^2L}f(x)|\,dr \leq
\int_t^{\infty}|tLe^{-r^2L}f(x)|\,dr. \label{eq7.10}\end{eqnarray}

\noindent As for the second part,

\begin{eqnarray}\nonumber
&&\left|\int_{1/4}^\infty \frac{e^{-u}}{\sqrt
u}\int_{t/\sqrt{4u}}^t 2rLe^{-r^2L}f(x)\,dr du\right|\\[4pt]
&&\quad= \int_0^t\left|\int_{t^2/(4r^2)}^\infty
\frac{e^{-u}}{\sqrt u}\,du\right|\, |2rLe^{-r^2L}f(x)|\,dr \leq
C\int_0^t|(r^2/t)\, Le^{-r^2L}f(x)|\,dr.
\label{eq7.11}\end{eqnarray}

\noindent Inserting the results into (\ref{eq7.9}), we get

\begin{eqnarray}\nonumber
g_P^{aux}f(x)&\leq & C\left(\int_0^\infty t^2 \left(\int_t^\infty
|rLe^{-r^2L}f(x)|\,\frac{dr}{r}\right)^{2}\frac{dt}{t}\right)^{1/2}\\[4pt]
\nonumber &&\qquad + C\left(\int_0^\infty \frac{1}{t^2}
\left(\int_0^t
|r^2Le^{-r^2L}f(x)|\,dr\right)^{2}\frac{dt}{t}\right)^{1/2}\\[4pt]
\nonumber &\leq & C\left(\int_0^\infty \int_t^\infty
|rLe^{-r^2L}f(x)|^2\,drdt\right)^{1/2}+ C\left(\int_0^\infty
\int_0^t
|r^2Le^{-r^2L}f(x)|^2\,dr\frac{dt}{t^2}\right)^{1/2}\\[4pt]
\nonumber &\leq& C\left(\int_0^\infty
|r^2Le^{-r^2L}f(x)|^2\,\frac{dr}{r}\right)^{1/2}=Cg_hf(x).
\end{eqnarray}

\noindent This finishes the argument. \ep

\begin{theorem}\label{t7.4}
Let $\eps>0$ and $M>n/4$.  Then for every representation
$\sum_{i=0}^\infty \lambda_i m_i$, where $\{m_i\}_{i=0}^\infty$ is
a family of $(2,\eps,M)$-molecules  and $\sum_{i=0}^\infty
|\lambda_i|<\infty$, the series $\sum_{i=0}^\infty \lambda_i m_i$
converges in $H^1_{{\cal N}_P}(\RR^n)$ and

\begin{equation}\label{eq6.45}
\Bigl\|\sum_{i=0}^\infty \lambda_i m_i\Bigr\|_{H^1_{{\cal
N}_P}(\RR^n)}\leq C\sum_{i=0}^\infty |\lambda_i|.
\end{equation}
\end{theorem}

\bp Let $2<p<\widetilde{p}_L$, $\epsilon>0$, $M>\frac n4$. Similarly to
Theorem~\ref{t6.3}, it is enough to prove that

\begin{equation}\label{eq7.13}
\|{\cal N}_P^* m\|_{L^1(\RR^n)}\leq C,
\end{equation}

\noindent for every $(p,\epsilon, M)$-molecule $m$, where

\begin{equation}\label{eq7.14} {\cal N}_P^* f(x): = \sup_{t>0}
\left( \frac{1}{t^n} \int_{|x-y|<t} |e^{-t\sqrt L} f(y) |^2
dy\right)^{1/2},\quad x\in\RR^n, \quad f\in L^2(\RR^n).
\end{equation}

\noindent To this end, by the standard dyadic annular decomposition and H\"{o}lder's inequality,
it will suffice to establish the estimate

\begin{equation}\label{eq7.15}
\|{\cal N}_P^* m\|_{L^p(S_j(Q))}\leq C (2^jl(Q))^{n\left(\frac
1p-1\right)}2^{-j\gamma}, \qquad j\in\NN\cup\{0\},
\end{equation}

\noindent where $Q$ is a cube associated to the molecule $m$ and
$\gamma$ is some fixed positive number.

Fix some $a$ such that $n\left(\frac n2 +2M\right)^{-1}<a<1$. Then

\begin{eqnarray}\nonumber
&&\left\|\sup_{t\geq 2^{aj}l(Q)} \left( \frac{1}{t^n} \int_{|\cdot
-y|<t} |e^{-t\sqrt L} m(y) |^2
dy\right)^{1/2}\right\|_{L^p(S_j(Q))}\\[4pt]\nonumber
&&\quad = \left\|\sup_{t\geq 2^{aj}l(Q)}
\left(\frac{l(Q)}{t}\right)^{2M}\left( \frac{1}{t^n}
\int_{|\cdot-y|<t} |(t^2L)^M e^{-t\sqrt L} (l(Q)^2L)^{-M} m(y) |^2
dy\right)^{1/2}\right\|_{L^p(S_j(Q))}\\[4pt]\nonumber
&&\quad \leq C
\left(\frac{1}{2^{aj}}\right)^{2M}\left(\frac{1}{2^{aj}l(Q)}\right)^{n/2}
\left\|\sup_{t>0} \left( \int_{\RR^n} |(t^2L)^M e^{-t\sqrt L}
(l(Q)^2L)^{-M} m(y) |^2
dy\right)^{1/2}\right\|_{L^p(S_j(Q))}.\nonumber
\end{eqnarray}

\noindent Resting on Lemma~\ref{l5.1}, one can prove that for
$M>n/4$ the family of operators $(t^2L)^M e^{-t\sqrt L}$ is
uniformly bounded in $L^2(\RR^n)$. Also,

\begin{equation}\label{eq7.16}
\|(l(Q)^2L)^{-M} m\|_{L^2(\RR^n)}\leq C |Q|^{-1/2}, \qquad
j\in\NN,
\end{equation}

\noindent by the definition of molecule and H\"older's inequality.
Then

\begin{eqnarray}\nonumber
&& \left\|\sup_{t\geq 2^{aj}l(Q)} \left( \frac{1}{t^n}
\int_{|\cdot -y|<t} |e^{-t\sqrt L} m(y) |^2
dy\right)^{1/2}\right\|_{L^p(S_j(Q))}\\[4pt]
&&\qquad \leq C 2^{aj(-2M-n/2)}2^{jn/p}l(Q)^{n/p-n}= C
(2^jl(Q))^{n\left(\frac 1p-1\right)}2^{-j\eps_1}, \label{eq7.17}
\end{eqnarray}

\noindent for $\eps_1=a(2M+n/2)-n>0$ by the assumptions on $a$.

Turning to the case $t\leq 2^{aj}l(Q)$, we follow a suggestion of P. Auscher, and split

\begin{eqnarray}\nonumber
&& \sup_{t\leq 2^{aj}l(Q)} \left( \frac{1}{t^n} \int_{|x -y|<t}
|e^{-t\sqrt L} m(y) |^2
dy\right)^{1/2}\\[4pt]
&&\qquad \leq \sup_{t\leq 2^{aj}l(Q)} \left( \frac{1}{t^n}
\int_{|x -y|<t} |(e^{-t\sqrt L}-e^{-t^2L}) m(y) |^2
dy\right)^{1/2}+{\cal N}_hm(x). \label{eq7.18}
\end{eqnarray}

\noindent We remark that Auscher has observed \cite{PAPersonal} that
this splitting yields $L^2$ boundedness of ${\cal N}_P$;
a similar idea has appeared previously in the work of Stein \cite{StPNAS}.
An argument analogous to Theorem~\ref{t6.3} shows that
${\cal N}_hm$ satisfies the desired estimate, so we will
concentrate on the first term on the right hand side of (\ref{eq7.18}). Observe that

\begin{eqnarray}\nonumber
&& t|(e^{-t\sqrt L}-e^{-t^2L}) m(y) |^2 =\left|\int_0^t
\partial_s \left(s^{1/2}(e^{-s\sqrt L}-e^{-s^2L}) m(y)\right)\,ds
\right|^2\\[4pt]\nonumber
&&\quad \leq\left|\int_0^t s^{1/2}
\partial_s (e^{-s\sqrt L}-e^{-s^2L}) m(y)\,ds +\frac 12 \int_0^t s^{-1/2}
(e^{-s\sqrt L}-e^{-s^2L}) m(y)\,ds \right|^2 \\[4pt]\nonumber
&&\quad \leq C t \Bigl(\int_0^t |(e^{-s\sqrt L}-e^{-s^2L})
m(y)|^2\,\frac{ds}{s} \\[4pt]
&&\qquad\qquad   +\int_0^t |s\sqrt L e^{-s\sqrt L}
m(y)|^2\,\frac{ds}{s} +\int_0^t |s^2 L e^{-s^2 L}
m(y)|^2\,\frac{ds}{s} \Bigr).\label{eq7.19}
\end{eqnarray}

\noindent Given Lemmas~\ref{l7.2} and \ref{l7.3}, this allows to
control the first term in (\ref{eq7.18}) by

\begin{equation}\label{eq7.20} \sup_{t\leq 2^{aj}l(Q)} \left( \frac{1}{t^n}
\int_{|x -y|<t} \int_0^\infty |s^2 L e^{-s^2 L}
m(y)|^2\,\frac{ds}{s} \,dy\right)^{1/2},\qquad x\in\RR^n.
\end{equation}

Much as for $t\geq 2^{aj}l(Q)$, we have

\begin{eqnarray}\nonumber
&&\sup_{t\leq 2^{aj}l(Q)} \left( \frac{1}{t^n} \int_{|x -y|<t}
\int_{2^{aj}l(Q)}^\infty |s^2 L e^{-s^2 L}
m(y)|^2\,\frac{ds}{s} \,dy\right)^{1/2}\\[4pt]
&&\quad \leq \left({\cal M} \left( \int_{2^{aj}l(Q)}^\infty
\left(\frac{l(Q)}{s}\right)^{4M} |(s^2L)^M e^{-s^2 L}
(l(Q)^2L)^{-M} m(y)|^2\,\frac{ds}{s}
\right)(x)\right)^{1/2},\label{eq7.21}
\end{eqnarray}

\noindent where ${\cal M}$ denotes the Hardy-Littlewood maximal operator.  Thus,

\begin{eqnarray}\nonumber
&&\left\|\sup_{t\leq 2^{aj}l(Q)} \left( \frac{1}{t^n} \int_{|\cdot
-y|<t} \int_{2^{aj}l(Q)}^\infty |s^2 L e^{-s^2 L}
m(y)|^2\,\frac{ds}{s} \,dy\right)^{1/2}\right\|_{L^p(S_j(Q))}\\[4pt]\nonumber
&&\quad\leq  C \left( \int_{2^{aj}l(Q)}^\infty
\left(\frac{l(Q)}{s}\right)^{4M} \left(\int_{\RR^n} |(s^2L)^M
e^{-s^2 L} (l(Q)^2L)^{-M} m(x)|^p\,dx\right)^{2/p}\frac{ds}{s}
\right)^{1/2}\\[4pt]\nonumber
&&\quad\leq  C \left( \int_{2^{aj}l(Q)}^\infty
\left(\frac{l(Q)}{s}\right)^{4M} s^{-2\left(\frac n2-\frac
np\right)}\frac{ds}{s} \right)^{1/2}\|(l(Q)^2L)^{-M}
m\|_{L^2(\RR^n)}\\[4pt]
&&\quad\leq  C 2^{aj(-2M-n/2+n/p)}l(Q)^{n\left(\frac
1p-1\right)}=C (2^jl(Q))^{n\left(\frac
1p-1\right)}2^{-j\eps_2},\label{eq7.22}
\end{eqnarray}

\noindent where $\eps_2=a(2M+n/2)-n+(1-a)n/p>\eps_1>0$ by our
assumptions on $a$.

It remains to estimate

\begin{equation}\label{eq7.23} \left\|\sup_{t\leq 2^{aj}l(Q)} \left( \frac{1}{t^n}
\int_{|\cdot -y|<t} \int_0^{2^{aj}l(Q)}|s^2 L e^{-s^2 L}
m(y)|^2\,\frac{ds}{s} \,dy\right)^{1/2}\right\|_{L^p(S_j(Q))}.
\end{equation}

\noindent Consider first the case $j\geq 10$. Observe that for
$x\in S_j(Q)$, $j\geq 10$, and $|x-y|<t$ we have $y\in U_j(Q)$, a slightly
fattened version of $S_j(Q)$
(see (\ref{eq6.53})). Then, in the notation of (\ref{eq6.53}),

\begin{eqnarray}\nonumber
&&\left\|\sup_{t\leq 2^{aj}l(Q)} \left( \frac{1}{t^n} \int_{|\cdot
-y|<t} \int_0^{2^{aj}l(Q)}|s^2 L e^{-s^2 L}
(m\chi_{R_j(Q)})(y)|^2\,\frac{ds}{s}
\,dy\right)^{1/2}\right\|_{L^p(S_j(Q))}\\[4pt]\nonumber
&&\qquad \leq C\left\|\left({\cal M} \Bigl(\int_0^{2^{aj}l(Q)}|s^2
L e^{-s^2 L} (m\chi_{R_j(Q)})|^2\,\frac{ds}{s}
\Bigr)\right)^{1/2}\right\|_{L^p(S_j(Q))} \\[4pt] &&\qquad \leq C\left\|g_h
(m\chi_{R_j(Q)})\right\|_{L^p(\RR^n)}\leq C
\left\|m\right\|_{L^p(R_j(Q))}\leq C (2^jl(Q))^{n\left(\frac
1p-1\right)}2^{-j\epsilon}, \label{eq7.24}
\end{eqnarray}

\noindent where the next-to-the-last inequality follows from
$L^p$-boundedness of $g_h$ for $p_L<p<\widetilde{p}_L$ (see
\cite{AuscherSurvey}) and the last inequality follows from the
definition of molecule. On the other hand,

\begin{eqnarray}\nonumber
&&\left\|\sup_{t\leq 2^{aj}l(Q)} \left( \frac{1}{t^n} \int_{|\cdot
-y|<t} \int_0^{2^{aj}l(Q)}|s^2 L e^{-s^2 L}
(m\chi_{E_j(Q)})(y)|^2\,\frac{ds}{s}
\,dy\right)^{1/2}\right\|_{L^p(S_j(Q))}\\[4pt]\nonumber
&&\qquad \leq C\left\|\left({\cal M}
\Bigl(\chi_{U_j(Q)}\int_0^{2^{aj}l(Q)}|s^2 L e^{-s^2 L}
(m\chi_{E_j(Q)})|^2\,\frac{ds}{s}
\Bigr)\right)^{1/2}\right\|_{L^p(S_j(Q))} \\[4pt]
&&\qquad\leq C\left(\int_0^{2^{aj}l(Q)}\|s^2 L e^{-s^2 L}
(m\chi_{E_j(Q)})\|^2_{L^p(U_j(Q))}\,\frac{ds}{s} \right)^{1/2}\leq
C 2^{(a-1)jN}\|m\|_{L^p(\RR^n)}, \nonumber
\end{eqnarray}

\noindent where $N$ is any natural number and the last inequality
follows from the Gaffney estimates. Clearly, we can take $N$ large
enough to bound the expression above by $C (2^jl(Q))^{n\left(\frac
1p-1\right)}2^{-j\epsilon}$.

Finally, in the case $j\leq 10$ following (\ref{eq7.24}) we show

\begin{eqnarray}\nonumber
&&\left\|\sup_{t\leq 2^{aj}l(Q)} \left( \frac{1}{t^n} \int_{|\cdot
-y|<t} \int_0^{2^{aj}l(Q)}|s^2 L e^{-s^2 L} m(y)|^2\,\frac{ds}{s}
\,dy\right)^{1/2}\right\|_{L^p(S_j(Q))}\\[4pt]
 &&\qquad \leq C\left\|g_h
m\right\|_{L^p(\RR^n)}\leq C \left\|m\right\|_{L^p(\RR^n)}\leq C
l(Q)^{n\left(\frac 1p-1\right)}, \label{eq7.25}
\end{eqnarray}

\noindent as desired.

Collecting all the terms, we arrive at (\ref{eq7.15}) with
$\gamma=\min\{\eps_1,\epsilon\}$. \ep

\begin{corollary}\label{c7.5}
$H^1_L(\RR^n)=H^1_{{\cal N}_P}(\RR^n)$, in particular,
$\|f\|_{H^1_L(\RR^n)}\approx\|{\cal N}_P f\|_{L^1(\RR^n)}+\| f\|_{L^1(\RR^n)}$.
\end{corollary}

\bp The Corollary follows from Theorems~\ref{t7.1}, \ref{t7.4} and
\ref{t5.6}.\ep

\section{$BMO_L(\RR^n)$:\,duality with Hardy spaces.}
\setcounter{equation}{0}

We start with an auxiliary lemma that gives an equivalent
characterization of $BMO_L(\RR^n)$ using
the resolvent in place of the heat semigroup. In the sequel we shall
frequently use the characterization below as the definition of
$BMO_L(\RR^n)$ without additional comments.  In addition, by the results of Section 4, we
are at liberty to choose the molecular parameters $\eps > 0$ and $M > n/4$ at our convenience.
In the sequel, we shall use this fact without further comment.

\begin{lemma}\label{l8.1}
An element $f \in \cap_{\eps > 0}({\bf M}^{2,\eps,M}_0(L^*))^*\equiv ({\bf M}^{2,M}_0(L^*))^*$ belongs to $BMO_L(\RR^n)$ if and only if

\begin{equation}\label{eq8.1}
\sup_{Q\subset\RR^n}\left(\frac{1}{|Q|}\int_Q
\left|(I-(1+l(Q)^2L)^{-1})^Mf(x)\right|^2\,dx\right)^{1/2}<\infty,
\end{equation}

\noindent where $M>n/4$ and $Q$ stands for a cube in $\RR^n$.
\end{lemma}

\bp For brevity in this proof we shall distinguish (\ref{eq8.1})
as $\|f\|_{BMO_L^{res}(\RR^n)}$. In the rest of the paper both the
norm based on the heat semigroup and the one based on the
resolvent will be denoted by $\|f\|_{BMO_L(\RR^n)}$.

\noindent {\bf Step I}. Let us start with the ``$\leq$"
inequality. To this end, we split

\begin{equation}\label{eq8.2}
f=(I-(1+l(Q)^2L)^{-1})^Mf+\left[I-(I-(1+l(Q)^2L)^{-1})^M\right]f.
\end{equation}

\noindent For every $Q\subset\RR^n$

\begin{eqnarray}
&&\hskip -1cm \left(\frac{1}{|Q|}\int_Q
\left|(I-e^{-l(Q)^2L})^M(I-(1+l(Q)^2L)^{-1})^Mf(x)\right|^2\,dx\right)^{1/2}
\nonumber\\[4pt]
&&\hskip -1cm\qquad \leq C \sum_{k=0}^M\sum_{j=0}^\infty
\left(\frac{1}{|Q|}\int_Q \left|e^{-kl(Q)^2L}\bigl[\chi_{S_j(Q)}
(I-(1+l(Q)^2L)^{-1})^Mf\bigr](x)\right|^2\,dx\right)^{1/2}
\nonumber\\[4pt]
&&\hskip -1cm\qquad \leq C
\|f\|_{BMO_L^{res}(\RR^n)}+\nonumber\\[4pt]
&&\hskip -1cm\qquad\qquad +C\sum_{j=2}^\infty
e^{-\frac{(2^jl(Q))^2}{c l(Q)^2}} \left(\frac{1}{|Q|}\int_{S_j(Q)}
\left|(I-(1+l(Q)^2L)^{-1})^Mf(x)\right|^2\,dx\right)^{1/2},\label{eq8.3}
\end{eqnarray}

\noindent where we used Lemmas~\ref{l2.3}, \ref{l2.4} for the second
inequality. Now one can cover $S_j(Q)$ by approximately $2^{jn}$
cubes of the sidelength $l(Q)$, this allows to bound the second term
in the expression above by

\begin{equation}\label{eq8.4}
C\sum_{j=2}^\infty e^{-c\,2^{2j}}2^{jn/2}\,
\|f\|_{BMO_L^{res}(\RR^n)}\leq C \|f\|_{BMO_L^{res}(\RR^n)},
\end{equation}

\noindent as desired.

As for the remaining term, observe that

\begin{eqnarray}\nonumber
&&\hskip
-1cm\frac{I-(I-(1+l(Q)^2L)^{-1})^M}{(I-(1+l(Q)^2L)^{-1})^{M}}=(I-(1+l(Q)^2L)^{-1})^{-M}-I
=\left(\frac{1+l(Q)^2L}{l(Q)^2L}\right)^M-I\\[4pt]
&&\hskip -1cm\qquad =
\left(1+(l(Q)^2L)^{-1}\right)^M-I=\sum_{k=1}^M\frac{{M}!}{({M}-k)!\,k!}\,(l(Q)^2L)^{-k},\label{eq8.5}
\end{eqnarray}

\noindent and therefore,

\begin{eqnarray}
&&\hskip -0.7cm \left(\frac{1}{|Q|}\int_Q
\left|(I-e^{-l(Q)^2L})^M[I-(I-(1+l(Q)^2L)^{-1})^M]f(x)\right|^2\,dx\right)^{1/2}
\nonumber\\[4pt]
&&\hskip -0.7cm \qquad \leq C \sum_{k=1}^M \Bigl(\frac{1}{|Q|}\int_Q
\Bigl| (I-e^{-l(Q)^2L})^{M-k}\Bigl(-
\int_0^{l(Q)}\partial_\tau e^{-\tau^2L}\,d\tau \Bigr)^{k} (l(Q)^2L)^{-k} \nonumber\\[4pt]
&&\hskip -0.7cm \qquad\qquad\qquad\qquad\qquad \times
(I-(1+l(Q)^2L)^{-1})^Mf(x)\Bigr|^2\,dx\Bigr)^{1/2}
\nonumber\\[4pt]
&&\hskip -0.7cm \qquad \leq C \sum_{k=1}^M \Bigl(\frac{1}{|Q|}\int_Q
\Bigl| (I-e^{-l(Q)^2L})^{M-k}\Bigl(
\int_0^{l(Q)}\frac{\tau}{l(Q)^2}\, e^{-\tau^2L}\,d\tau \Bigr)^{k}\nonumber\\[4pt]
&&\hskip -0.7cm \qquad\qquad\qquad\qquad\qquad \times
(I-(1+l(Q)^2L)^{-1})^Mf(x)\Bigr|^2\,dx\Bigr)^{1/2}. \label{eq8.6}
\end{eqnarray}

\noindent Having this at hand, we obtain the required estimate
changing the order of integration above and using the annular
decomposition and Gaffney estimates, much as in
(\ref{eq8.3})--(\ref{eq8.4}).

\noindent {\bf Step II}. Let us now consider the ``$\geq$" part of
(\ref{eq8.1}). For every $x\in\RR^n$

\begin{eqnarray}
&&\hskip -1cm f(x)= 2^M\left(l(Q)^{-2}\int_{l(Q)}^{\sqrt 2
l(Q)}s\,ds\right)^Mf(x)
\nonumber\\[4pt]
&&\hskip -1cm \quad = 2^M\, l(Q)^{-2}\int_{l(Q)}^{\sqrt
2
l(Q)}s_1(I-e^{-s_1^2L})^M\,ds_1\,\left(l(Q)^{-2}\int_{l(Q)}^{\sqrt
2 l(Q)}s\,ds\right)^{M-1} f(x)
\nonumber\\[4pt]
&&\hskip -1cm \qquad + \sum_{k=1}^M C_{k,M}
l(Q)^{-2}\int_{l(Q)}^{\sqrt 2 l(Q)}s_1
e^{-ks_1^2L}\,ds_1\,\left(l(Q)^{-2}\int_{l(Q)}^{\sqrt 2
l(Q)}s\,ds\right)^{M-1} f(x),\label{eq8.7}
\end{eqnarray}

\noindent where $C_{k,M}\in\RR$ are some constants depending on
$k$ and $M$ only. However, $\partial_s e^{-ks^2L}=-2kLs
e^{-ks^2L}$ and therefore,

\begin{eqnarray}
&& 2kL\int_{l(Q)}^{\sqrt 2 l(Q)}s e^{-ks^2L}\,ds =
e^{-kl(Q)^2L}-e^{-2kl(Q)^2L}=e^{-kl(Q)^2L}(I-e^{-kl(Q)^2L})\nonumber\\[4pt]
&&\qquad = e^{-kl(Q)^2L}(I-e^{-l(Q)^2L})\sum_{i=0}^{k-1}
e^{-il(Q)^2L}.\label{eq8.8}
\end{eqnarray}

\noindent Applying the procedure outlined in
(\ref{eq8.7})--(\ref{eq8.8}) $M$ times, we arrive at the following
formula

\begin{equation}\label{eq8.9}
f(x)=\sum_{i=1}^{(M+1)^M}\,l(Q)^{-2M}L^{-N_i} \prod_{k=1}^M
C(i,k,M)\,p_{i,k}f(x),
\end{equation}

\noindent where $0\leq N_i\leq M$ and for all $i,k$ as above either

\begin{equation}\label{eq8.10}
p_{i,k}= \int_{l(Q)}^{\sqrt 2 l(Q)}s(I-e^{-s^2L})^M\,ds
\end{equation}

\noindent or $p_{i,k}$ is of the form (\ref{eq8.8}).

Fix some $Q\subset\RR^n$ and $x\in Q$ and consider
$(I-(1+l(Q)^2L)^{-1})^Mf(x)$ with $f$ represented in the form
(\ref{eq8.9}). The negative powers of $L$ can be handled writing

\begin{equation}\label{eq8.11}
(I-(1+l(Q)^2L)^{-1})^{N_i}l(Q)^{-2N_i}L^{-N_i}=(1+l(Q)^2L)^{-N_i}.
\end{equation}

\noindent Then the new expression for $(I-(1+l(Q)^2L)^{-1})^Mf(x)$
is a linear combination of terms, with the property that each term
contains

\begin{equation}\label{eq8.12}
\mbox{ either }\quad l(Q)^{-2}\int_{l(Q)}^{\sqrt 2
l(Q)}s(I-e^{-s^2L})^M\,ds\quad\mbox{ or }\quad (I-e^{-l(Q)^2L})^M,
\end{equation}

\noindent and a finite number of factors (almost) in the form of
resolvent or heat semigroup corresponding to $t\approx l(Q)$. One
can now build an argument similar to Step I,
(\ref{eq8.3})--(\ref{eq8.4}), using dyadic annular decomposition
and Gaffney estimates, to single out $(I-e^{-l(Q)^2L})^M$ or
$(I-e^{-s^2L})^M$, $s\approx l(Q)$, and obtain the desired
estimate. We leave the details to the interested reader. \ep

\begin{theorem}\label{t8.2}
Let $f\in BMO_{L^*}(\RR^n)$ for some $M\in\NN$. Then the linear
functional given by

\begin{equation}\label{eq8.13}
l(g)=\langle f,g \rangle ,
\end{equation}

\noindent initially defined on the dense subspace of finite linear
combinations of $(2,\eps, M)$-molecules, $\eps>0$, via the pairing of ${\bf M}^{2,\eps,M}_0$ with its dual,
has a unique bounded extention to $H^1_L(\RR^n)$ with

\begin{equation}\label{eq8.14}
\|l\|\leq C\|f\|_{BMO_{L^*}(\RR^n)}.
\end{equation}
\end{theorem}

\bp Let us prove first that for every $(2,\eps, M)$-molecule $m$

\begin{equation}\label{eq8.15}
\left|\langle f,m \rangle \right|\leq C\|f\|_{BMO_{L^*}(\RR^n)}.
\end{equation}

\noindent By definition, $f \in ({\bf M}^{2,M}_0(L))^*,$ so in particular $(I-(1+l(Q)^2L^*)^{-1})^Mf \in
L^2_{loc}$ (see the discussion preceding (\ref{eq1.def})).  Thus, we may write

\begin{eqnarray}\nonumber
\langle f,m \rangle
&=&\int_{\RR^n}(I-(1+l(Q)^2L^*)^{-1})^Mf(x)\overline{m(x)}\,dx\\[4pt]\nonumber
&&\qquad\qquad\qquad+\,\left\langle \Bigl[I-(I-(1+l(Q)^2L^*)^{-1})^M\Bigr]f,m \right\rangle\\[4pt]
&&\qquad =: I_1+I_2\label{eq8.16},
\end{eqnarray}

\noindent where $Q$ is the cube associated to $m$. Then

\begin{eqnarray}\nonumber
|I_1|&\leq & \sum_{j=0}^\infty
\left(\int_{S_j(Q)}\left|(I-(1+l(Q)^2L^*)^{-1})^Mf(x)\right|^2\,dx\right)^{1/2}
\left(\int_{S_j(Q)}|m(x)|^2\,dx\right)^{1/2} \\[4pt]\nonumber
&\leq & \sum_{j=0}^\infty 2^{-j\eps}
\left(\frac{1}{(2^jl(Q))^n}\int_{S_j(Q)}
\left|(I-(1+l(Q)^2L^*)^{-1})^Mf(x)\right|^2\,dx\right)^{1/2}\\[4pt]
& \leq & C\|f\|_{BMO_{L^*}(\RR^n)}\label{eq8.17},
\end{eqnarray}

\noindent where we used (\ref{eq1.8}) for the second inequality,
and the third one follows by covering $S_j(Q)$ by $C2^{jn}$ cubes
of the sidelength $l(Q)$.

To analyze $I_2$ recall (\ref{eq8.5}) (with $L^*$ in place of $L$), and write

\begin{eqnarray}\nonumber
|I_2|&\leq & C \sum_{k=1}^M\left|
\int_{\RR^n}(I-(1+l(Q)^2L^*)^{-1})^Mf(x)\,\overline{(l(Q)^2L)^{-k}m(x)}\,dx\right|\\[4pt]
&\leq & C\sum_{k=1}^M\sum_{j=0}^\infty
\left(\int_{S_j(Q)}\left|(I-(1+l(Q)^2L^*)^{-1})^Mf(x)\right|^2\,dx\right)^{1/2}\times
\nonumber\\[4pt]
&&\qquad\qquad
\times\left(\int_{S_j(Q)}|(l(Q)^2L)^{-k}m(x)|^2\,dx\right)^{1/2}.\label{eq8.18}
\end{eqnarray}

\noindent We finish as in (\ref{eq8.17}) using (\ref{eq1.9}).  Thus, (\ref{eq8.15})
is now established.

Having at hand (\ref{eq8.15}), our goal is to show that for every $N\in\NN$
and for every  $g=\sum_{j=0}^N \lambda_j m_j$, where
$\{m_j\}_{j=0}^N$ are $(2,\eps, M')$-molecules, and $M'>n/4$ is chosen large enough relative to $M$,
we have

\begin{equation}\label{eq8.19}
\left|\langle f,g \rangle \right| \leq C\|g\|_{H^1_L(\RR^n)}
\|f\|_{BMO_{L^*}(\RR^n)}.
\end{equation}

\noindent Since the space of finite linear combinations of
$(2,\eps,M')$-molecules is dense in $H^1_L(\RR^n)$, the linear
functional $l$ will then have a unique bounded extension to $H^1_L(\RR^n)$
defined in a standard fashion by continuity.  We point out that this extension by continuity
depends on having a bound in terms of $\|g\|_{H^1_L(\RR^n)}$ in (\ref{eq8.19}), as opposed to  $\sum_{j=0}^N
|\lambda_j|$. The latter bound is immediately obtainable from (\ref{eq8.15}) (since in particular, a $(2,\eps,M')$-molecule is a $(2,\eps,M)$-molecule whenever $M'\geq M$), but may be much larger
than the $H^1_L$ norm.  To obtain the sharper bound (\ref{eq8.19}) will be somewhat delicate.
In the classical setting, the same issue arises, but may be handled in a fairly routine fashion
by truncating the BMO function so that it may be approximated in $(H^1)^*$ by bounded functions
(see, e.g. \cite{St2}, pp. 142-143).  This avenue is not available to us, as we cannot expect that any
$L^\infty$ truncation will interact well with our operator $L$.  Instead, we should seek a ``truncation"
in $L^p, p\in(p_L,\widetilde{p}_L).$  In fact, approximating by $L^2$ functions will be most convenient, and this is what we shall do.
We note that it is to deal with this difficulty that we have been forced to introduce the equivalent
norm $\|\cdot\|_{\widetilde{H}^1_L}$. The reason for our doing so will become apparent in the sequel.

We shall require some rather extensive preliminaries.
In particular, we shall use the ``tent space" approach of Coifman-Meyer-Stein \cite{CMS}.  Let us now recall the basic theory.

For some $F:\RR^{n+1}_+\to \CC$,
$\RR^{n+1}_+=\RR^n\times(0,\infty)$, consider the square function
$SF:=S^1F$, where $S^\alpha F$, $\alpha>0$, was defined in
(\ref{eq2.6}) and

\begin{equation}\label{eq8.22}
CF(x):=\sup_{B:\,x\in B}\left(\frac{1}{|B|}\dint_{\widehat
B}|F(y,t)|^2\,\frac{dydt}{t}\right)^{1/2}, \qquad x\in\RR^n,
\end{equation}

\noindent where $B$ stands for a ball in $\RR^n$ and

\begin{equation}\label{eq8.23}
\widehat B:=\{(x,t)\in\RR^n\times(0, \infty):\,{\rm
dist}(x,\,^cB)\geq t\},
\end{equation}

\noindent is the tent region above ball $B$. Define the tent spaces

\begin{equation}\label{eq8.24}
T^1(\RR^{n+1}_+):=\{F:\RR^{n+1}_+\longrightarrow
\CC;\,\|F\|_{T^1(\RR^{n+1}_+)}:=\|SF\|_{L^1(\RR^n)}<\infty\},
\end{equation}

\noindent and

\begin{equation}\label{eq8.25}
T^\infty(\RR^{n+1}_+):=\{F:\RR^{n+1}_+\longrightarrow
\CC;\,\|F\|_{T^\infty(\RR^{n+1}_+)}:=\|CF\|_{L^\infty(\RR^n)}<\infty\},
\end{equation}

\noindent and recall from \cite{CMS} that
$(T^1(\RR^{n+1}_+))^*=T^\infty(\RR^{n+1}_+)$.

We now prove the following analogue of a classical estimate of \cite{FeSt}.
\begin{lemma}\label{l8.3}  The operator $$f \mapsto C\left((t^2L)^Me^{-t^2L}f\right)$$
maps $BMO_L(\RR^n) \to T^\infty (\RR_+^{n+1})$;  i.e.,
$$\sup_{B}\left(\frac{1}{|B|}\dint_{\widehat
B}|(t^2L)^Me^{-t^2L}f(x)|^2\,\frac{dxdt}{t}\right)^{1/2}\leq C
\|f\|_{BMO_L(\RR^n)}.$$
\end{lemma}

\bp For every cube $Q\subset\RR^n$

\begin{eqnarray}
&&\left(\frac{1}{|Q|}\int_0^{l(Q)}\int_{Q}
\left|(t^2L)^Me^{-t^2L}f(y)\right|^2\,\frac{dydt}{t}\right)^{1/2}\nonumber\\[4pt]
&& \qquad \quad = \left(\frac{1}{|Q|}\int_0^{l(Q)}\int_{Q}
\left|(t^2L)^Me^{-t^2L}(I-(1+l(Q)^2L)^{-1})^M
f(y)\right|^2\,\frac{dydt}{t}\right)^{1/2}+\nonumber\\[4pt]
&& \qquad \qquad +\left(\frac{1}{|Q|}\int_0^{l(Q)}\int_{Q}
\left|(t^2L)^Me^{-t^2L}\bigl[I-(I-(1+l(Q)^2L)^{-1})^M\bigr]
f(y)\right|^2\,\frac{dydt}{t}\right)^{1/2}\nonumber\\[4pt]
&& \qquad \quad =:I_1+I_2.\label{eq9.7}
\end{eqnarray}

\noindent Then

\begin{eqnarray}
I_1 &\leq &\sum_{j=0}^\infty
\left(\frac{1}{|Q|}\int_0^{l(Q)}\int_{Q}
\left|(t^2L)^Me^{-t^2L}\bigl[\chi_{S_j(Q)}(I-(1+l(Q)^2L)^{-1})^M
f\bigr](y)\right|^2\,\frac{dydt}{t}\right)^{1/2}\nonumber\\[4pt]
&\leq & \frac{1}{|Q|^{1/2}}\,\left\|g_h^M\bigl(\chi_{2Q}
(I-(1+l(Q)^2L)^{-1})^M f\bigr)\right\|_{L^2(\RR^n)}+\nonumber\\[4pt]
&& \, +\sum_{j=2}^\infty
\frac{C}{|Q|^{1/2}}\left(\int_0^{l(Q)}\,e^{-\frac{(2^jl(Q))^2}{ct^2}}\,\frac{dt}{t}\right)^{1/2}
\, \|(I-(1+l(Q)^2L)^{-1})^M f\|_{L^2(S_j(Q))},\label{eq9.8}
\end{eqnarray}

\noindent where

\begin{equation}\label{eq9.9}
g_h^Mf(x):=\left(\int_0^\infty|(t^2L)^{M}e^{-t^2 L}f(x)|^2\, \frac
{dt}{t}\right)^{1/2}, \qquad M\in\NN,\quad x\in\RR^n,
\end{equation}

\noindent is bounded in $L^2(\RR^n)$ according to \cite{ADM}.
Therefore, for every $N\in\NN$

\begin{eqnarray}
I_1 &\leq & \frac{C}{|Q|^{1/2}}\,\left\|(I-(1+l(Q)^2L)^{-1})^M f\right\|_{L^2(2Q)}+\nonumber\\[4pt]
&& \qquad +\sum_{j=2}^\infty 2^{-jN} \frac{C}{(2^{j}l(Q))^{n/2}}\,
\|(I-(1+l(Q)^2L)^{-1})^M f\|_{L^2(S_j(Q))} \nonumber\\[4pt]
&\leq &C\|f\|_{BMO_L(\RR^n)}.\label{eq9.10}
\end{eqnarray}

\noindent To estimate $I_2$ we use (\ref{eq8.5}) and write

\begin{eqnarray*}
&& \hskip -0.7 cm I_2 \leq C\sup_{1\leq k\leq
M}\left(\frac{1}{|Q|}\int_0^{l(Q)}\int_{Q}
\left|(t^2L)^Me^{-t^2L}(l(Q)^2L)^{-k}(I-(1+l(Q)^2L)^{-1})^M
f(y)\right|^2\,\frac{dydt}{t}\right)^{1/2}\nonumber\\[4pt]
&& \hskip -0.7 cm \,\,\leq C\sup_{1\leq k\leq
M}\left(\frac{1}{|Q|}\int_0^{l(Q)}\left(\frac{t}{l(Q)}\right)^{2k}\int_{Q}
\left|(t^2L)^{M-k}e^{-t^2L}(I-(1+l(Q)^2L)^{-1})^M
f(y)\right|^2\,\frac{dydt}{t}\right)^{1/2},
\end{eqnarray*}

\noindent the rest of the argument is similar to
(\ref{eq9.8})--(\ref{eq9.10}). This finishes the proof of
Lemma~\ref{l8.3}. \ep

We shall also require
\begin{lemma}\label{l8.4} Suppose that $f\in ({\bf M}^{2,M}_0(L))^*$ satisfies the ``controlled growth estimate"
\begin{equation}\label{eq9.5}
\int_{\RR^n}\frac{|(I-(1+L^*)^{-1})^Mf(x)|^2}{1+|x|^{n+\eps_1}}\,dx<\infty ,
\end{equation}
for some $\eps_1 > 0$ (in particular, this holds trivially for every $\eps_1 > 0$ if $f\in BMO_{L^*}$).
Then for every $g \in H^1_{L}$ that can be represented as a
finite linear combination of $(2,\eps,M')$ molecules, with $\eps,M'$ sufficiently large compared to
$\eps_1,M$, we have
\begin{equation}\label{eq9.11}
\langle f,g \rangle =C_M \dint_{\RR^{n+1}_+}
(t^2L^*)^Me^{-t^2L^*}f(x)\,\overline{t^2Le^{-t^2L}g(x)}\,\frac{dtdx}{t}.
\end{equation}
\end{lemma}

\bp For $\delta,R>0$ consider

\begin{eqnarray}
&& \hskip -.7 cm\int_{\RR^n}\int_\delta^R
(t^2L^*)^Me^{-t^2L^*}f(x)\,\overline{t^2Le^{-t^2L}g(x)}\,\frac{dtdx}{t}
=\left\langle f,\left(\int_\delta^R
\bigl(t^2L\bigr)^{M+1}e^{-2t^2L}g\,\frac{dt}{t}\right)\right\rangle
\nonumber\\[4pt]
&&\quad = C_M^{-1}\,\langle f,g \rangle
\,-\,\left\langle f,\left(C_M^{-1}\,g -\int_\delta^R
\bigl(t^2L\bigr)^{M+1}e^{-2t^2L}g\,\frac{dt}{t}\right)\,\right\rangle .
 \label{eq9.12}
\end{eqnarray}
We will now write $f$ in the following way

\begin{eqnarray}\nonumber
f&=&(I-(1+L^*)^{-1}+(1+L^*)^{-1})^Mf\\ \nonumber &=&\sum_{k=0}^M\frac{{M}!}{({M}-k)!\,k!}
\,(I-(1+L^*)^{-1})^{M-k} (1+L^*)^{-k}f\\
&=&\sum_{k=0}^M\frac{{M}!}{({M}-k)!\,k!} \,(L^*)^{-k}(I-(1+L^*)^{-1})^{M}f
\label{eq9.13}
\end{eqnarray}
Thus, the last expression in (\ref{eq9.12}) equals $ \sum_{k=0}^M C_{k,M}$ times
\begin{eqnarray}\nonumber
\left\langle (I-(1+L^*)^{-1})^M f,\left( C_M^{-1} L^{-k} g - \int_\delta^R
\bigl(t^2L\bigr)^{M+1}e^{-2t^2L}L^{-k} g\,\frac{dt}{t}\right)\,\right\rangle \,=
&& \\ \nonumber \left\langle (I-(1+L^*)^{-1})^M f, \int_0^\delta
\bigl(t^2L\bigr)^{M+1}e^{-2t^2L}L^{-k} g\,\frac{dt}{t}\,\right\rangle \qquad \qquad &\qquad \qquad & \qquad \\ \qquad \qquad \qquad \,+\,\left\langle (I-(1+L^*)^{-1})^M f, \int_R^\infty \bigl(t^2L\bigr)^{M+1}e^{-2t^2L}L^{-k}g\,\frac{dt}{t}\,\right\rangle,&& \label{eq8.split}
\end{eqnarray}
since for $L^{-k} g\in L^2(\RR^n)$ the Calder\'{o}n reproducing
formula is valid.  The last term in (\ref{eq8.split}) is bounded by a constant times
\begin{eqnarray}
&&\hskip -.3cm
\left(\int_{\RR^n}\frac{|(I-(1+L^*)^{-1})^{M}f(x)|^2} {1+|x|^{n+\eps_1}}\,dx\right)^{1/2}\nonumber\\[4pt]
&&\hskip -.3cm \qquad \qquad \times \sup_{0\leq k \leq M}
\left(\int_{\RR^n}\Bigl|\int_R^\infty
\bigl(t^2L\bigr)^{M+1}e^{-2t^2L}L^{-k}g(x)\,\frac{dt}{t}\Bigr|^2
(1+|x|^{n+\eps_1})\,dx\right)^{1/2}\nonumber\\[4pt]
&& \hskip -.3cm \leq C \Upsilon \sup_{0\leq k \leq M}\sum_{j=0}^\infty
2^{j(n+\eps_1)/2} \left(\int_{S_j(Q_0)}\Bigl|\int_R^\infty
\bigl(t^2L\bigr)^{M+1}e^{-2t^2L}L^{-k}g(x)\,\frac{dt}{t}\Bigr|^2
\,dx\right)^{1/2}\nonumber\\[4pt]
&&\hskip -.3cm \leq C \Upsilon \sup_{0\leq k \leq M}\sum_{j=0}^\infty
2^{j(n+\eps_1)/2} \int_R^\infty \Bigl\|
\bigl(t^2L\bigr)^{M+M'-k+1}e^{-2t^2L}
L^{-M'}g\Bigr\|_{L^2(S_j(Q_0))}\,\frac{dt}{t^{2(M'-k)+1}},\nonumber\\
\label{eq9.14}
\end{eqnarray}

\noindent where $\Upsilon$ is the finite quantity defined in (\ref{eq9.5}), and $Q_0$ is the cube centered at $0$ with the
sidelength 1. Then the expression under the $\sup$ sign above is
bounded modulo multiplicative constant by

\begin{eqnarray}
&&\hskip -.7 cm \int_R^\infty \Bigl\|
\bigl(t^2L\bigr)^{M+M'-k+1}e^{-2t^2L}
L^{-M'}g\Bigr\|_{L^2(4Q_0)}\,\frac{dt}{t^{2(M'-k)+1}} \nonumber\\[4pt]
&&\hskip -.7 cm  \quad + \sum_{j=3}^\infty 2^{j(n+\eps_1)/2}
\int_R^\infty \Bigl\| \bigl(t^2L\bigr)^{M+M'-k+1}e^{-2t^2L}
\left[\chi_{\RR^n\setminus
2^{j-2}Q_0}L^{-M'}g\right]
\Bigr\|_{L^2(S_j(Q_0))}\,\frac{dt}{t^{2(M'-k)+1}} \nonumber\\[4pt]
&&\hskip -.7 cm \quad + \sum_{j=3}^\infty 2^{j(n+\eps_1)/2}
\int_R^\infty \Bigl\| \bigl(t^2L\bigr)^{M+M'-k+1}e^{-2t^2L}
\left[\chi_{2^{j-2}Q_0}L^{-M'}g\right]
\Bigr\|_{L^2(S_j(Q_0))}\,\frac{dt}{t^{2(M'-k)+1}}\nonumber\\[4pt]
&&\hskip -.7 cm \,\, \leq
\frac{C}{R^{2(M'-k)}}\,\Bigl\|L^{-M'}g\,\Bigr\|_{L^2(\RR^n)}
+\frac{C}{R^{2(M'-k)}}\sum_{j=3}^\infty 2^{j(n+\eps_1)/2}
\Bigl\|L^{-M'}g\,\Bigr\|_{L^2(\RR^n\setminus
2^{j-2}Q_0)}\nonumber\\[4pt]
&&\hskip -.7 cm \quad + \sum_{j=3}^\infty 2^{j(n+\eps_1)/2}
\left(\int_R^\infty
e^{-\frac{2^{2j}}{ct^2}}\,\frac{dt}{t^{2(M'-k)+1}}\right)\,
\left\|L^{-M'}g\right\|_{L^2(\RR^n)}. \label{eq9.15}
\end{eqnarray}

\noindent However,

\begin{equation}\label{eq9.16}
\int_R^\infty e^{-\frac{2^{2j}}{ct^2}}\,\frac{dt}{t^{2(M'-k)+1}}\leq
C \frac{1}{2^{2(M'-k)j}} \int_{R/2^j}^\infty
e^{-\frac{1}{s^2}}\,\frac{ds}{s^{2(M'-k)+1}}\leq C
\frac{1}{2^{2(M'-k)j}} (2^j/R)^{\epsilon'}
\end{equation}

\noindent for every $\epsilon'>0$. Also, $g$ is a finite linear
combination of $(2,\eps,M')$- molecules, therefore for large $j$

\begin{equation}\label{eq9.17}
\Bigl\|L^{-k'}g\,\Bigr\|_{L^2(S^j(Q_0))}\leq C
2^{-j(n/2+\eps)},\qquad 0\leq k'\leq M',
\end{equation}

\noindent which allows to estimate the second term in
(\ref{eq9.15}). Without loss of generality we can assume that
$\eps>\eps_1/2$ and $M'>\frac{n+\eps_1}{4}+M$. Then there exists
$\epsilon_0>0$ such that the quantity in (\ref{eq9.15}), and hence
the one in (\ref{eq9.14}), does not exceed $C/R^{\epsilon_0}$.

We now turn to the integral over $(0,\delta)$.  For convenience of notation, we set
$$\tilde{f} \equiv (I-(1+L^*)^{-1})^Mf .$$ Since
$-2tLe^{-t^2L}=\partial_t e^{-t^2L}$, we may write write

\begin{eqnarray}
&& \left|\left\langle \tilde{f},\int_0^\delta
\bigl(t^2L\bigr)^{M+1}e^{-2t^2L}L^{-k}g(x)\,\frac{dt}{t}\right\rangle\right|
=C\left|\left\langle \tilde{f},\int_0^\delta
\bigl(t^2L\bigr)^{M}\partial_t e^{-2t^2L}L^{-k}g(x)\,dt \right\rangle \right|\nonumber\\[4pt]
&& \quad \leq C \left|\left\langle \tilde{f},\int_0^\delta
\bigl(t^2L\bigr)^{M}e^{-2t^2L}L^{-k}g(x)\,\frac{dt}{t} \right\rangle \right|
+C\left|\left\langle \tilde{f},
\bigl(\delta^2L\bigr)^{M}e^{-2\delta^2L}L^{-k}g(x)\right\rangle\right|\nonumber\\[4pt]
\nonumber\\[4pt]
&& \quad \leq C\sum_{k=1}^M\left|\left\langle \tilde{f},
\bigl(\delta^2L\bigr)^{k}e^{-2\delta^2L}L^{-k}g(x)\right\rangle \right|
+C\left|\left\langle \tilde{f},
\bigl(e^{-2\delta^2L}-I\bigr)L^{-k}g(x)\right\rangle \right|,\nonumber
\end{eqnarray}

\noindent repeatedly integrating by parts in $t$. Therefore, as in (\ref{eq9.14}),

\begin{eqnarray}
 && \left|\left\langle \tilde{f},\int_0^\delta
\bigl(t^2L\bigr)^{M+1}e^{-2t^2L}L^{-k}g(x)\,\frac{dt}{t}\right\rangle \right|\nonumber\\[4pt]
&&\qquad \leq C \Upsilon \sup_{0\leq k'\leq M}\sum_{k=1}^M \sum_{j=0}^\infty 2^{j(n+\eps_1)/2}
\Bigl\| \bigl(\delta^2L\bigr)^{k}e^{-2\delta^2L}
L^{-k'}g\Bigr\|_{L^2(S_j(Q_0))}\nonumber\\[4pt]
&&\qquad\qquad \qquad +C \Upsilon \sup_{0\leq k'\leq M}\sum_{j=0}^\infty
2^{j(n+\eps_1)/2} \Bigl\| \bigl(e^{-2\delta^2L}-I\bigr)
L^{-k'}g\Bigr\|_{L^2(S_j(Q_0))}. \label{eq9.18}
\end{eqnarray}

\noindent Now let us split $L^{-k'}g=\chi_{R_j}L^{-k'}g+\chi_{\,^cR_j}L^{-k'}g$ where

\begin{eqnarray*}
R_j&=&2^{j+2}Q_0, \quad \mbox{if} \quad j=0,1,2,\\[4pt]
R_j&=&2^{j+2}Q_0\setminus 2^{j-2}Q_0, \quad \mbox{if} \quad
j=3,4,...,
\end{eqnarray*}

\noindent and start with the part of (\ref{eq9.18}) corresponding
to $\chi_{R_j}L^{-k'}g$. Fix some $\eta>0$. Then for $N\in\NN$
and for all $0\leq k'\leq M$

\begin{eqnarray}
&&C \sum_{k=1}^M \sum_{j=N}^\infty 2^{j(n+\eps_1)/2}\Bigl( \,\Bigl\|
\bigl(\delta^2L\bigr)^{k}e^{-2\delta^2L}
(\chi_{R_j}L^{-k'}g)\Bigr\|_{L^2(S_j(Q_0))}\\[4pt]
&&\qquad\qquad\qquad\qquad +\,\,\Bigl\|
\bigl(e^{-2\delta^2L}-I\bigr)
(\chi_{R_j}L^{-k'}g)\Bigr\|_{L^2(S_j(Q_0))}\,\Bigr)\nonumber\\[4pt]
&&\qquad \leq C \sum_{j=N}^\infty 2^{j(n+\eps_1)/2}
\|L^{-k'}g\|_{L^2(R_j)}\leq C \sum_{j=N}^\infty
2^{j(n+\eps_1)/2} 2^{-j(n/2+\eps)}, \label{eq9.20}
\end{eqnarray}

\noindent where the last inequality uses (\ref{eq9.17}). Recall that
$\eps>\eps_1/2$. Then choosing $N\approx -\ln \eta$, we can control
the expression above by $\eta$. As for the remaining part, for
$\delta$ small enough

\begin{eqnarray}
&&C \sum_{k=1}^M \sum_{j=0}^N 2^{j(n+\eps_1)/2}\Bigl( \Bigl\|
\bigl(\delta^2L\bigr)^{k}e^{-2\delta^2L}
(\chi_{R_j}L^{-k'}g)\Bigr\|_{L^2(S_j(Q_0))}\nonumber\\[4pt]
&&\qquad\qquad + \Bigl\| \bigl(e^{-2\delta^2L}-I\bigr)
(\chi_{R_j}L^{-k'}g)\Bigr\|_{L^2(S_j(Q_0))}\Bigr)\leq C \eta
\label{eq9.21}
\end{eqnarray}

\noindent using that
$\bigl(\delta^2L\bigr)^{k}e^{-2\delta^2L}\to 0$ and
$e^{-2\delta^2L}-I\to 0$ in the strong operator topology as
$\delta\to 0$.

The integral corresponding to $\chi_{\,^cR_j}L^{-k'}g$ is
analyzed similarly, with the only difference that the Gaffney
estimates instead of $L^2$-decay of $L^{-k'}g$ are used to
control an analogue of (\ref{eq9.20}).

We have proved that the second term in (\ref{eq9.12}) vanishes as
$\delta \to 0$ and $R\to\infty$. Therefore, the formula
(\ref{eq9.11}) is justified for every $g$ belonging to the space of
finite linear combinations of molecules.

\ep

We return now to the proof of (\ref{eq8.19}).
We shall approximate $f$ by $$f_K \equiv \int_{1/K}^K t^2L^* e^{-t^2L^*}
\left(\chi_{B_K} (t^2L^*)^M e^{-t^2L^*} f \right) \frac{dt}{t},$$
where $B_K \equiv \{x \in \RR^n: |x| < K \}$.
We claim that $f_K \in L^2$, and that
\begin{equation}\label{eq8.b} \sup_K \|f_K\|_{BMO_{L^*}(\RR^n)} \leq
C \|(t^2L^*)^M e^{-t^2L^*} f\|_{T^\infty(\RR^{n+1}_+)} \leq C
\|f\|_{BMO_{L^*}(\RR^n)}.
\end{equation}
We note that the second inequality in (\ref{eq8.b}) is just Lemma \ref{l8.3},
so the key issue is the first inequality.  Let us take the claim for granted momentarily.
Since $g$ is a finite linear
combination of $(2,\eps,M')$-molecules, in particular we have that $g \in \widehat{H}^1_L$.
Consequently, there is a $\delta>0$ and a $\delta$-representation
$g=\sum \lambda_i m_i$, converging in $L^2$, with
$\sum |\lambda_i| \approx \|g\|_{\widehat{H}^1_L} \approx \|g\|_{H^1_L}$ (by Theorem \ref{t4.1}). Thus, for $f_K \in L^2$, we have that
\begin{eqnarray}|\langle f_K,g \rangle | &=& \left|\sum \lambda_i \langle f_K,m_i \rangle \right|\nonumber\\&
\leq &
C \sum |\lambda_i| \|f_K\|_{BMO_{L^*}} \leq C \|f\|_{BMO_{L^*}} \|g\|_{H^1_L},\label{eq8.c}\end{eqnarray}
where we have used the claim (\ref{eq8.b}).  Now, we also have that
$$\langle f_K,g \rangle \to \dint_{\RR^{n+1}_+} (t^2L^*)^M e^{-t^2L^*} f(x) \,
\overline{t^2Le^{-t^2L}g(x)}\, \frac{dx dt}{t},$$ by a dominated
convergence argument which uses Lemma \ref{lmain} (ii), Lemma
\ref{l8.3}, and the duality of $T^1$ and $T^\infty$ \cite{CMS}.
But by Lemma \ref{l8.4},  the last expression equals $\langle
f,g\rangle $; i.e., $$\langle f_K,g\rangle \to \langle
f,g\rangle,$$  so that (\ref{eq8.19}) follows from (\ref{eq8.c}).

To complete the proof of Theorem \ref{t8.2}, it remains only to establish the claims concerning $f_K$.  To see that $f_K \in L^2(\RR^n)$, it suffices by
Lemma \ref{l8.3} to observe that for all $F_t \in T^\infty (\RR^{n+1}_+)$, we have
\begin{eqnarray}\left\|\int_{1/K}^K t^2L^* e^{-t^2L^*} \left( \chi_{B_K} F_t \right)
\frac{dt}{t}\right\|_{L^2(\RR^n)}
&\leq & C\int_{1/K}^K \|F_t\|_{L^2 (B_K) } \frac{dt}{t}\nonumber\\
&\leq & C_K \left( \int_0^K \int_{B_K} |F_t(x)|^2 dx \frac{dt}{t}\right)^{1/2} \nonumber \\
&\leq & C_K \|F_t\|_{T^\infty(\RR^{n+1}_+)} |B_K|^{1/2}.
\label{eq8.d}\end{eqnarray} To prove the claim (\ref{eq8.b}),
again by Lemma \ref{l8.3} it suffices to prove the following
\begin{lemma}\label{l8.5}  Suppose that $F_t \in T^\infty(\RR^{n+1}_+),$ and set
$$f_K \equiv \int_{1/K}^K t^2L^* e^{-t^2L^*} \left( \chi_{B_K} F_t \right) \frac{dt}{t}.$$
Then $$ \sup_K \|f_K\|_{BMO_{L^*}(\RR^n)} \leq C \|F_t\|_{T^\infty(\RR^{n+1}_+)}.$$
\end{lemma}

\bp
We need to prove that for every cube $Q\subset\RR^n$

\begin{equation}\label{eq8.35}
\hskip -.3cm \left(\frac{1}{|Q|}\int_Q
\left|(I-(1+l(Q)^2L^*)^{-1})^M\int_{1/K}^K
tL^*e^{-t^2L^*}\left(\chi_{B_K}F_t\right)\,dt\right|^2\,dx\right)^{\frac
12}\leq C\|F_t\|_{T^\infty(\RR^{n+1}_+)},
\end{equation}
uniformly in $K$.  To this end, we
split the integral in $t$ in (\ref{eq8.35}) into two integrals
over $(1/K,l(Q)]$ and $(l(Q), K)$ (these are of course vacuous if $\ell(Q) < 1/K$, or if
$\ell(Q) > K$, respectively), and consider first the case
$t\leq l(Q)$. Let $h\in L^2(\RR^n)$ such that ${\rm
supp}\,h\subset Q$ and $\|h\|_{L^2(\RR^n)}=1$.   The left hand side
of (\ref{eq8.35}), restricted to $t\leq l(Q)$, is bounded by the supremum over all such $h$ of
the following:

\begin{eqnarray}\nonumber
&& \frac{1}{|Q|^{1/2}}\left|\int_Q
(I-(1+l(Q)^2L^*)^{-1})^M\int_{1/K}^{l(Q)}
tL^*e^{-t^2L^*}\left(\chi_{B_K}F_t\right)(x)\,dt\,\overline{h(x)}\,dx\right| \\[4pt]\nonumber
&& \qquad \leq
\frac{C}{|Q|^{1/2}}\left|\int_{1/K}^{l(Q)}\int_{\RR^n}
\left(\chi_{B_K}F_t(x)\right)\,\overline{t^2Le^{-t^2L}(I-(1+l(Q)^2L)^{-1})^M
h(x)}\, \frac{dxdt}{t}\right|
\\[4pt]\nonumber
&& \qquad \leq C\sum_{j=0}^\infty
\left(\frac{1}{|Q|}\int_0^{l(Q)}\int_{S_j(Q)}|F_t(x)|^2
\,\frac{dxdt}{t}\right)^{1/2}\\[4pt]\nonumber
&&\qquad \qquad \times \left(\int_0^{l(Q)}\int_{S_j(Q)}
\Bigl|t^2Le^{-t^2L}(I-(1+l(Q)^2L)^{-1})^M
h(x)\Bigr|^2\,\frac{dxdt}{t}\right)^{1/2}\\[4pt]
&& \qquad \leq C\sum_{j=0}^\infty 2^{jn/2}
\|F_t\|_{T^\infty(\RR^{n+1}_+)}\nonumber\\[4pt]
&&\qquad \qquad \times  \left(\int_0^{l(Q)}\int_{S_j(Q)}
\Bigl|t^2Le^{-t^2L}(I-(1+l(Q)^2L)^{-1})^M
h(x)\Bigr|^2\,\frac{dxdt}{t}\right)^{1/2} \label{eq8.36},
\end{eqnarray}

\noindent where we majorized the integral over $S_j(Q)\times
(0,l(Q))$ by the integral over $\widehat B$ for some ball $B$ with
size comparable to $(2^jl(Q))^n$ in the last inequality.

If $j=0,1$

\begin{eqnarray}\nonumber
&& \left(\int_0^{l(Q)}\int_{S_j(Q)}
\Bigl|t^2Le^{-t^2L}(I-(1+l(Q)^2L)^{-1})^M
h(x)\Bigr|^2\,\frac{dxdt}{t}\right)^{1/2}\\[4pt]
&&\qquad \leq C\sup_{0\leq k\leq M} \|g_h
(1+l(Q)^2L)^{-k}h\|_{L^2(\RR^n)}\leq C, \label{eq8.37}
\end{eqnarray}

\noindent since $g_h$ is bounded on $L^2(\RR^n)$ (see \cite{ADM})
and $(1+l(Q)^2L)^{-1}$ is uniformly bounded on $L^2(\RR^n)$ (see
Lemma~\ref{l2.4}).

Assume now that $j\geq 2$. Then

\begin{eqnarray}\nonumber
&& \left(\int_0^{l(Q)}\int_{S_j(Q)}
\Bigl|t^2Le^{-t^2L}(I-(1+l(Q)^2L)^{-1})^M
h(x)\Bigr|^2\,\frac{dxdt}{t}\right)^{1/2}\\[4pt]
&&\qquad \leq C\sup_{0\leq k\leq M}
\left(\int_0^{l(Q)}\int_{S_j(Q)}
\Bigl|t^2Le^{-t^2L}(1+l(Q)^2L)^{-k}
h(x)\Bigr|^2\,\frac{dxdt}{t}\right)^{1/2}. \label{eq8.38}
\end{eqnarray}

\noindent When $k=0$,

\begin{equation}
\left(\int_0^{l(Q)}\int_{S_j(Q)} \Bigl|t^2Le^{-t^2L}
h(x)\Bigr|^2\,\frac{dxdt}{t}\right)^{1/2}\leq
C\left(\int_0^{l(Q)}e^{-\frac{(2^jl(Q))^2}{ct^2}}
\,\frac{dt}{t}\right)^{1/2}\leq C2^{-jN}, \label{eq8.39}
\end{equation}

\noindent for every $N\in\NN$. Here we used Gaffney estimates and
the fact that ${\rm supp}\,h\subset Q$, $\|h\|_{L^2(\RR^n)}=1$.

When $1\leq k\leq M$, the quantity under the $\sup$ sign in
(\ref{eq8.38}) can be rewritten as

\begin{eqnarray}\nonumber
&&\left(\int_0^{l(Q)}\int_{S_j(Q)}\frac{t^4}{l(Q)^4}\,
\Bigl|e^{-t^2L}\bigl[(l(Q)^2L)(1+l(Q)^2L)^{-1}\bigr](1+l(Q)^2L)^{-k+1}
h(x)\Bigr|^2\,\frac{dxdt}{t}\right)^{1/2}\\[4pt]\nonumber
&&\quad\leq C \left(\int_0^{l(Q)}\int_{S_j(Q)}\frac{t^4}{l(Q)^4}\,
\Bigl|e^{-t^2L}(I-(1+l(Q)^2L)^{-1})(1+l(Q)^2L)^{-k+1}
h(x)\Bigr|^2\,\frac{dxdt}{t}\right)^{1/2}
\\[4pt]\nonumber
&&\quad\leq C\left(\int_0^{l(Q)}\frac{t^4}{l(Q)^4}\,
e^{-\frac{(2^jl(Q))^2}{ct^2}} \,\frac{dt}{t}\right)^{1/2}+C
\left(\int_0^{l(Q)}\frac{t^4}{l(Q)^4}\,
e^{-\frac{(2^jl(Q))^2}{cl(Q)^2}} \,\frac{dt}{t}\right)^{1/2},
\end{eqnarray}

\noindent where the first term above comes from the case $k=1$ and
we use Lemmas~\ref{l2.4} and \ref{l2.3}. The last sum in
(\ref{eq8.40}) is bounded by $C2^{-jN}$ for every $N\in\NN$, and
combining (\ref{eq8.36})--(\ref{eq8.40}) we deduce the desired
estimate for (\ref{eq8.36}) when $t\leq l(Q)$.

As for the case $t\in (l(Q),K)$,

\begin{eqnarray}\nonumber
&& \frac{1}{|Q|^{1/2}}\int_Q
(I-(1+l(Q)^2L^*)^{-1})^M\int_{l(Q)}^K
tL^*e^{-t^2L^*}\left(\chi_{B_K}F_t\right)(x)\,dt\,\overline{h(x)}\,dx \\[4pt]\nonumber
&& \qquad \leq C\sum_{k=0}^\infty\sum_{j=0}^\infty
\left(\frac{1}{|Q|}\int_{2^kl(Q)}^{2^{k+1}l(Q)}\int_{S_j(2^kQ)}|F_t(x)|^2
\,\frac{dxdt}{t}\right)^{1/2}\times\\[4pt]\nonumber
&&\qquad \qquad \times
\left(\int_{2^kl(Q)}^{2^{k+1}l(Q)}\int_{S_j(2^kQ)}
\Bigl|t^2Le^{-t^2L}(I-(1+l(Q)^2L)^{-1})^M
h(x)\Bigr|^2\,\frac{dxdt}{t}\right)^{1/2}\\[4pt]
&& \qquad \leq C\sum_{k=0}^\infty\sum_{j=0}^\infty 2^{(j+k)n/2}
\|F_t\|_{T^\infty(\RR^{n+1}_+)}\times\nonumber\\[4pt]\nonumber
&&\qquad \qquad \times
\left(\int_{2^kl(Q)}^{2^{k+1}l(Q)}\int_{S_j(2^kQ)}
\left(\frac{l(Q)}{t}\right)^{4M}
\Bigl|(t^2L)^{M+1}e^{-t^2L}(1+l(Q)^2L)^{-M}
h(x)\Bigr|^2\,\frac{dxdt}{t}\right)^{1/2}\nonumber\\[4pt]
&& \qquad \leq C\sum_{k=0}^\infty\sum_{j=0}^\infty
2^{jn/2}2^{k(n/2-2M)}
\|F_t\|_{T^\infty(\RR^{n+1}_+)}\times\nonumber\\[4pt]
&&\qquad \qquad \times
\left(\int_{2^kl(Q)}^{2^{k+1}l(Q)}\int_{S_j(2^kQ)}
\Bigl|(t^2L)^{M+1}e^{-t^2L}(1+l(Q)^2L)^{-M}
h(x)\Bigr|^2\,\frac{dxdt}{t}\right)^{1/2}. \label{eq8.40}
\end{eqnarray}

From this point the argument is essentially the same as the one
for small $t$. For $j=0,1$ the expression in the parentheses above
is bounded by $C\|h\|_{L^2(\RR^n)}^2\leq C$, and the sum in $k$
converges for $M>n/4$. For $j\geq 2$ we use Gaffney estimates to
bound the quantity in (\ref{eq8.40}) by

\begin{equation}
\sum_{k=0}^\infty\sum_{j=2}^\infty 2^{jn/2}2^{k(n/2-2M)}
\|F_t\|_{T^\infty(\RR^{n+1}_+)} \left(\int_{2^kl(Q)}^{2^{k+1}l(Q)}
e^{-\frac{(2^{j+k})^2}{ct^2}}\,\frac{dt}{t}\right)^{1/2}\leq C
\|F_t\|_{T^\infty(\RR^{n+1}_+)}. \label{eq8.41}
\end{equation}

Therefore, (\ref{eq8.35}) is valid for all $M>n/4$.
This concludes the proof of Lemma \ref{l8.5}, and thus also that of Theorem \ref{t8.2}. \ep

Next, we prove the converse:
\begin{theorem}\label{t8.3} Suppose $M>n/4$, $\eps>0$,
and that $l$ is a bounded linear functional  on $H^1_L(\RR^n)$. Then in fact,
$l\in BMO_{L^*}(\RR^n)$ and for all $g\in H^1_L(\RR^n)$
which can be represented as finite linear combinations of $(2,\eps,M)$-molecules,
we have

\begin{equation}\label{eq8.20}
l(g)=\langle l,g \rangle,
\end{equation}
where the latter pairing is in the sense of ${\bf M}^{2,\eps,M}_0(L)$ and its dual.
Moreover,

\begin{equation}\label{eq8.21}
\|l\|_{BMO_{L^*}(\RR^n)}\leq C\|l\|.
\end{equation}
\end{theorem}
We observe that the combination of Theorems~\ref{t8.2} and \ref{t8.3} gives
Theorem~\ref{t1.2}.

\bigskip
\bp By Theorem \ref{t4.1} and its proof, we have in particular
that for any $(2,\eps,M)$-molecule $m$, $\|m\|_{H^1_L} \leq C.$
Thus, $$l(m) \leq C\|l\|$$ for every $(2,\eps,M)$-molecule $m$. In
particular, $l$ defines a linear functional on ${\bf
M}^{2,\eps,M}_0(L)$ for every $\eps >0, M > n/4.$  Thus, $(I - (I
+ t^2 L^*)^{-1})^M l$ is well defined and belongs to $L^2_{loc}$
for every $t>0$ (recall (\ref{eq1.loc}) and the related
discussion).  Fix a cube $Q$, and let $\varphi \in L^2(Q)$, with
$\|\varphi\|_{L^2(Q)} \leq |Q|^{-1/2}.$  As we have observed above
(\ref{eq1.mol}),
$$ \tilde{m} \equiv (I - (I + \ell(Q)^2 L)^{-1})^M \varphi $$
is (up to a harmless multiplicative constant) a $(2,\eps,M)$-molecule for every $\eps > 0.$
Thus, \begin{eqnarray*} |\langle (I - (I + t^2 L^*)^{-1})^M l, \varphi \rangle | &\equiv &
|\langle  l, (I - (I + t^2 L)^{-1})^M\varphi \rangle |\\
&=& |\langle  l, \tilde{m} \rangle | \leq C \|l\|.
\end{eqnarray*}
Taking a supremum over all such $\varphi$ supported in $Q$, we obtain that
$$\frac{1}{|Q|} \int_Q |(I - (I + t^2 L^*)^{-1})^M l (x)|^2 dx \leq C \|l\|^2.$$
Since $Q$ was arbitrary, the conclusion of the theorem follows.

\ep

\begin{corollary}\label{c8.4} The operator
$(L^*)^{-1/2}\,{\rm div}=\left(\nabla L^{-1/2}\right)^*$ is
bounded from $L^\infty(\RR^n)$ to \break $BMO_{L^*}(\RR^n)$.
\end{corollary}

\bp The corollary follows from Theorem~\ref{t3.2} and
Theorem~\ref{t1.2} or can be proved directly. The argument is
standard, we leave the details to the interested reader. \ep

We conclude this section with the following consequence of Theorems \ref{t8.2} and \ref{t8.3},
and Corollary \ref{c4.3}.

\begin{corollary}\label{c8.5} For every $M > n/4$, the spaces $BMO_L(\RR^n)$
defined by the norms (\ref{eq1.18})
are equivalent.
\end{corollary}

\section{$BMO_L(\RR^n)$:\,connection with Carleson measures.}
\setcounter{equation}{0}

A Carleson measure is a positive measure $\mu$ on $\RR^{n+1}_+$ such
that

\begin{equation}\label{eq9.1}
\|\mu\|_{\cal C}:=\sup_{B}\frac{1}{|B|}\,\mu\bigl(\widehat
B\bigr)<\infty,
\end{equation}

\noindent where $B$ denotes a ball in $\RR^n$ and $\widehat B$ is
a tent over $B$ (see (\ref{eq8.23})). Recall the definition of the
operator $C$ in (\ref{eq8.22}) and observe that

\begin{equation}\label{eq9.2}
\|CF\|_{L^\infty(\RR^n)}^2=
\left\||F(y,t)|^2\,\frac{dydt}{t}\right\|_{\cal C}.
\end{equation}

\begin{theorem}\label{t9.1} Assume that $M\in\NN, M > n/4$. Then for every $f\in BMO_L(\RR^n)$

\begin{equation}\label{eq9.3}
\mu_f:=\left|(t^2L)^Me^{-t^2L}f(y)\right|^2\,\frac{dydt}{t}
\end{equation}

\noindent is a Carleson measure and

\begin{equation}\label{eq9.4}
\left\|\mu_f\right\|_{\cal C}\leq C\|f\|_{BMO_L(\RR^n)}^2.
\end{equation}

Conversely, if $f\in ({\bf M}^{2,M}_0(L^*))^*$ satisfies the controlled growth bound
(\ref{eq9.5}) (with $L$ in place of $L^*$) for some $\eps_1>0$, and if $\mu_f$ defined in (\ref{eq9.3})
is a Carleson measure, then $f\in BMO_L(\RR^n)$ and

\begin{equation}\label{eq9.6}
\|f\|_{BMO_L(\RR^n)}^2\leq C\left\|\mu_f\right\|_{\cal C}.
\end{equation}
\end{theorem}
\bp The direction $BMO_L$ implies (\ref{eq9.4}) is just a
restatement of Lemma \ref{l8.3}.

For the converse we use the duality of the tent spaces. More
precisely, for $f$ satisfying (\ref{eq9.5}) and every $g\in
H^1_{L^*}(\RR^n)$\, that can be represented as a finite linear
combination of $(2,\eps,M')$-molecules, $\eps>\eps_1/2$ and $M'>n/4$ large enough compared to $M$,
we have by Lemma \ref{l8.4} that

\begin{equation}
\langle f,g \rangle=C_M\dint_{\RR^{n+1}_+}
(t^2L)^Me^{-t^2L}f(x)\,\overline{t^2L^*e^{-t^2L^*}g(x)}\,\frac{dtdx}{t}.
\end{equation}\nonumber

Now according to Theorem~1 in \cite{CMS}

\begin{eqnarray}
&&\dint_{\RR^{n+1}_+}\left|
(t^2L)^Me^{-t^2L}f(x)\,\overline{t^2L^*e^{-t^2L^*}g(x)}\right|\,dx\nonumber\\[4pt]
&&\qquad\leq C\int_{\RR^n}
C((t^2L)^Me^{-t^2L}f)(x)\,S(t^2L^*e^{-t^2L^*}g)(x)\,dx\nonumber\\[4pt]
&&\qquad\leq C\,\|C((t^2L)^Me^{-t^2L}f)\|_{L^\infty(\RR^n)}
\,\|S(t^2L^*e^{-t^2L^*}g)\|_{L^1(\RR^n)}. \label{eq9.22}
\end{eqnarray}

\noindent Then using (\ref{eq9.2}) and Theorem~\ref{t4.1}, we have

\begin{equation}\label{eq9.23}
\left| \langle f,g \rangle \right|\leq
C\left\|\Bigl|(t^2L)^Me^{-t^2L}f(y)\Bigr|^2\,\frac{dydt}{t}\,\right\|_{\cal
C}^{1/2} \,\|g\|_{H^1_{L^*}(\RR^n)},
\end{equation}

\noindent for every $g\in H^1_{L^*}(\RR^n)$. By
Theorem~\ref{t8.3} this gives the desired conclusion
(\ref{eq9.6}). \ep

\noindent {\it Remark.}\, Fix some $p\in (p_L,2)$. Using the
finite linear combinations of $(p',\eps,M')$- molecules in the
proof of Lemma \ref{l8.4}, $\frac
1p+\frac{1}{p'}=1$, we can prove that the condition (\ref{eq9.5})
can be replaced by

\begin{equation}\label{eq9.24}
\int_{\RR^n}\frac{|(I-(1+L)^{-1})^Mf(x)|^p}{1+|x|^{n+\eps_1}}\,dx<\infty
\end{equation}

\noindent for some $\eps_1>0$. We will use this fact later in conjunction
with the fact that every $f\in BMO_L^p(\RR^n)$ satisfies
(\ref{eq9.24}).

\section{John-Nirenberg inequality}
\setcounter{equation}{0}

We start with the following auxiliary result, which is a
modification of Lemma~2.14 in \cite{acta}.

\begin{lemma}\label{l10.1} Suppose there exist numbers $0<\alpha<1$
and $0<N<\infty$ such that for some function $F\in
L^2_{loc}((0,\infty)\times \RR^n)$, some $a\in \RR$ and every cube
$Q\subset \RR^n$

\begin{equation}\label{eq10.1}
\left|\left\{x\in Q:\,\Bigl(\dint_{|x-y|<3at<3al(Q)}
|F(t,y)|^2\,\frac{dydt}{t^{n+1}}\Bigr)^{1/2}> N\right\}\right|\leq
\alpha|Q|.
\end{equation}

\noindent Then there exists $C>0$ such that

\begin{equation}\label{eq10.2}
\sup_{Q\subset\RR^n}\,\frac{1}{|Q|}\int_Q
\left(\dint_{|x-y|<at<al(Q)}
|F(t,y)|^2\,\frac{dydt}{t^{n+1}}\right)^{p/2}\,dx\leq C,
\end{equation}

\noindent for all $p\in(1,\infty)$.
\end{lemma}

\bp Denote the set on the left-hand side of (\ref{eq10.1}) by
$\Omega$, so that $|\Omega|\leq \alpha|Q|$, and let $\cup Q_j$ be
a Whitney decomposition of $\Omega$. Also,
$$
M(\delta):=\sup_{Q'\subseteq Q} \frac{1}{|Q'|}\int_{Q'}
\left(\dint_{\stackrel{|x-y|<a(t-\delta)}{\delta<t<l(Q')}}
|F(t,y)|^2\,\frac{dydt}{t^{n+1}}\right)^{p/2}\,dx,
$$

\noindent where the integral is set to be zero whenever
$l(Q')\leq\delta$. Then

\begin{eqnarray}\nonumber
&&\hskip -.7cm \int_Q
\left(\dint_{\stackrel{|x-y|<a(t-\delta)}{\delta<t<l(Q)}}
|F(t,y)|^2\,\frac{dydt}{t^{n+1}}\right)^{p/2}dx\\[4pt]\nonumber
&&\hskip -.7cm \quad \leq \int_{Q\setminus\Omega}
\left(\dint_{\stackrel{|x-y|<3at}{t<l(Q)}}
|F(t,y)|^2\,\frac{dydt}{t^{n+1}}\right)^{p/2}dx
+\sum_{j:\,l(Q_j)>\delta} \int_{Q_j}
\left(\dint_{\stackrel{|x-y|<a(t-\delta)}{\delta<t<l(Q_j)}}
|F(t,y)|^2\,\frac{dydt}{t^{n+1}}\right)^{p/2}dx\\[4pt]\nonumber
&&\hskip -.7cm \quad \qquad+ \sum_j \int_{Q_j}
\left(\dint_{\stackrel{|x-y|<a(t-\delta)}{\max\{l(Q_j),\delta\}<t<l(Q)}}
|F(t,y)|^2\,\frac{dydt}{t^{n+1}}\right)^{p/2}dx\\[4pt]\label{eq10.3}
&&\hskip -.7cm \quad \leq N^p|Q|+M(\delta) \alpha |Q|+ \sum_j
\int_{Q_j}
\left(\dint_{\stackrel{|x-y|<a(t-\delta)}{\max\{l(Q_j),\delta\}<t<l(Q)}}
|F(t,y)|^2\,\frac{dydt}{t^{n+1}}\right)^{p/2}dx.
\end{eqnarray}

\noindent By the properties of the Whitney decomposition ${\rm
dist}\,(x',Q_j)\approx l(Q_j)$ for some $x'\in Q\setminus\Omega$.
Therefore, without loss of generality we can assume that for every
$x\in Q_j$ there exists $x'\in Q\setminus\Omega$ such that

\begin{equation}\label{eq10.4}
\left\{y:\,|x-y|<a(t-\delta),\,\max\{l(Q_j),\delta\}<t<l(Q)\right\}\,\subset\,
\left\{y:\,|x'-y|<3at,\,t<l(Q)\right\}.
\end{equation}

\noindent Then the last term in (\ref{eq10.3}) can be bounded by

\begin{equation}\label{eq10.5}
|Q| \sup_{x'\in
Q\setminus\Omega}\left(\dint_{\stackrel{|x'-y|<3at}{t<l(Q)}}
|F(t,y)|^2\,\frac{dydt}{t^{n+1}}\right)^{p/2}\leq N^p|Q|.
\end{equation}

\noindent Repeating the procedure in
(\ref{eq10.3})--(\ref{eq10.5}) for every cube $Q'\subset Q$, we
arrive at

\begin{equation}\label{eq10.6}
(1-\alpha) M(\delta)\leq C N^p,
\end{equation}

\noindent and the argument can be finished taking limit as
$\delta\to 0$. \ep

\vskip 0.08in

\noindent {\it Proof of Theorem~\ref{t1.3}}.\,\,\, {\bf Step I.}\,
It is an easy consequence of H\"older's inequality that

\begin{equation}\label{eq10.7}
\|f\|_{BMO_L^p(\RR^n)}\leq C \|f\|_{BMO_L(\RR^n)}\leq C
\|f\|_{BMO_L^q(\RR^n)} \quad \mbox{for}\quad p_L<p<2<q<\widetilde{p}_L,
\end{equation}

\noindent so we will concentrate on the reverse estimates.

\vskip 0.08in \noindent {\bf Step II.}\, In this part we will show
that $f\in BMO_L^{p_0}(\RR^n)$ implies $f\in BMO_L(\RR^n)$ with
$p_L<p_0<2$.  Let us first prove that whenever $f\in
BMO_L^{p_0}$ the inequality (\ref{eq10.2}) holds with $p=p_0$,
$a=1$. We split $f$ as in (\ref{eq8.2}). Then the contribution of
the first part is handled as follows. Making the dyadic annular
decomposition,

\begin{eqnarray}
&&\hskip -0.7cm \left(\frac{1}{|Q|}\int_Q \left(\dint_{|x-y|<t<l(Q)}
|(t^2L)^M
e^{-t^2L}(I-(1+l(Q)^2L)^{-1})^Mf(y)|^2\,\frac{dydt}{t^{n+1}}\right)^{\frac{p_0}{2}}
\,dx\right)^{\frac{1}{p_0}}
\nonumber\\[4pt]
&&\hskip -0.7cm\leq C \sum_{j=0}^\infty\Bigl(\frac{1}{|Q|}\int_Q
\Bigl(\dint_{|x-y|<t<l(Q)} |(t^2L)^M e^{-t^2L}\times
\nonumber\\[4pt]
&&\hskip -0.7cm \qquad\qquad\times
\bigl[\chi_{S_j(Q)}(I-(1+l(Q)^2L)^{-1})^Mf\bigr](y)|^2\,
\frac{dydt}{t^{n+1}}\Bigr)^{\frac{p_0}{2}}dx\Bigr)^{\frac{1}{p_0}} \nonumber\\[4pt]
&&\hskip -0.7cm\leq C\frac{1}{|Q|^{1/p}}\left(\int_{4Q}
|(I-(1+l(Q)^2L)^{-1})^Mf(x)|^{p_0}\,dx\right)^{\frac{1}{p_0}} \nonumber\\[4pt]
&&\hskip -0.7cm\quad +\,
C\sum_{j=3}^\infty\frac{1}{|Q|^{1/2}}\Bigl(\dint_{2Q\times(0,l(Q))}
|(t^2L)^M e^{-t^2L}\times
\nonumber\\[4pt]
&&\hskip -0.7cm
\qquad\qquad\times\bigl[\chi_{S_j(Q)}(I-(1+l(Q)^2L)^{-1})^Mf\bigr](x)|^2\,
\frac{dxdt}{t}\Bigr)^{\frac{1}{2}}.\label{eq10.8}
\end{eqnarray}

\noindent Here we used $L^{p_0}(\RR^n)$ boundedness of the conical
square function (see Lemma~\ref{l2.6}) for the first term above.
To get the second one we applied H\"older's inequality to pass
from $L^{p_0}$ to $L^2$ norm and then Lemma~\ref{l2.1}. Using
$L^{p_0}-L^2$ off-diagonal estimates (Lemma~\ref{l2.5}) the second
term can be further bounded by

\begin{equation}
\hskip -.06 cm
C\sum_{j=3}^\infty\frac{1}{|Q|^{1/{2}}}\left(\int_0^{l(Q)}\!
e^{-\frac{(2^jl(Q))^{2}}{ct^2}}
\Bigl(\int_{S_j(Q)}\!|(I-(1+l(Q)^2L)^{-1})^Mf(x)|^{p_0}\,
dx\Bigr)^{\frac{2}{p_0}}\!\!
\frac{dt}{t^{1-\frac{n}{2}+\frac{n}{p_0}}}\right)^{\frac
12}\!.\label{eq10.9}
\end{equation}

\noindent Then covering $S_j(Q)$ by approximately $2^{jn}$ cubes
of the sidelength $l(Q)$ and integrating in $t$ we control
(\ref{eq10.9}) and hence (\ref{eq10.8}) by
$\|f\|_{BMO_L^{p_0}(\RR^n)}$.

The contribution of

\begin{equation}\left[I-(I-(1+l(Q)^2L)^{-1})^M\right]f= \sum_{k=1}^M C_{k,M}
\,(l(Q)^2L)^{-k}(I-(1+l(Q)^2L)^{-1})^Mf \label{eq10.10}
\end{equation}

\noindent (cf. (\ref{eq8.5})) can be estimated in the same way as
(\ref{eq10.8})--(\ref{eq10.9}), first combining $L^{-k}$, $1\leq
k\leq M$, with $L^M$.

We just proved that $F$ given by $F(t,y)=t^2Le^{-t^2L}f(y)$
satisfies (\ref{eq10.2}) for some $p=p_0$ and $a=1$. Then by
Chebyshev's inequality it satisfies (\ref{eq10.1}) with $a=1/3$.
Hence (\ref{eq10.2}) holds for $p=2$, $a=1/3$ and $F$ as above by
Lemma~\ref{l10.1}. The latter fact implies that $f\in
BMO_L(\RR^n)$ using Theorem~\ref{t9.1} and the Remark after
Theorem~\ref{t9.1}.

\vskip 0.08in \noindent {\bf Step III.}\, Let us consider the
estimate

\begin{equation}\label{eq10.11}
\|f\|_{BMO_L^q(\RR^n)}\leq C \|f\|_{BMO_L(\RR^n)} \quad
\mbox{for}\quad 2<q<\widetilde{p}_L.
\end{equation}

Fix a cube $Q$, and let $\varphi \in L^2(Q)\supseteq L^p(Q),$
where $p=q'$, i.e., $$(\widetilde{p}_L)' \equiv \widetilde{p}_L/(\widetilde{p}_L -1) < p
\equiv q/(q-1) <2,$$ and suppose that $$\|\varphi\|_p \leq |Q|^{1/p - 1}.$$
We claim that for some harmless constant $C_0$,
$$m \equiv \frac{1}{C_0}\left(I - (I+\ell(Q)^2 L^*)^{-1}\right)^M \varphi$$
is a $(p,\eps,M)$-molecule, for the operator $L^*$,  adapted to $Q$,
for every $\eps >0$. Indeed, by a simple duality argument, we have
that $$(\widetilde{p}_L)' = p_{L^\ast},\,\,\, p_L' =
\widetilde{p}_{L^*}.$$ Thus, for $(\widetilde{p}_L)'=p_{L^\ast} < p
< 2$, the resolvent kernel $(I + t^2 L^*)^{-1}$ satisfies the
$L^p-L^p$ off-diagonal estimates by Lemma~\ref{l2.5}. Taking $E=Q$,
and $F=S_i(Q)$, the reader may then readily verify that $m$ is a
molecule as claimed. We omit the details.

Now, suppose that $f \in BMO_L$.  Then by Theorem \ref{t8.2}, $f
\in (H^1_{L^*})^*$.  Thus, since $\|m\|_{H^1_{L^*}} \leq C$, we
have that $$\frac{1}{C_0}\left|\left\langle \left(I - (I+\ell(Q)^2
L)^{-1}\right)^M f,\varphi \right\rangle \right| \equiv |\langle
f,m \rangle | \leq C \|f\|_{BMO_L(\RR^n)}.$$ Thus, taking a
supremum over all $\varphi$ as above, we obtain (\ref{eq10.11}).

\ep

\vskip 0.20in \noindent --------------------------------------
\vskip 0.40in

\noindent {\tt Steve Hofmann}

\noindent Department of Mathematics

\noindent University of Missouri at Columbia

\noindent Columbia, MO 65211, USA

\noindent {\tt e-mail}: {\it hofmann\@@math.missouri.edu}

\vskip 0.15in

\noindent {\tt Svitlana Mayboroda}

\noindent Department of Mathematics

\noindent The Ohio State University

\noindent 231 W 18th Avenue

\noindent Columbus, OH 43210 USA

\noindent {\tt e-mail}: {\it svitlana\@@math.ohio-state.edu}


\begin{thebibliography}{999}

\bibitem{ADM} D.\,Albrecht, X.\,Duong, A.\,McIntosh, {\it Operator theory
and harmonic analysis}, Instructional Workshop on Analysis and
Geometry, Part III (Canberra, 1995), 77--136, Proc. Centre Math.
Appl. Austral. Nat. Univ., 34, Austral. Nat. Univ., Canberra, 1996.

\bibitem{PAPersonal} P.\,Auscher, {\it Personal communication}.

\bibitem{A} P.\,Auscher, {\it On $L^p$-estimates for square roots of second order elliptic
operators on $\RR^n$}, Publ. Mat. 48 (2004), 159-186.

\bibitem{AuscherSurvey} P.\,Auscher, {\it On necessary and sufficient
conditions for $L^p$-estimates of Riesz transforms associated to
elliptic operators on $\RR^n$ and related estimates}, preprint.

\bibitem{ACT} P.\,Auscher, T.\,Coulhon, Ph.\,Tchamitchian, {\it Absence de principe
du maximum pour certaines \'{e}quations paraboliques complexes}, Coll.
Math. 171 (1996), 87–95.

\bibitem{KatoMain} P.\,Auscher, S.\,Hofmann, M.\,Lacey, A.\,McIntosh,
Ph.\,Tchamitchian, {\it The solution of the Kato square root
problem for second order elliptic operators on ${\Bbb R}\sp n$},
Ann. of Math. (2) 156 (2002), no. 2, 633--654.

\bibitem{acta} P.\,Auscher, S.\,Hofmann, J.L.\,Lewis, Ph.\,Tchamitchian,
{\it Extrapolation of Carleson measures and the analyticity of
Kato's square-root operators}, Acta Math. 187 (2001), no. 2,
161--190.

\bibitem{AR} P.\,Auscher, E.\,Russ, {\it Hardy spaces and divergence
operators on strongly Lipschitz domains of $\Bbb R\sp n$}, J. Funct.
Anal. 201 (2003), no. 1, 148--184.

\bibitem{AT1} P.\,Auscher, Ph.\,Tchamitchian, {\it Calcul fonctionnel pr{\'e}cis{\'e} pour des
op{\'e}rateurs elliptiques complexes en dimension un (et
applications {\`a} certaines {\'e}quations elliptiques complexes en
dimension deux)}, Ann. Inst. Fourier 45 (1995), 721--778.

\bibitem{AT} P.\,Auscher,\, Ph. Tchamitchian, {\it Square root problem for divergence
operators and related topics}, Ast\'erisque, no. 249 (1998).

\bibitem{BK} S.\,Blunck, P.\,Kunstmann, {\it Weak-type $(p, p)$ estimates for Riesz
transforms}, Math. Z. 247 (2004), no. 1, 137--148.


\bibitem{BK2} S.\,Blunck, P.\,Kunstmann, {\it Calderón-Zygmund theory for non-integral operators and the $H\sp \infty$ functional calculus},  Rev. Mat. Iberoamericana  19  (2003),  no. 3, 919--942.


\bibitem{Co} R.R.\,Coifman, {\it A real variable characterization of $H\sp{p}$},
 Studia Math. 51 (1974), 269--274.

\bibitem{CMS} R.R.\,Coifman, Y.\,Meyer, E.M.\,Stein, {\it Some
new function spaces and their applications to harmonic analysis}, J.
Funct. Anal. 62 (1985), no. 2, 304--335.

\bibitem{CW} R.R.\,Coifman, G.\,Weiss, {\it Extensions
of Hardy spaces and their use in analysis}, Bull. Amer. Math. Soc.
83 (1977), no. 4, 569--645.

\bibitem{Davies} E.B.\,Davies, {\it Limits on $L^p$ regularity of selfadjoint elliptic operators},
J. Differential Equations 135, no. 1 (1997), 83--102.

\bibitem{DL} X.T.\,Duong and L.X. Yan, {\it Duality of Hardy and
BMO spaces associated with operators with heat kernel bounds}, J.
Amer. Math. Soc. 18 (2005), 943-973.

\bibitem{DL2} X.T.\,Duong and L.X. Yan, {\it New function spaces of BMO type,
the John-Nirenberg inequality, interpolation, and applications},
Comm. Pure Appl. Math. 58 (2005), no. 10, 1375-1420.

\bibitem{FeSt} C.\,Fefferman, E.M.\,Stein, {\it $H\sp{p}$ spaces
of several variables}, Acta Math. 129 (1972), no. 3-4, 137--193.

\bibitem{GaCu} J.\,Garcia-Cuerva, J.L.\,Rubio de Francia, {\it
Weighted norm inequalities and related topics}, North-Holland
Mathematics Studies, 116. North-Holland Publishing Co., Amsterdam,
1985.

\bibitem{HoMa} S.\,Hofmann, J.M.\,Martell, {\it $L\sp p$ bounds
for Riesz transforms and square roots associated to second order
elliptic operators}, Publ. Mat. 47 (2003), no. 2, 497--515.

\bibitem{JN} F.\,John, L.\,Nirenberg, {\it
On functions of bounded mean oscillation}, Comm. Pure Appl. Math.
14 (1961), 415--426.

\bibitem{KePi} C.E.\,Kenig, J.\,Pipher, {\it The Neumann problem for
elliptic equations with nonsmooth coefficients}, Invent. Math. 113
(1993), no. 3, 447--509.

\bibitem{La} R.H.\,Latter, {\it
A characterization of $H^{p}(\RR^n)$ in terms of atoms}, Studia
Math. 62 (1978), no. 1, 93--101.

\bibitem{MNP} V. G.\,Maz'ya, S. A.\,Nazarov and B. A.\,Plamenevski\u\i,
{\it Absence of a De Giorgi-type theorem for strongly elliptic
equations with complex coefficients}, Boundary value problems of
mathematical physics and related questions in the theory of
functions, 14. Zap. Nauchn. Sem. Leningrad. Otdel. Mat. Inst.
Steklov. (LOMI) 115 (1982), 156--168, 309.

\bibitem{St} E.M.\,Stein, {\it Singular integrals and differentiability
properties of functions}, Princeton Mathematical Series, No. 30
Princeton University Press, Princeton, N.J. 1970.

\bibitem{St2} E.M.\,Stein, {\it Harmonic analysis: real-variable methods,
orthogonality, and oscillatory integrals}, Princeton Mathematical
Series, 43. Monographs in Harmonic Analysis, III. Princeton
University Press, Princeton, NJ, 1993.

\bibitem{StPNAS} E.M.\,Stein, {\it Maximal functions. II. Homogeneous curves,}  Proc.
Nat. Acad. Sci. U.S.A.  73  (1976), no. 7, 2176--2177.

\bibitem{StWe} E.M.\,Stein, G.\,Weiss, {\it On the theory of harmonic
functions of several variables. I. The theory of $H\sp{p}$-spaces},
Acta Math. 103 (1960), 25--62.

\bibitem{TW} M. Taibleson, and G. Weiss, {\it The molecular characterization of certain Hardy spaces.  Representation theorems for Hardy spaces},  pp. 67--149, Ast\'{e}risque, 77, Soc. Math. France, Paris, 1980.

\bibitem{Wilson} J. M. Wilson, {\it On the atomic decomposition for Hardy
spaces}, Pacific J. Math. 116 (1985), no. 1, 201--207.

\end{thebibliography}
\end{document}